\documentclass[oneside, 12pt]{amsart}
\usepackage{amscd}
\usepackage{amssymb}
\usepackage{amsfonts}
\usepackage{remark}
\usepackage{enumerate}
\usepackage[all]{xy}

\newtheorem{thm}{Theorem}[section]
\newtheorem{theorem}[thm]{Theorem}

\newtheorem{lemma}[thm]{Lemma}
\newtheorem{corollary}[thm]{Corollary}
\newtheorem{proposition}[thm]{Proposition}
\newremark{definition}[thm]{Definition}
\newremark{ques}[thm]{Question}
\newremark{remark}[thm]{Remark}
\newremark{example}[thm]{Example}
\newremark{emp}[thm]{\kern -4pt}

\newcommand{\nice}{$\operatorname{WC}(0)$}

\newcommand\Supp{\operatorname{Supp}}
\newcommand{\Proj}{\operatorname{Proj}}
\newcommand{\Spec}{\operatorname{Spec}}
\newcommand{\Pic}{\mbox{\rm Pic}\kern 1pt}
\newcommand{\Br}{\mbox{\rm Br}\kern 1pt}
\newcommand{\charr}{\operatorname{char}}
\newcommand{\length}{\ell}
\newcommand{\cok}{{\mathcal O}_K}
\newcommand{\cO}{{\mathcal O}}

\newcommand{\X}{{\mathcal X}}
\newcommand{\Y}{{\mathcal Y}}

\newcommand{\m}{{\mathfrak m}}
\newcommand{\n}{{\mathfrak n}}
\newcommand{\Ht}{\operatorname{ht}}
\newcommand{\p}{{\mathfrak p}}
\newcommand{\q}{{\mathfrak q}}
\newcommand{\Ass}{\operatorname{Ass}}
\newcommand{\reg}{\mathrm{reg}}
\newcommand{\red}{\mathrm{red}}
\newcommand{\CH}{\mathcal A}
\newcommand{\CHU}{{\mathcal A}(X,U)}
\newcommand{\SCH}{{\mathcal S}}
\newcommand{\Ann}{\operatorname{Ann}}
\newcommand{\dv}{\operatorname{div}}
\newcommand{\Coker}{\operatorname{Coker}\kern 1pt}
\newcommand{\Ker}{\operatorname{Ker}}
\newcommand{\ord}{\operatorname{ord}}
\newcommand{\Gr}{\operatorname{Gr}}

\newcommand{\ol}{\overline} 
\newcommand{\Frac}{\operatorname{Frac}}
\newcommand{\FA}{\operatorname{FA} }
\newcommand{\deltadelta}{\overline{\delta}}
\newcommand{\codim}{{\operatorname{codim}}}

\newcommand{\sep}{\mathrm{sep}}

\begin{document}

\date{\today}
\thanks{D.L. was supported by NSF Grant 0902161.}
\begin{abstract}
Let $K$ be the field of fractions of a Henselian discrete valuation
ring ${\cO_K}$. Let $X_K/K$ be a smooth proper geometrically connected 
scheme admitting a regular model $X/{\cO_K}$. We show that the index 
$\delta(X_K/K)$ of  $X_K/K$ can be explicitly computed using data pertaining only
to the special fiber $X_k/k$ of the model $X$. 

We give two proofs of this theorem,
using two moving lemmas. One moving lemma pertains to horizontal $1$-cycles on a regular projective scheme
$X $ over the spectrum of a semi-local Dedekind domain, and the second moving lemma can be applied to $0$-cycles on an $\FA$-scheme $X$ which need not be regular.

The study of the local algebra needed to prove these moving lemmas led us to introduce
an invariant $\gamma(A)$ of a singular local ring $(A, \m)$: the greatest common divisor of all the 
Hilbert-Samuel multiplicities $e(Q,A)$, over all $\m$-primary ideals $Q$ in $\m$. 
We relate this invariant $\gamma(A)$ to the index of the exceptional divisor
in a resolution of the singularity of $\Spec A$, and we give a new way of computing the 
index of a smooth subvariety $X/K$ of ${\mathbb P}^n_K$ over any field $K$,  using the invariant $\gamma$ of the 
local ring at the vertex of a cone over $X$. 

\vspace*{.3cm}
\noindent 
\begin{tiny}KEYWORDS. \end{tiny} Index of a variety, separable index, moving lemma, $1$-cycles, $0$-cycles, rationally equivalent,
Hilbert-Samuel multiplicity,
resolution of singularities, cone.

\vspace*{.3cm}
\noindent 
\begin{tiny}MATHEMATICS SUBJECT CLASSIFICATION: 13H15, 14C25, 14D06, 14D10,  14G05, 14G20  \end{tiny}
\end{abstract}

\title
{
The index of an algebraic variety}  
\author{Ofer Gabber}
\author{Qing Liu}
\author{Dino Lorenzini}
\address{IH\'ES, 35 route de Chartres, 
91440 Bures-sur-Yvette, France}
\email{gabber@ihes.fr}
\address{Universit\'e de Bordeaux 1, Institut de Math\'ematiques de 
Bordeaux, 33405 Talence, France} 
\email{Qing.Liu@math.u-bordeaux1.fr}
\address{Department of mathematics, University of Georgia, 
Athens, GA 30602, USA} 
\email{lorenzin@uga.edu}

\maketitle

Let $W$ be a non-empty scheme of finite type over a field $F$. 
Let $\mathcal D(W/F)$ denote the set of all degrees of closed points of $W$.
The {\it index} $\delta(W/F)$ of $W/F$ is the greatest common
divisor of the elements of $\mathcal D(W/F)$.
  The index 
is also the smallest positive integer occurring as the degree of a 
$0$-cycle on $W$. When $W$ is integral, let $W^{\reg}$ denote 
the regular locus of $W$, open in $W$. We note in 
\ref{pro.deltareg} that $\delta(W^{\reg}/F)$ 
is a birational invariant of $W/F$.

Let now $K$ be the field of fractions of a discrete valuation
ring ${\cO_K}$ with residue field    $k$. Let $S:= \Spec {\cO_K}$.
Let $X \to S$ be a proper flat morphism, with $X$ regular and
irreducible. Let $X_K/K$ be the generic fiber of $X/S$. 
Write the special fiber $X_k$, viewed as a divisor on $X$, as $\sum_{i=1}^n r_i \Gamma_i$,
where for each $i=1,\dots, n$, $\Gamma_i$ is irreducible, of multiplicity 
$r_i$ in $X_k$. Using the intersection of Cartier divisors with 
$1$-cycles on the regular scheme $X $, 
we easily find that $\gcd_i\{r_i\delta(\Gamma_i/k)\}$ divides 
$\delta(X_K/K)$ (see \ref{defX/S}).  
Our theorem below strengthens this divisibility, and shows that when ${\cO_K}$ is Henselian,
the index of the generic fiber can be computed using only data pertaining to the special fiber.

\medskip
\noindent {\bf Theorem \ref{dx-ds}} {\em  Keep the above assumptions on 
$X/S$. 
\begin{enumerate}[\rm (a)]
\item 
Then $\gcd_i\{ r_i \delta(\Gamma^{\reg}_{i}/k)\}$  divides  $\delta(X_K/K)$.
\item  When ${\cO_K}$ is Henselian, then 
$ \delta(X_K/K) = \gcd_i\{ r_i\delta(\Gamma^{\reg}_{i}/k)\}$. 
\end{enumerate}
}
\medskip
Theorem \ref{dx-ds} answers positively a question of Clark 
(\cite{Cla}, Conj. 16). This theorem is known already when 
$k$ is a finite field (\cite{Bosch-Liu}, 1.6, see also \cite{C-S}, 3.1),
or when $k$ is algebraically closed (same proof 
as in \cite{Bosch-Liu}),
or when $X_K/K$ is a curve with 
semi-stable reduction (\cite{Cla}, Thm. 9). 

We give two proofs of Theorem \ref{dx-ds}, using two different 
moving lemmas which may be of 
independent interest. The first proof uses the 
Moving Lemma \ref{mv-slocal} stated below. A slightly strengthened version is proved in the text.
The definition and main properties of the notion of {\it rational equivalence of cycles} are recalled 
in section \ref{section2}.

\medskip
\noindent
{\bf Theorem \ref{mv-slocal}} {\em Let $R$ be a semi-local Dedekind domain, and let $S:=\Spec(R)$.
Let $X/S$ be flat and quasi-projective, with $X$ regular. 
Let $C$ be a $1$-cycle on $X$, closure in $X$ of a closed point of the generic fiber of $X$.
Let $F$ be a closed subset of $X$ such that for all $s \in S$, $F \cap X_s$ has codimension at least $1$ in $X_s$. 
Then $C$ is rationally equivalent to 
a cycle $C'$ on $X$ whose support does not meet $F$.}
\medskip

Our second proof of Theorem \ref{dx-ds} uses the Moving Lemma \ref{mvsimple} below, which allows some moving of a multiple 
of a cycle 
on a scheme $X$ which need not be regular. Recall (\ref{FA-g}) that 
$X$ {\it is an $\FA$-scheme} 
if every finite subset of $X$ is contained in  an 
affine open subset of $X$.

\medskip
\noindent
{\bf Theorem \ref{mvsimple}} {\em 
Let $X$ be a noetherian $\FA$-scheme.
Let $F$ be a closed subset of $X$ of positive codimension in $X$.
Let $x_0\in X$.  Let $Q$ be a $\m_{X,x_0}$-primary ideal of 
$\cO_{X,x_0}$, with Hilbert-Samuel multiplicity $e(Q, \cO_{X,x_0})$. Then 
the cycle $e(Q)[\overline{\{ x_0 \}}]$ is rationally equivalent in $X$ 
to a cycle $Z$ such that 
no irreducible cycle occurring in 
$Z$ is contained in $F$.}
\medskip

Theorem \ref{mvsimple} is a consequence of a local analysis of the noetherian local
ring $\cO_{X,x_0}$ found in section \ref{section-local}, and in particular in Theorem \ref{local-mv}.
Our investigation of the local algebra needed to prove 
Theorem \ref{mvsimple} led us to introduce the following local invariant in \ref{gammaA}.
Let $(A, \m)$ be any noetherian local ring. 
Let $\mathcal E(A)$ denote the set 
of all Hilbert-Samuel multiplicities $e(Q,A)$, for all $\m$-primary ideals $Q$ of $A$.
Define $\gamma(A)$ {\it to be the greatest common divisor of the elements of} $\mathcal E(A)$.
Theorem \ref{mvsimple} and the definition of  $\gamma(A)$ show that:

\medskip
\noindent
{\bf Corollary  
\ref{mv-degree}.} {\em   
Let $X/k$ be a reduced scheme of finite type over a field $k$
and let $x_0 \in X$ be a closed point. 
Then 
$\delta(X^{\reg}/k)$ divides $ \gamma(\cO_{X,x_0})\deg_k(x_0)$.
}

\medskip 
This statement is slightly strengthened when $\cO_{X,x_0}$ is not equidimensional in \ref{pro.deltareg2}.
Recall that the Hilbert-Samuel multiplicity $e(\m,A)$ is the smallest element in the set $\mathcal E(A)$,
and that it
is a measure of the singularity of the ring $A$: 
if the completion of $A$ is a domain (or, more generally, is unmixed \cite{HIO}, 6.8 and 6.9), then $A$ is regular if and 
only if $e(\m,A)=1$. 
The invariant $\gamma(A)$ is also related to the singularity of the ring $A$:
our next theorem shows that $\gamma(A)$ is equal to the index of the exceptional divisor of a desingularization
of $\Spec A$.

\medskip
\noindent
{\bf Theorem \ref{thm.na} and \ref{cor.compute-n}.} {\em 
Let $A$ be an excellent noetherian equidimensional local ring of positive dimension.
Let $X := \Spec A$, with closed point $x_0$.
Let $f : Y\to X$ be a proper birational morphism such that $Y$ is regular.
Let $E:=f^{-1}(x_0)$. 
Then 
$\gamma(A)=\delta(E/k(x_0))$.
 }
\medskip

Note that in the above theorem, the set of degrees $\mathcal D(E)$ and the set
of Hilbert-Samuel multiplicities $\mathcal E(A)$ need not be equal, 
and neither needs to contain the greatest common divisors of its elements.
The proof that $\gamma(A)=\delta(E/k(x_0))$ involves a third set $\mathcal N$ of integers
attached to $A$, which is an ideal in $\mathbb Z$, so that the greatest common divisor $n(A)$
of the elements of $\mathcal N$ belongs to $\mathcal N$ (\ref{def-na}). The proof shows 
that  
 $\gamma(A)=n(A)$ and  $n(A) =\delta(E/k(x_0))$. The properties of the invariant $n(A)$
 are further studied in section \ref{perspective}.

As an application of Theorem \ref{cor.compute-n}, we obtain a  new description of the index of a 
projective variety.

\medskip
\noindent
{\bf Theorem \ref{newdescriptionIndex}.} {\em 
Let $K$ be any field. Let $V/K$ be a regular 
closed integral subscheme of ${\mathbb P}^n_K$.
Denote by $W$ a 
cone over $V$ in ${\mathbb P}_K^{n+1}$. Let $w_0 \in W$
denote the vertex of the cone. 
Then $\delta(V/K) = \gamma(\cO_{W,w_0})$.
}

\medskip
The variant of Hilbert's Tenth Problem, which asks whether there exists an algorithm which decides
given a geometrically irreducible variety $V/{\mathbb Q}$,
whether $V/\mathbb Q$ has a ${\mathbb Q}$-rational point, is an open question to this date (\cite{Poo}, p.\ 348). 
So is the possibly weaker 
question of the existence of an algorithm which decides, given $V/{\mathbb Q}$,
whether $\delta(V/{\mathbb Q})=1$ (\ref{algorithm}). In view of Theorem \ref{newdescriptionIndex},
we may also ask whether there exists an algorithm which decides, given a ${\mathbb Q}$-rational
 point $w_0$ on a scheme of finite type $W/{\mathbb Q}$, whether $\gamma(\cO_{W,w_o})=1$.
 
 In the last section of this article, we settle a question of Lang and Tate
 \cite{LT}, page 670, when the ground field $K$ is imperfect, and prove: 

\medskip
\noindent
{\bf Theorem  \ref{sep.index}.} {\it  Let $X$ be a regular and generically smooth non-empty scheme of finite type over a
field $K$. Then the index $\delta(X/K)$ is equal to the separable index $\delta_{sep}(X/K)$.
}

 \medskip
 We thank the referee for a  meticulous reading of the article and for many useful comments.

\tableofcontents

\begin{section}{Rational equivalence} \label{section2}

We review below the basic notation needed to state our moving lemmas.
Let $X$ be a noetherian scheme. Let ${\mathcal Z}(X)$ denote the free 
abelian group on the set of closed integral subschemes of $X$.
An element of ${\mathcal Z}(X)$ is called a cycle, and if $Y$ is an integral closed subscheme of $X$,
we denote by $[Y]$ the associated element in ${\mathcal Z}(X)$.

Let ${\mathcal K}_X$ denote the sheaf of 
meromorphic functions 
on a noetherian scheme $X$ (see \cite{K}, top of page 204 or 
\cite{Liubook}, Definition 7.1.13). Let $f\in \mathcal K_X^*(X)$. 
Its associated  
principal Cartier divisor is denoted by $\dv(f)$ and defines a cycle on $X$:
$$[\dv(f)]=\sum_{x} \ord_x(f_x)[\overline{\{x\}}]$$ 
where $x$ ranges through the points of codimension $1$ in $X$ (i.e., the points $x$ such that the closure $\overline{\{x\}}$
has codimension $1$ in $X$; this latter condition is equivalent to the condition $\dim \cO_{X,x} =1$).
The function $\ord_x: \mathcal K_{X,x}^* \to \mathbb Z$ is defined, for a regular element of $g \in \cO_{X,x}$,
to be the length of the $\cO_{X,x}$-module $\cO_{X,x}/(g)$. 
 
A cycle $Z$ is \emph{rationally equivalent to $0$} (\cite{K2}, \S 2),
or \emph{rationally trivial}, 
if there are finitely many integral closed subschemes $Y_i$ and
principal Cartier divisors $\dv(f_i)$ on $Y_i$, such that 
$Z=\sum_i [\dv(f_i)]$. Two cycles $Z$ and $Z'$ are 
{rationally equivalent} in $X$ if $Z-Z'$ is rationally 
equivalent to $0$. We denote by $\CH(X)$ the quotient of $\mathcal Z(X)$ by
the subgroup of rationally trivial cycles.

\begin{emp} \label{fact-1} We will need the following facts. 
Given a  ring $A$ and an $A$-module $M$, we denote by $\ell_A(M)$ the {\it length} of $M$.
Let $(A, \m)$ be a noetherian local ring of dimension $1$. Let
$\p_1, \dots, \p_t$ be its minimal prime ideals. Let 
$\Frac(A)$ be the total ring of fractions of $A$ (with $\Frac(A)= \mathcal K_X(X)$ for $X=\Spec A$)
and denote by $\ord_A: \Frac(A)^* \to \mathbb Z$ the associated order function. Then  
\begin{enumerate}
\item Let $f\in \Frac(A)^*$, and let  $f_i$ denote the image of $f$ in $\Frac(A/{\p_i})^*$.
Using \cite{BLR}, Lemma 9.1/6, we obtain that
  $$\ord_A(f)=\sum_{1\le i\le t} \length_{A_{\p_i}}(A_{\p_i})\ord_{A/\p_i}(f_i).$$
\item If $A$ is reduced, then the canonical homomorphism 
$$\Frac(A)\longrightarrow \oplus_i A_{\p_i}= \oplus_i \Frac(A/\p_i)$$ 
is an isomorphism. 
\item Let $A \subseteq B \subseteq \Frac(A)$ be a subring such that $B/A$ is finite. 
Let $\mathfrak n_1, \dots, \mathfrak n_r$ be the maximal ideals of $B$, and let $b\in \Frac(A)^*$.
Then
$$\ord_A(b)=\sum_{1\le i\le r} [B/\mathfrak n_i : A/\m] 
\ord_{B_{\mathfrak n_i}}(b).$$
Indeed, our hypothesis implies that the $A$-module $B/A$ has finite length, so 
for a regular element $b \in A$, $\length_A(A/bA)= \length_A(B/bB)$.
Conclude using \cite{FulBook}, A.1.3, and the isomorphism  $B/bB \to \prod_{i=1}^r B_{\n_i}/bB_{\n_i}$.
\end{enumerate}
\end{emp} 

\begin{remark} \label{principal-div}
Let $\Gamma_1,\dots, \Gamma_r$ be the irreducible components of $X$, endowed with the reduced structure.
Let $\xi_1,\dots,\xi_r,$ denote their generic points.
Let $f\in \mathcal K_X^*(X)$, and 
let $f|_{\Gamma_i}$ be the meromorphic function restricted to $\Gamma_i$. 
If for all $i\le r$, $[\dv(f|_{\Gamma_i})]$ only involves codimension 
$1$ points of $X$, then $[\dv(f)]$ is rationally equivalent
to $0$ on $X$, since  $[\dv(f)] = \sum_{i=1}^r[\dv((f|_{\Gamma_i})^{ \length( \cO_{X,\xi_i})})]$ (use \ref{fact-1} (1)).

However, in general $[\dv(f)]$ \emph{is not} rationally 
equivalent to $0$. Consider for instance the  
projective variety $X$ over a field $k$, union of 
$\mathbb P^1_k$ and $\mathbb P^2_k$ intersecting 
transversally at a single point $\infty\in \mathbb P^1_k$. 
Let $x$ be a coordinate function on $\mathbb P^1_k$ with $[\dv_{\mathbb P^1_k}(x)]= [0] - [\infty]$.
Let $f$ be a rational function
on $X$ which restricts to $x$ on $\mathbb P^1_k$ and
is equal to $1$ on $\mathbb P^2_k$. Then 
$[\dv_X(f)]$ is the $0$-cycle $[0]$, since the point $\infty$ does not have codimension $1$ in $X$.
It is clear however that $[0]$ is not rationally trivial in $X$. 
This shows that the implication (1) $\Rightarrow$ (3) in 
the proposition in \cite{Ful79}, \S 1.8, does not hold in general. 
\end{remark} 

A proper morphism of schemes $\pi : Y\to X$ induces by {\it push forward of cycles} 
a group homomorphism $\pi_*: {\mathcal Z}(Y) \to
{\mathcal Z}(X)$. If $Z $ is any closed integral subscheme of $Y$, then 
$\pi_*([Z]):= [k(Z):k(\pi(Z))] [\pi(Z)]$, with the convention 
that  $[k(Z):k(\pi(Z))]=0$ if 
the extension $k(Z)/k(\pi(Z))$ is not algebraic.
It is known (\cite{K2}, 
\cite{Th}, and \ref{emp.proper} below) that
in general further assumptions are needed for a proper morphism $\pi$ to induce a group 
homomorphism $\pi_*: \CH(Y) \to
\CH(X)$. This is illustrated by our next example, also used later in \ref{ex.univcat}, \ref{rem.non-uc}, and  \ref{rem.non-uc2}.
(This example contradicts \cite{FulBook}, Example 20.1.3.)

\begin{example}\label{non-uc} 
We exhibit below a finite birational morphism $\pi : Y\to X$ of 
affine integral noetherian schemes with $Y$ regular, and a closed point 
$y_1 \in Y$ of codimension $1$ with $[y_1]$ rationally equivalent 
to $0$ on $Y$, but such 
that $\pi_*([y_1])$ is not rationally equivalent to $0$ on $X$.  The key  
feature in this example is that $\pi$ maps the point $y_1$ of codimension 
$1$ in $Y$ to a point of codimension $2$ in $X= \Spec A$. It turns out that $A$ 
is not universally catenary. Our example is
similar to that of \cite{EGA}, IV.5.6.11. The idea 
of the construction of a ring that is not universally catenary by gluing 
two closed points of distinct codimensions 
is due to Nagata (see \cite{Mat}, 14.E).

Let $k_0$ be any field. Let $k:=k_0(t_{\alpha})_{\alpha \in {\mathbb N}}$ be the field of 
rational functions with countably many variables. Consider the polynomial ring in one variable $k[S]$ and
the discrete valuation ring $R:=k[S]_{(Sk[S])}$. Let $Y:=\Spec R[T]$. Let $P(T)\in k[T]$ be an irreducible polynomial
of degree $d\ge 1$. Let $y_0\in Y$ be the closed point corresponding 
to $\p:=(P(T), S)$ and let $y_1$ be the closed point corresponding to 
$\q:=(ST-1)R[T]$. Then $\dim\ol{\{y_i\}}=0$ and $\dim\cO_{Y,y_i}=2-i$. 
The residue field $k(y_0):=k[T]/(P(T))$ is a finite
extension of $k$ of degree $d$, and $k(y_1)=\Frac(R)=k(S)$.  

Choose a field isomorphism  $\varphi : k\to k(S)$. Let 
$X:=\Spec A$ be the scheme obtained by identifying $y_1$ and $y_0$ via $\varphi$ (see \cite{Ped}, Teorema 1,
\cite{Sch}, 3.4, or \cite{Fer}, 5.4): 
$$A:=\{ f \in R[T] \mid f(y_0)\in k, \ \varphi(f(y_0))=f(y_1) \}.$$ 
By definition, $A$ is the pre-image of the field 
$\{ (\lambda, \varphi(\lambda)) \mid \lambda\in k\}$ under the canonical 
surjective homomorphism $R[T]\to k(y_0)\oplus k(y_1)$. The ideal
$\m:=\p\cap\q$ of $R[T]$ is then a maximal ideal of $A$, defining
a closed point $x_0\in X$ whose residue field $k(x_0)$ 
is isomorphic to $\{ (\lambda, \varphi(\lambda)) \mid \lambda\in k\}$. 
The inclusion $A\to R[T]$ induces a
morphism $\pi : Y\to X$.

It is easy to see that $R[T]/\m$ is finitely generated over $A/\m$. Since $\m\subseteq A$,
we can thus produce a finite system of generators for the $A$-module $R[T]$. More precisely, we have
$R[T]=A+T(TS-1)A+\dots+T^{d-1}(TS-1)A+TSA$. Therefore, $\pi$ is finite and, hence, 
$A$ is noetherian by Eakin-Nagata's theorem.  The ring $A$ 
has dimension 2 and, thus, is catenary.
The induced morphism $\pi: Y \setminus \{y_0,y_1\} \to X \setminus \{x_0\}$ 
is an isomorphism. Indeed, for any special open subset $D(h) \subseteq  X \setminus \{x_0\}$ (i.e., $h\in\m \setminus \{0\})$, 
we have $hR[T]\subseteq \m\subseteq A$. So any fraction $g/h^n$ with $g \in R[T]$ is equal to $gh/h^{n+1}$ with $gh \in A$, and  $R[T]_h=A_h$. 

Fix now $d\ge 2$. 
Let $f:=ST-1\in\q$. Then $[\dv(f)]=[y_1]$. By construction, $\pi$ induces 
an isomorphism $k(x_0)\simeq k(y_1)$, so that $\pi_*([y_1])=[x_0]$. 
We claim that $[x_0]$ is not rationally trivial on 
$X$. Indeed, let $C$ be a closed integral subscheme of $X$ containing $x_0$ as a point of codimension $1$. 
As $\cO_{X,x_0}$ and $X$ are of dimension $2$, we must have $\dim C =1$. 
Let $\tilde{C}$ be the schematic closure of $\pi^{-1}(C \setminus \{x_0\})$ in $Y$, 
and let $\rho : \tilde{C}\to C$ be the restriction of  $\pi$.
Then $\rho$ is a finite birational morphism of integral noetherian
schemes of dimension $1$. 
The point $y_1$ cannot belong to $\tilde{C}$, since otherwise the prime ideal 
$\q$ would properly contain the prime ideal of height $1$ corresponding to the generic point of 
$\tilde{C}$.
Hence, $\rho^{-1}(x_0)=\{y_0\}$. 

Now let 
$\dv(g)$ 
be a principal 
Cartier divisor on $C$. Then, using \ref{fact-1} (3),
$\ord_{x_0}(g)=[k(y_0) : k(x_0)]\ord_{y_0}(g)=d\ord_{y_0}(g)$. 
Therefore, if $n[x_0]$ is rationally equivalent 
to $0$, then $d\mid n$. It follows that 
$[x_0]$ is not rationally equivalent to $0$ when $d \geq 2$; in fact, $[x_0]$ has
order $d$ in the group $\CH(X)$. The same proof shows
that $[x_0]$ has order $d$ in the  group
$\CH(\Spec\cO_{X,x_0})$. 
\end{example}
\begin{emp} \label{grading} 
For general noetherian schemes, 
Thorup introduced a notion of rational equivalence
depending on a grading $\delta_X$ on $X$, which  turns the quotient $\CH(X, \delta_X)$ of ${\mathcal Z}(X)$ 
by this equivalence into a 
covariant functor for proper morphisms and a contravariant functor for flat 
equitranscendental morphisms (\cite{Th}, Proposition 6.5). 

We  briefly recall Thorup's theory  below. 
A \emph{grading} on a non-empty scheme $X$ is a map $\delta_X : X\to \mathbb Z$ such that
if $x\in \ol{\{ y\}}$, then $\Ht(x/y)\le \delta_X(y)-\delta_X(x)$ (\cite{Th}, 3.1). 
A grading $\delta_X$ is \emph{catenary} if the above inequality 
is always an equality (\cite{Th}, 3.6). An example of a grading on $X$ is the
canonical grading $\delta_{\mathrm{can}}(x):=-\dim\cO_{X,x}$. This 
grading is catenary if and only if $X$ is catenary and every local ring is 
equidimensional\footnote{
Recall that  a  ring $A$ of finite Krull dimension is equidimensional 
if $\dim A/\p=\dim A$ for every minimal prime ideal $\p$ of $A$.
A point $x \in X$ is equidimensional if ${\mathcal O}_{X,x}$ is.}
(\cite{Th}, p.\ 266)
at every point. 

Let $Y$ be an integral closed subscheme of $X$ with
generic point $\eta$, and let $f\in k(Y)^*$. Denote by $[\dv(f)]^{(1)}$ the cycle 
$[\dv(f)]$ where we discount all components $\ol{\{ x \}}$ such
that $\delta_X(x)< \delta_X(\eta)-1$. One defines the 
\emph{graded rational equivalence} on $\mathcal Z(X)$ using the subgroup generated
by the cycles 
$[\dv(f)]^{(1)}$, for all closed integral subschemes of $X$. If 
$\delta_X$ is catenary, then the graded rational equivalence is
the same as the usual (ungraded) one (\cite{Th}, Note 6.6). Denote by 
$\CH(X, \delta_X)$ the (graded) Chow group defined by the graded
rational equivalence. 

Let $f: Y\to X$ be a morphism essentially of finite type. Let $\delta_X$
be a grading on $X$. Then $f$ induces a grading $\delta_f$ on $Y$ 
defined in \cite{Th} (3.4), by 
$$\delta_f(y):=\delta_X(f(y))+{\mathrm{trdeg}}(k(y)/k(f(y))).$$
If $f$ is proper, then $f$ induces a homomorphism 
$f_{*}: \CH(Y, \delta_f)\to\CH(X, \delta_X)$ (\cite{Th}, Proposition 6.5). 
If $X$ is universally catenary and equidimensional at every point, and 
$\delta_X=\delta_{\mathrm{can}}$, then $\delta_f$ is a catenary 
grading on $Y$ (\cite{Th}, 3.11). It is also true that if $X$ is universally catenary
and $\delta_X$ is a catenary grading, then $\delta_f$ is a catenary 
grading on $Y$ (\cite{Th}, p. 266, second paragraph). 
\end{emp}
\begin{emp} \label{emp.proper}
In particular, 
assume that both $X/S$ and $Y/S$ are schemes of finite type over 
a noetherian scheme $S$ which is universally catenary and equidimensional at every point,
and $f:Y\to X$ is a proper morphism of $S$-schemes. Let $C$ and $C'$ be 
two cycles on $Y$ (classically) rationally equivalent. 
Then $f_*(C)$ and $f_*(C')$ are (classically) rationally equivalent on $X$.

In Example \ref{non-uc},  endow $X$ with the canonical grading, 
and $Y$ with the grading $\delta_\pi$. Then $\delta_X$ is catenary
but $\delta_\pi$ is not, because $y_1$ has virtual codimension $2$. 
Computations show that $\CH(Y, \delta_\pi)=\mathbb Z \oplus \mathbb Z$,
generated by the classes of $[y_1] $ and $[Y]$. The group 
$\CH(Y)$ is isomorphic to $\mathbb Z$, generated by the class of $[Y]$. 
The group
$\CH(X)$ is isomorphic to $\CH(X,\delta_{\mathrm{can}})=(\mathbb Z/d\mathbb Z) \oplus \mathbb Z$,
generated by the classes of $[x_0]$ and $[X]$, with the former of order $d$. 
\end{emp}
\begin{emp}
Let $S$ be a separated integral noetherian regular scheme of 
dimension at most $1$. Let $\eta $ denote its generic point. Endow $S$ 
with the catenary grading $1+\delta_{\mathrm{can}}$ (which is 
also the usual topological grading). Let $f : X\to S$ be a morphism 
of finite type, and endow $X$ with the grading $\delta_f$. This 
grading is catenary (\ref{grading}).  

Let $n\ge 1$ and let $x\in X$ be such that $\delta_f(x)=n$. 
Then $[\overline{\{x\}}]$ is an {\it $n$-cycle} on $(X,\delta_f)$. 
If $\dim(S)=1$ and $f(x)$ is a closed point $s\in S$, then $\ol{\{x\}}$ is a 
subscheme of dimension $n$ of the fiber $X_s$. If $f(x)=\eta$, 
then $\ol{\{x\}}\to S$ is dominant and $\dim\ol{\{x\}}_{\eta}=n-1$. 
In the latter case, $\dim\ol{\{x\}}=n-1$ if and only if $S$ is 
semi-local and $\ol{\{x\}}$ is contained in $X_{\eta}$. Otherwise, 
$\dim\ol{\{x\}}=n$. 

In particular, the irreducible $1$-cycles on 
$(X,\delta_f)$ are of two types: the integral closed subschemes 
$C$ of $X$ of dimension $1$ such that $C$ meets at least one closed 
fiber, and the closed points of $X$ contained in $X_{\eta}$ (in which case
$S$ must be semi-local). We say that a $1$-cycle is {\it horizontal}
if its support is quasi-finite over $S$, and that it is {\it vertical} 
if its support is not dominant over $S$.
\end{emp}

\end{section}

\begin{section}{Moving Lemma for $1$-cycles on regular $X/S$ with $S$ semi-local}

Let $X$ be a quasi-projective scheme of pure dimension $d$
a field $k$.
Let $X^{\mathrm{sing}}$ denote the non-smooth locus of $X$.
The classical Chow's Moving Lemma \cite{Rob} and its generalization 
(\cite{Con}, II.9, assuming $k$ algebraically closed)
immediately imply the following statement:

\begin{emp} \label{Chow} 
{\it 
Let $0 \leq r \leq d$. Let $Z$ be a $r$-cycle on $X$ with $\Supp(Z) \cap X^{\mathrm{sing}} = \emptyset$.
Assume that $\dim(X^{\mathrm{sing}}) < d-r$.
Let $F$ be a closed subset of $X$ of codimension at least $r+1$ in $X$. 
Then there exists an $r$-cycle $Z'$ on $X$, rationally equivalent to $Z$, 
and such that $\Supp(Z') \cap (F \cup X^{\mathrm{sing}}) = \emptyset$.
}
\end{emp}

Our goal in this section is to prove a variant of this statement for a scheme $X$
over a semi-local Dedekind  base $S=\Spec R$. 
An application of such a relative moving lemma is given in Theorem \ref{dx-ds}.

\begin{emp}\label{def-FA} \label{FA-g} Let $X$ be a scheme. We say that 
$X$ {\it is an $\FA$-scheme}, or simply that $X$ is $\FA$,  
if every finite subset of $X$ is contained in  an 
affine open subset of $X$. In particular, an $\FA$-scheme is separated. 
The following examples of $\FA$-schemes
are well-known:

\smallskip
(1) Any  affine scheme is $\FA$. Any 
quasi-projective scheme over an affine scheme is $\FA$
 (\cite{Liubook}, Proposition 3.3.36). 
More generally, a scheme admitting an ample invertible 
sheaf is $\FA$ (\cite{EGA}, II.4.5.4). 

\smallskip
(2) If $X$ is $\FA$, then any closed subscheme of $X$ is clearly
$\FA$. The same holds for any open subset $U$ of $X$. Indeed, 
let $F$ be a finite subset of $U$, then $F$ is 
contained in an affine open subset $V$ of $X$. Hence, 
$F\subseteq U\cap V$ with $U \cap V$ quasi-affine. By (1), $F$ is contained
in an affine open subset of $U\cap V$. 

\smallskip
(3) More generally, if $Y$ is $\FA$ and $f : X\to Y$ is a morphism  
of finite type admitting a relatively ample invertible sheaf, then 
$X$ is $\FA$. Indeed, any finite subset of $X$ has finite image in
$Y$, so we can suppose that $Y$ is affine. Then $X$ admits an
ample invertible sheaf \cite{EGA}, II.4.6.6, and we are reduced
to the case (1).

\smallskip
(4) A noetherian separated scheme  of dimension $1$ is $\FA$ (\cite{RayFA}, Prop. VIII.1).
\end{emp} 
Suppose $k$ is an algebraically closed field, and that $X/k$ is a regular 
$\FA$-scheme of finite type. Let $S/k$ be a separated scheme of finite type. Then any proper $k$-morphism $X \to S$ 
is projective (\cite{K66}, Cor. 2).

For the purpose of our next theorem, we will call a noetherian integral domain
$R$ a {\it Dedekind domain} if it is integrally closed of dimension $0$ or $1$.
A version of this theorem where $R$ is not assumed to be semi-local is proved
in \cite{GLL2}, 7.2.

\begin{theorem} \label{mv-slocal} Let $S$ be the spectrum 
of a semi-local Dedekind domain $R$.
Let $f:X \to S$ be a separated morphism 
of finite type, with $X$ regular and $\FA$. 
Let $C$ be a horizontal $1$-cycle on $X$ with 
$\Supp(C)$ finite
over $S$.
Let $F$ be a closed subset of $X$ such that for every $s \in S$, 
any irreducible component of $F \cap X_s$ that meets $C$ is not an irreducible 
component of $X_s$. 
Then there exists a 
horizontal $1$-cycle $C'$ on $X$ with $f_{|C'}: \Supp(C') \to S$ finite, 
rationally equivalent to $C$,
and 
such that 
 $\Supp(C')\cap F = \emptyset$.

In addition, since  $S$ is semi-local, $C$ consists of finitely
many points, and since $X$ is $\FA$, there exists an affine open subset $V$ of $X$ which contains $C$.
 Then, for any such open subset $V$, the 
horizontal $1$-cycle $C'$ can be chosen to be contained in $V$, 
and to be such that if $g:Y \to S$ is any separated morphism of finite type with an open embedding
$V \to Y$ over $S$, then $C$ and $C'$ are closed  and rationally equivalent on $Y$.
\end{theorem} 

\proof It suffices to prove the theorem when $C$ is irreducible and $\Supp(C) \cap F \not = \emptyset$.
Choose an affine open subset $V$ of $X$ 
containing $C$. Since $C$ is closed in $V$, it is affine.

Proposition \ref{normalization-L-global} 
shows the existence of a finite birational morphism
 $D\to C$ such that the composition $D\to C \to S$ is a local complete intersection morphism (l.c.i).
 Clearly, when $S$ is excellent, we can take $D\to C$ to be the 
normalization morphism, in which case $D$ is even regular, and \ref{normalization-L-global}  is not needed.
Since $C$ is affine and $D \to C$ is finite,  
there exists for some $N \in {\mathbb N}$ a closed immersion 
$D \to  C\times_S \mathbb A^N_S\subseteq V\times_S \mathbb A^N_S$.

Let $U:= V\times_S \mathbb A^N_S$.  We claim that it suffices to prove  the theorem 
for the $1$-cycle $D$ and the closed subset ${\bf F}:= F \times_S \mathbb A^N_S$ in the affine scheme $f':U\to S$.
Indeed, let $D'$ be a horizontal $1$-cycle whose existence is asserted by the theorem in this case.
In particular, $\Supp(D') \cap {\bf F} = \emptyset$.
Let $V \to Y$ be any open immersion over $S$. 
Consider the associated open immersion $U \to Y\times_S \mathbb P^N_S$ and 
the projection $p : Y\times_S \mathbb P^N_S\to Y$. By hypothesis,
$D$ and $D'$ are closed and rationally equivalent in $Y\times_S \mathbb P^N_S$.
One easily checks that $p_*(D) = C$ because $D\to C$ is birational.
It follows from 
\ref{emp.proper} that $p_*(D) = C$ is rationally equivalent to $C':= p_*(D')$ on $Y$.
Moreover, $\Supp(C') \cap F = \emptyset$.

The existence of $D'$ with the required properties follows from  
Proposition \ref{mv-d} below. Indeed, 
first note that since $D/S$ is l.c.i., each local ring $\cO_{D,x}$, $x \in D$,
is an absolute complete intersection ring (\cite{EGA}, IV.19.3.2). 
It follows that the closed immersion $D \to  U$ is a regular immersion (\cite{EGA}, IV.19.3.2).

Let $d:=\codim(D, U)$. We note that $d>0$ since $\Supp(C) \cap F \not = \emptyset$
and for each point in $\Supp(C) \cap F$ over $s \in S$, $F \cap V_s$ 
is not an irreducible component of $V_s$. 
Let $x\in D$ be a closed point, and let $s:=f'(x)$.
Then $\dim\cO_{D,x}=\dim(S)$, $\dim\cO_{U,x}=d+\dim(S)$,  and $\dim\cO_{U_s,x}=d$.
Our assumption on $F$ implies that 
the irreducible components of ${\bf F}\cap U_s$ passing through
$x$ have dimension at most $ d-1$. We can thus apply \ref{mv-d} below to conclude the proof of \ref{mv-slocal}.
\qed

\begin{proposition}\label{mv-d} Let $S$ be any semi-local affine noetherian scheme.
Let $U \to S$ be a morphism of finite type with $U$ affine.
Let $C$ be an integral closed subscheme of $U$, 
of codimension $d\ge 1$, and finite over $S$. 
Suppose that the closed immersion $C\to U$ is regular.
Let $F$ be a closed subset of $U$ such that 
for all closed points $s\in S$, 
the irreducible components of $F\cap U_s$ that intersect 
$C$
all have dimension at most $ d-1$. 
Then there exists a cycle $C'$ on $U$ rationally equivalent to $C$ and such that: 
\begin{enumerate}[\rm (1)] 
\item The support 
of $C'$ is finite
over $S$ and 
 does not meet $F \cup C$. 
Moreover, for any closed point $s \in S$, 
$\Supp(C')$ does not contain any irreducible component of $U_s$.
\item Suppose that $S$ is universally catenary. 
Let $Y\to S$ be any separated morphism of finite type
and let 
 $h: U\to Y$ be any $S$-morphism. Then $h_*(C)$ is rationally equivalent 
to $h_*(C')$ on $Y$. 
\end{enumerate}
\end{proposition} 

\proof {\it Reduction to the case $d=1$}. Write $U:=\Spec A$ and $C:=V(J)$.
Suppose $d\ge 2$. 
By hypothesis, the $(A/J)$-module $J/J^2$ 
is locally free, hence free of rank $d$. 
Now lift a basis of $J/J^2$ to elements 
$f_1,\dots, f_d\in J$. For all $\p\in C$, 
we have $J_\p=f_1A_\p+\dots+f_dA_\p$ by Nakayama's Lemma. 
As $J_\p$ is generated by a regular sequence by hypothesis, 
Lemma \ref{reg-seq} implies that $f_1,\dots, f_d$ is a regular 
sequence in $A_\p$. 

Let $\Gamma_1, \dots, \Gamma_n$ denote 
the irreducible components of $F\cap U_s$ 
that intersect $C$, with $s$ ranging through the finitely many closed points of $S$.
By hypothesis, $\dim \Gamma_i \leq d-1$ for all $i$.
Apply Lemma \ref{cut-d} to $A$, $J$,  
$\Gamma_1, \dots, \Gamma_n$, and $f_1,\dots, f_{d-1}$ as above.
We obtain the existence of $g_1,\dots, g_{d-1}$, such that 
$C \subseteq V(g_1, \dots, g_{d-1})$, and such that
every irreducible component of 
$V(g_1,\dots, g_{d-1})\cap \Gamma_i$ either has dimension $0$ 
or is contained in $C$. 

Lemma \ref{reg-seq} 
implies that $g_1,\dots, g_{d-1},f_d$ is a regular sequence at the points
of $C$, so the immersion $C \to V(g_1,\dots, g_{d-1})$ is regular.
It is easy to check that the proposition is proved
if it can be proved for the closed subsets $C$ and $F \cap V(g_1,\dots, g_{d-1})$ 
inside the affine scheme $V(g_1,\dots, g_{d-1})$.
Moreover, we note that now
\begin{enumerate}[\rm (a)]
\item  at any point $\p$ of $C$, $C$ is defined in $V(g_1,\dots, g_{d-1})$
by the regular element $f_d$, and 
\item   
any irreducible component of $F\cap V(g_1,\dots, g_{d-1})_s$ is either contained in $C$ 
or disjoint from $C$.
\end{enumerate}

We make one further reduction if  $V:=V(g_1,\dots, g_{d-1})$ is not integral.
For every $x \in C$, $f_d$ is a regular element of $\cO_{V,x}$, and $\cO_{C,x} \simeq \cO_{V,x}/(f_d)$.
This easily implies that  $\cO_{V,x}$ is a domain. Thus, there is a unique irreducible component $W$ of $V$
through $x$. Clearly, $W$ is independent of $x$, and in a neighborhood of $C$, $V$ coincides with $W$ endowed with the reduced subscheme structure. The morphism $C \to W$ is still regular. We may thus replace $V$ with the integral subscheme $W$.

\smallskip 

{\it We assume henceforth that $U$ is integral and that $d =1$}. 
Fix now an open $S$-immersion $U \to X $, with $X/S$ projective and $X$ integral.
Let $\overline{F}$ be the Zariski closure of $F$ in $X$.
Let $Z$ denote the  closed subset of $X$ consisting in the finite union of 
the following closed sets:
\begin{enumerate}[{\rm (a)}]
\item  $X\setminus U$;
\item  
All 
irreducible components of 
$\overline{F}\cap X_{s}$ 
which do not intersect $C$, for each closed point $s \in S$; and
\item 
One closed point of $\Gamma$ which does 
not belong to $C$, for each
irreducible component $\Gamma$  of $U_{s}$ which is not contained in $C$, and for each closed point $s \in S$.
\end{enumerate}
By construction, $Z \cap C = \emptyset$.

Since $C$ is proper over $S$, it is  closed in $X$.
Since $d=1$, $C$ is in fact a Cartier divisor on $U$, and since it is closed 
in $X$, we can extend it to a Cartier divisor on $X$. Let $\mathcal J$ be the sheaf of ideals on $X$ 
defining $C$. This is then an invertible sheaf, and as usual we let $\cO_X(nC):= {\mathcal J}^{-n}$.

Let $\mathcal I$ be the sheaf 
of ideals on $X$ defining the reduced induced structure on 
$Z$. Let $\mathcal I(nC)$ denote the image of 
$\mathcal I\otimes \cO_{X}(nC)$ in $\cO_{X}(nC)$. 
For all $n\ge 2$, we have a natural exact sequence 
$$ 0 \to \mathcal I((n-1)C)\to \mathcal I(nC)\to \mathcal F_n \to 0,$$  
where $\mathcal F_n$ is a coherent sheaf annihilated by 
$\mathcal J$ and, hence, supported on $C$.  
Applying $H^1(X,-)$ to the above exact sequence, we get 
$$H^1(X, \mathcal I((n-1)C))\to H^1(X, \mathcal I(nC))\to 0$$ 
because $C$ is affine. Since $X/S$ is projective, we have obtained a system of finitely generated 
$\cO_S(S)$-modules with surjective transition maps. Since   $\cO_S(S)$ is noetherian, 
the transition maps are eventually isomorphisms. Hence, there exists $n_0 \geq 0$ such that
for all  $n\ge n_0$, 
$H^0(X, \mathcal I(nC))\to H^0(C, \mathcal F_n)$ 
is surjective. 

By hypothesis, the stalk $\mathcal J_x$ at each $x \in C$ is generated 
by a regular element. Since $C$ is semi-local and closed in the affine scheme $U$, we can find 
an affine open subset $V$ in $X\setminus Z$ containing $C$ such that ${\mathcal J}_{|V}$ is principal, say generated by a function $\varphi$.
For all $n \geq 0$, $\varphi^{-n}$ induces a generator $\overline{\varphi^{-n}} $ of $\cO(nC)_{|C}$. 
Since $Z \cap C = \emptyset$,
we find for all $n \geq 0$ that $({\mathcal F}_n)_{|C}$ is isomorphic to $\cO(nC)_{|C}$.
Thus, for all $n \geq n_0$, we can find a global section $f_n$ of 
$\mathcal I(nC)$ which lifts the generator $\overline{\varphi^{-n}}$ of $\cO_{X}(nC)_{|C}$. 

Fix $n \geq n_0$, and consider $g:=(1+f_{n+1})/(1+f_n)\in {\mathcal K}_{X}(X)$. 
Let 
$$C':=C+[\dv_X(g)].$$ 
Then $\Supp(C')$ is projective
over $S$ because it is closed in $X$.
By construction, $[\dv_{V}(g)] = -C$, and $[\dv_X(g)]$ has support disjoint from $Z$.
Thus,
$C\cap \Supp(C')=\emptyset$,
and since $Z$ contains $X\setminus U$, we find that $\Supp(C')\subseteq U$.
Since $U$ is affine, $\Supp(C')$ is finite over $S$.

Recall that
 for each closed point $s \in S$, $Z$ contains
all 
irreducible components of 
$\bar{F}\cap X_{s}$ 
which do not intersect $C$.
Recall also that by hypothesis (in the case $d=1$), 
for all closed points $s\in S$, 
the irreducible components of $F\cap U_s$ that intersect 
$C$
have dimension $0$ and are thus contained in $C$.
It follows that  $F\cap \Supp(C')= \emptyset$.

Finally, recall that for each closed point $s \in S$ and for each
irreducible component $\Gamma$  of $U_{s}$ which is not contained in $C$, 
then $Z$ contains a closed point of $\Gamma$ which does 
not belong to $C$. Then $\Supp(C')$ does not 
contain any irreducible component of $U_s$ which is not contained in $C$. This shows (1).

(2) Now suppose that $S$ is universally catenary. 
We start with the following reduction. Recall from the beginning of the proof 
the existence of a closed integral subscheme $W$ of $U$ such that $C \subset W$ is a regular embedding
and $C$ has codimension $1$ in $W$. Let $T$ denote the schematic closure of the image of $W$ in $S$.
It suffices to prove (2) for  $W \to T$, and the morphism $h': W \to Y\times_{S} T$. 
We are thus reduced to the case where both $U$ and  $S$ are integral. In particular, $S$ is equidimensional at every point
and   universally catenary, and the theory recalled in \ref{grading} applies.

Fix as in (1) an open $S$-immersion $U \to X $, with $X/S$ projective and $X$ integral.
Let $g$ be as in (1).
Let $\Gamma\subseteq X\times_S Y$ be the schematic closure of the graph 
of the rational map $X\dasharrow Y$ induced by $h:U \to Y$. Let 
$p: \Gamma \to X $ and $q: \Gamma\to Y$ be the associated projection maps over $S$. 
Since $Y/S$ is separated, the graph of $h:U \to Y$ is closed in $U \times_S Y$. 
Hence, $p : p^{-1}(U)\to U$ is an isomorphism.
Since $\Gamma$ is integral and its generic point maps to the generic 
point of $X$, 
the rational function $g$ on $X$ induces a rational function, again denoted by $g$, on $\Gamma$.
As $p : p^{-1}(U)\to U$ is an isomorphism,
we let $p^*(C)$ and $p^*(C')$ denote the preimages of $C$ and $C'$ in $p^{-1}(U)$; they are 
closed subschemes of $\Gamma$.
Since $g$ is an  invertible function in a neighborhood of $X \setminus U$,
$[\dv_\Gamma(g)] =p^*(C)-p^*(C')$,
and $p^*(C)$ and $p^*(C')$ are rationally equivalent on $\Gamma$. 
Fix a catenary grading on $S$ and define gradings on schemes 
of finite type over $S$ accordingly (\ref{grading}).  Then, as $q$ is 
proper, $q_*p^*C$ and $q_*p^*C'$ are rationally equivalent in $Y$.
Since $h_*C =q_*p^*C$ and $h_*C'= q_*p^*C'$, (2) follows. 
\qed

\begin{lemma} \label{cut-d} 
Let $U=\Spec A$ be a noetherian affine scheme. Let $C:=V(J)$ be 
a closed subset 
of $U$. Let $\Gamma_1,\dots, \Gamma_n$ be irreducible closed subsets 
of $U$. 
 Let $ f_1,\dots, f_{\delta} \in  J$. Then there exist $g_1,\dots, g_{\delta}\in J$ such that 
 $g_i\in f_i+J^2$ for all $i=1,\dots, \delta$, and such that the following property holds.
 Let $i\le \delta$ and $j\le n$. Then any irreducible 
component of $\Gamma_j\cap V(g_1,\dots, g_{i})$ not 
contained in $C$ has codimension $i$ in $\Gamma_j$ and, hence, dimension at most $ \dim \Gamma_j -i$. 
\end{lemma}

\proof 
For $j=1, \dots,n$, let $\q_j$ be the prime ideal of $A$ such that $\Gamma_j=V(\q_j)$.
If $\q_j$ contains $J$ for all $j\le n$, we set $g_i := f_i$ for all $i=1,\dots, \delta$, 
and the lemma is proved.  
Suppose that for some $j$, $\q_j$ does not contain $J$.
Upon renumbering if necessary, assume that the ideals $\q_1, \dots, \q_m$ do not contain $J$,
and $\q_{m+1},\dots, \q_n$ contain $J$. The lemma is proved
if we can prove it for the sets $\Gamma_1, \dots, \Gamma_m$.
We may thus assume that none of the $\Gamma_j$'s is contained in $C$ or, in other words,
that none of the $\q_j$'s contain $J$.

We proceed by induction on $\delta$. When $\delta=1$, 
we find that $f_1A+J^2 \not\subseteq \q_j$ for all $j\le n$. 
Then there exists $a_1 \in J^2$ 
such that \ $g_1:=f_1+a_1\notin\cup_{1\le j\le n} \q_j$  (\cite{B-H}, Lemma 1.2.2 or 
\cite{Kp}, Theorem 124, page 90). Suppose that $\Theta$ is an
irreducible component of $\Gamma_j\cap V(g_1)$. Then  $\Theta$ has codimension $1$ in $\Gamma_j$ and
$\dim\Theta\le \dim \Gamma_j\cap V(g_1) 
\le \dim \Gamma_j-1$, since $A/\q_j$ is a domain, and $g_1 \notin \q_j$.
 
If $\delta\ge 2$, we apply the induction hypothesis 
to the sequence $f_1,\dots, f_{\delta-1}$ to obtain the desired $g_1,\dots, g_{\delta-1}$.
We then apply the case $\delta=1$ to the 
ring $A/(g_1, \dots, g_{\delta-1})$, the ideal $J/(g_1, \dots, g_{\delta -1})$, the image
of $f_{\delta}$ in $A/(g_1, \dots, g_{\delta-1})$, and to the irreducible components
of the $\Gamma_j\cap V(g_1, \dots, g_{\delta-1})$'s which are not contained in $C$.
We find then an element $\bar{g}_{\delta}$ in $\bar{f}_{\delta}+ (J/(g_1, \dots, g_{\delta -1}))^2$, 
which we lift to $g'_{\delta}= f_{\delta}+g_1a_1+\dots +g_{\delta-1}a_{\delta-1}+j$ with $j \in J^2$ and $a_i \in A$.
Since the desired property is now achieved for the irreducible components
of $\Gamma_j \cap V(g_1,\dots,g_{\delta-1},g'_{\delta})$ not 
contained in $C$, we find that  the sequence $g_1,\dots,g_{\delta-1},g_{\delta}$, with $g_{\delta}:=f_{\delta}+j$,
satisfies the conclusion of the lemma. 
\qed 
\begin{lemma} \label{reg-seq} Let $A$ be a noetherian local ring. 
Let $I$ be a proper ideal of $A$ 
generated by a regular sequence 
$f_1,\dots, f_d$. Let $g_1,\dots, g_d\in I$. If the image of $\{ g_1,\dots, g_d\}$ in 
 $I/I^2$ is a basis of $I/I^2$ over $A/I$, then $g_1,\dots, g_d$ is a regular 
sequence. 
\end{lemma}

\proof
This is well-known, and follows from
the equivalence between quasi-regular sequences 
and regular sequences in noetherian local rings 
(\cite{Mat}, 15.B, Theorem 27). 
\qed 

\begin{remark} \label{ex.univcat} We  note that  the cycle $C'$ in Proposition \ref{mv-d} (1) is in the same graded component  of ${\mathcal Z}(U)$ as $C$,  for every catenary grading on $U$.
Note also that it may happen in  \ref{mv-d} (1) that the cycle $C'$ is the trivial cycle, with empty support.
Indeed, consider the case where $U \to S$ is a finite morphism. According to \ref{mv-d}  (1), the support of $C'$ does not meet $C$, and for any closed point $s \in S$, 
${\rm Supp}(C')$ does not contain any point in   $U_s \setminus C $. It follows that ${\rm Supp}(C')$ does not contain 
any closed point of $U$. Since $U$ is affine, we find that ${\rm Supp}(C')$ is empty.

We now modify Example \ref{non-uc}  to produce an example to show that the statement of Proposition \ref{mv-d} (2) does not hold in general 
if $S$ is not assumed to 
be universally catenary. Indeed, with the notation as in Example \ref{non-uc}, let $d>1$, let $Y'$ be the semi-localization of $Y$ at $\{y_0,y_1\}$, and let $X'$ the pinching of $Y'$ (that is, the scheme obtained by identifying $y_0$ and $y_1$ in $Y'$). Let again $\pi$ denote the 
natural finite morphism $Y' \to X'$. 
Then the cycle $[y_1]$ is rationally trivial on $Y'$, but the cycle $\pi_*([y_1])$ is   not rationally trivial on $X'$. For the desired example where the statement of Proposition \ref{mv-d} (2) does not hold, we take $S$ and $Y$ to be $X'$, $U$ to be  $Y'$, and $C$ to be $\{y_1\}$.

\end{remark}

\end{section}

\begin{section}{Inductive limits of l.c.i. algebras} 

We prove in this section a technical statement 
needed when considering base schemes $S$ which are not excellent,
as in the proof of \ref{mv-slocal}, or of 6.2 
 in \cite{GLL2}.

\begin{lemma}\label{normalization-L} 
Let $(R,(\pi))$ be a discrete valuation ring with
field of fractions $K$ of characteristic $p>0$. Let $L/K$
be a purely inseparable extension of degree $p$, and let $R_L$ be the
integral closure of $R$ in $L$. Then 
there exists a sequence $(\alpha_n)_{n \in {\mathbb N}^*}$ with $\alpha_n \in R_L$ and $L=K(\alpha_n)$, such that
the $R$-algebras
$B_n:=R[\alpha_n]$ are finite and  local complete intersections 
over $R$,  and $R_L = \cup_n B_n$. 
\end{lemma}

\proof Recall that since $[L:K]=p$, 
 $R_L$ is also a discrete valuation ring, with maximal ideal $\m$, and
either $e_{\m/(\pi)}f_{\m/(\pi)}=1$ and $R_L$ is not a finitely generated $R$-module,
or $e_{\m/(\pi)}f_{\m/(\pi)}=p$ and $R_L$ is a finitely generated $R$-module (see \cite{Bou}, VI.8.5, Theorem 2). Moreover, 
if $e_{\m/(\pi)}=p$ and $s \in \m \setminus \m^2$, then $R_L=R[s]$.
If $f_{\m/(\pi)}=p$ and $s \in R_L$ reduces modulo $\m$ to a $t$ such that $R_L/\m=(R/(\pi))(t)$, 
then $R_L=R[s]$.

The lemma is thus completely proved in the case where $R_L$ is a finitely generated $R$-module 
by setting $B_n=R_L$ for all $n$.
When $e_{\m/(\pi)}f_{\m/(\pi)}=1$, we can embed $R_L$ into the completion of $\widehat{R}$ of $R$ with 
respect to $(\pi)$.
Since $[L:K]=p$, we can choose $\alpha\in R_L$ such that $L=K(\alpha)$.  
For any $n\ge 1$, we can approximate the element $\alpha$ by an element  $r_n\in R$, such that
$\alpha=r_n+\pi^n \alpha_n$ for some $\alpha_n\in R_L$. Let $B_n:=R[\alpha_n]$. Then 
$(B_n)_{n\ge 1}$ is an increasing sequence of finite l.c.i. algebras
over $R$, and $L=K(\alpha_n)$ for all $n$. 

It remains to show that $R_L=\cup_n B_n$.
Let $b\in R_L$. Then  
$b=\sum_{j=0}^{ p-1} t_{j}\alpha^j$, with $t_j\in K$ for $j=0,\dots, p-1$. Replacing  $\alpha$
by $r_n+\pi^n \alpha_n$, we see that for $n> \max_{j \geq 1} |\ord_{\pi}(t_j)|$, we can write 
$b=c_n+\beta_n$ with $c_n=\sum_j t_j r_n^j\in K$ and 
$\beta_n \in R[\alpha_n]$.  It follows  
that $c_n\in R_L$ and, hence,  $c_n \in R=K\cap R_L$. Therefore, $b\in R[\alpha_n]$.
\qed 

\begin{proposition}\label{normalization-L-global} 
Let $A$ be a Dedekind domain,
with field of fractions $K$. Let $B$ be an integral domain 
containing $A$, and with field of fractions $L$.
Assume that $B$ is finite over $A$. Then there exists 
a domain $C$ with $B \subseteq C \subseteq L$  
such that $C$ is finite over $A$, and a local complete 
intersection over $A$.
\end{proposition}

\proof Let $K_{\sep}$ be the separable closure of $K$ in $L$ and
let $B_{\sep}$ be the integral closure of $A$ in $K_{\sep}$. Then
$B_{\sep}$ is a Dedekind domain finite over $A$. Let 
$B'$ be the sub-$A$-algebra of $L$ generated by $B$ and $B_{\sep}$.
Then it is finite over $B_{\sep}$. If we can find $C$ containing 
$B'$, finite and l.c.i. over $B_{\sep}$, then $C$
contains $B$ and is finite and l.c.i. over $A$. Therefore, 
it is enough to treat the case when $K=K_{\sep}$. We thus assume now that $L$ is purely 
inseparable over $K$, and we prove the proposition by induction on 
$[L:K]$. The case $[L:K]=1$ is trivial. Suppose $[L:K]>1$ and that 
the proposition is true for any Dedekind domain $A'$ and for 
any purely inseparable extension $L'$ of $\Frac(A')$ of 
degree strictly less than $[L:K]$. 

We will construct $C$ as the global sections of a coherent sheaf of
$A$-algebras $\mathcal C$ over $\Spec A$. 
Let $L'$ be a subextension of $L$ of index $p$. Consider the 
sub-$A$-algebra $B\cap L'$ of $L'$. It is finite over $A$ because
$B$ is finite over $A$. By induction 
hypothesis, it is contained in a finite l.c.i. $A$-algebra
$B'\subseteq L'$. Let $\alpha\in B$ be such that $L=L'[\alpha]$. 
We claim that $B'[\alpha]$ is l.c.i. over $B'$, 
so that   $B'[\alpha]$ is a finite l.c.i.\ $A$-algebra. Indeed,
the minimal polynomial of $\alpha$ over $L'$ is $x^p - \alpha^p \in B'[x]$,
since by hypothesis, $\alpha^p \in L' \cap B \subseteq B'$. Therefore, 
the map $B'[x]/(x^p - \alpha^p) \to B'[\alpha]$ is an isomorphism.
 
As $B'[\alpha]$ and $B$ are both finite over $A$ and have the 
same field of fractions, there exists $f\in A \setminus\{0\}$ such 
that $B'[\alpha]_f=B_f$ (where as usual $B_f$ denotes the localization of $B$
with respect to the multiplicative set $\{1,f,f^2,\dots\}$). The algebra $B_f$ is finite and l.c.i.\ over
$A_f$, and we define ${\mathcal  C}(D(f))$ to be $B_f$. 

Let now $\p$ be any maximal 
ideal of $A$, and let 
$\q$ be the unique maximal ideal of $B$ lying over $\p$. 
Let us show that 
$B_{\q}$ is contained in some finite l.c.i.\ $A_\p$-algebra ${\mathcal  C}_\p$ contained in $L$.  
First, note that 
the localization of $B$ at the multiplicative set $A\setminus \p$ 
is equal in $L$ to $B_{\q}$, 
because $x^{[L:K]} \in A \setminus \p$ for all $x \in B \setminus \q$.
In particular, $B_{\q}$ is finite over $A_{\p}$. 
Write $B_{\q}=\sum_{1\le i\le r}b_i A_\p$, with $b_i \in B$. Let $R$ be the integral 
closure of $A_\p$ in $L'$. 
It follows from Lemma \ref{normalization-L} that there exists
$\alpha\in L$, integral over $R$ with $L=L'(\alpha)$, 
and such that $b_i\in R[\alpha]$ for all $i\le r$.  
Then there also exists a finite sub-$A_\p$-algebra $D$ of $L'$,
with $\alpha^p \in D$ and  $b_i\in D[\alpha]$ for all $i\le r$, and such that
$\Frac(D)=L'$. 
Apply now the induction hypothesis to the extension 
$A_\p\subseteq D$: there
exists $D\subseteq E \subseteq R$ such that $E$ is finite and
l.c.i. over $A_\p$. Set ${\mathcal  C}_\p:=E[\alpha]$. We have 
$\Frac({\mathcal  C}_\p)=L'[\alpha]=L$. Since $\alpha^p \in D$, the 
minimal polynomial of $\alpha$ over $L'$ is in $D[x]$.  As earlier, we find that  $E[\alpha]$
is l.c.i. over $E$ and, therefore, also over $A_\p$.

Write $\Spec A = D(f) \cup \{\p_1, \ldots, \p_s\}$.  
For each $i=1,\dots, s$, there exists an open neighborhood $U_i$ of $\p_i$ 
with $U_i\setminus \{ \p_i\} \subseteq D(f)$, and a 
finite sub-$\cO(U_i)$-algebra $C_i$ of 
$L$ such that $C_i\otimes_{\cO(U_i)} A_{\p_i}={\mathcal  C}_{\p_i}$ and 
$$C_i\otimes_{\cO(U_i)} \cO(U_i\setminus \{ \p_i\})=
B_f\otimes_{A_f}\cO(U_i\setminus \{\p_i\}).$$ We then can glue the sheaves $(C_i)^{\sim}{}_{/U_i}$
and the sheaf $(B_f)^{\sim}{}_{/D(f)}$  
into a coherent sheaf $\mathcal C$ of $A$-algebras. The
global sections $C:={\mathcal  C}(\Spec A)$ is a finite sub-$A$-algebra 
of $L$ which contains  $B$ and is l.c.i. by construction. 
\qed 

\end{section}

\begin{section}{Hilbert-Samuel multiplicities} \label{section-local}

Let $(A, \m)$ be a noetherian local ring. Let $Q$ be an
$\m$-primary ideal of $A$
(equivalently, $Q$ is a proper ideal of $A$ containing 
some power of $\m$). 
Let $M$ be a non-zero finitely
generated $A$-module.
Recall that there exists a polynomial 
$f_Q(x)\in \mathbb Q[x]$ 
such that for all $n$ large enough, $f_Q(n) = \length_A(M/Q^nM)$. 
This polynomial has degree $d=\dim M:=\dim V(\Ann(M))$.

The {\it Hilbert-Samuel multiplicity} $e(Q, M)$ 
is the coefficient of $x^d$ in $f_Q(x)$ multiplied by $d!$. 
When there is no need to specify 
the ring $A$, we may write $e(Q, A)$ simply as $e(Q)$.
The integer $e(A):=e(\m, A)$ is called the Hilbert-Samuel  
multiplicity of $A$. If $A$ is regular, then $e(A)=1$.

We prove in Proposition \ref{eq-sp} below that the Hilbert-Samuel multiplicity $e(Q, M)$
can be expressed in terms of Hilbert-Samuel multiplicities $e(({\bf f}), M)$ of ideals generated 
by strict systems of parameters. This result is a key ingredient in the proof of the main result of this section,
Theorem \ref{local-mv}. 
A more global geometric version of Theorem \ref{local-mv} is given in the Generic Moving Lemma \ref{mvsimple}.

\begin{emp}
Let $(A, \m)$ be a noetherian local ring. 
An ordered sequence ${\bf f}=\{f_1,\dots, f_r\}$ of $r \geq 1$ elements of $\m$ 
is said to be \emph{strictly secant} if 
$f_1$ does not belong to any minimal prime ideal of 
$A$ and if,  for all $i\in \{ 2, \dots, r\}$, $f_i$ does not belong to any minimal prime ideal 
over ${(f_1,\dots, f_{i-1})}$. 
The sequence ${\bf f}$ is called a \emph{strict system of parameters} of $A$ if
$r=\dim A$. In general, the property of being
strictly secant depends on the order of the $f_i$'s. 
\end{emp}

A sequence $\{f_1,\dots, f_r \}$ of $\m$ is called {\it secant }
if $\dim A/(f_1,\dots,f_r) = \dim A - r$ (\cite{Bou8-9}, VIII.26, Definition 1).
The sequence $\{f_1\}$ is secant if and only $f_1$ does not belong to any 
minimal prime ideal $\p$ of $A$ such that $\dim A/\p = \dim A$ (\cite{Bou8-9}, VIII.27, Proposition 3).
By induction, one sees that a strictly secant sequence is secant.

In a noetherian local ring $(A, \m)$  of dimension $d$,
a \emph{system of parameters} of $A$ is a system $f_1,\dots,f_d$ of elements 
of $ \m$
such that the ideal $(f_1,\dots,f_d)$ contains a power of $\m$ or, equivalently, such that 
the $A$-module $A/(f_1,\dots, f_d)$ has finite length. 
Since a strict system of parameters ${\bf f}=\{f_1,\dots, f_d\}$ is a secant sequence, we find that 
 $\dim A/(f_1,\dots,f_d) = \dim A - d=0$, and ${\bf f}$ is a system of parameters.
 
\begin{emp}
\label{height}
Recall that any ideal in $\m$ generated by 
$i$ elements has height at most $i$. 
It follows that a sequence $\{ f_1, \dots, f_r\}$ is strictly secant
if and only if for all $i\le r$, all minimal prime 
ideals over $(f_1,\dots, f_i)$ have height $i$.
\end{emp}
 
\begin{lemma} 
Let $(A, \m)$ be a noetherian catenary equidimensional local ring. 
Then a secant sequence (resp. a system of parameters) 
is strictly secant (resp. a strict system of parameters).  
\end{lemma}
\proof  Let 
$f\in A$ be such that $\dim A/fA=\dim A -1$. Since $A$ is equidimensional,
$f$ does not belong to any minimal prime ideal of $A$ and, thus, the sequence $\{f\}$ is 
strictly secant. The quotient $A/(f)$ is catenary since $A$ is. We claim that 
$A/fA$ is also equidimensional. 
Indeed, let $\q$ be any minimal 
prime ideal of $A$ over $(f)$ and let $\p$ be a minimal prime ideal of $A$ contained
in $\q$. By Krull's principal ideal theorem, $\Ht(\q/\p)=1$. Hence, since $A$ is catenary,
$\Ht(\m/\q)=\Ht(\m/\p)-\Ht(\q/\p)=\dim A -1$ for any minimal prime $\q/(f)$ of $A/(f)$. Therefore, $A/(f)$ is
equidimensional. 

Let now ${\bf f}:=\{f_1,\dots, f_r\}$ be a secant sequence for $A$. If $r=1$, then ${\bf f}$
is strictly secant. If $r>1$, then $\{f_2,\dots, f_r\}$ is secant for $A/(f_1)$ and, by induction, strictly secant
for $A/(f_1)$. Then ${\bf f}$ is strictly secant for $A$. 
\qed

\begin{lemma} \label{strict-secant} 
Let $(A, \m)$ be a noetherian local ring. 
Let $Q$ be a $\m$-primary ideal, 
and let $I$ be a proper ideal of $A$. Let 
${\bf g}:=\{ g_1, \dots, g_r\}$ be a strictly secant sequence of elements
of $(Q+I)/I$. Then there exists a strict system of parameters 
$\{f_1, \dots, f_d\}$ in $Q$ such that for
$i\le r$, the map $A \to A/I$ sends $f_i$ to $g_i$, and   for $i > \dim A/I$,  $f_i\in I\cap Q$. 
\end{lemma}
\proof When $I=(0)$,  the lemma states that 
any strictly secant sequence of elements of $Q$ can always be 
completed into a strict system of parameters contained in $Q$. 
We leave the proof of this fact to the reader.

Assume that $I \neq (0)$. We first show by induction on $r$ that 
${\bf g}$ can be lifted to a strictly secant sequence in $Q$.
Let $f_1' \in Q$ be a lift of $g_1$.
Let $\p_1,\dots,\p_m$ be the minimal prime ideals of $A$.
If $f_1' \notin \p_i$ for all $i =1, \dots, m$, then by definition 
$\{f_1'\}$ is a strictly secant sequence. 

Suppose now that $f_1'$ belongs to some $\p_i$, and 
renumber these minimal primes so that there exists $m_0>0 $ such that
$f_1' \in \p_i$ if and only if $i\le m_0$. 
Let 
$\q_1,\dots, \q_n$ be the minimal prime ideals of $A$ over $I$. 
If $i\le m_0$, then $\q_j\not\subseteq \p_i$ 
(otherwise they would be equal, but
$f'_1\notin \q_j$) and $Q\not\subseteq \p_i$ because 
$\dim A\ge r\ge 1$. So for $i \leq m_0$, 
$Q\cap (\cap_{1\le j\le n}\q_j)\cap (\cap_{k\ge m_0+1} \p_k)$ is not contained in $ \p_i$ and,
therefore, there exists 
$$\alpha\in Q\cap I\cap (\cap_{k\ge m_0+1} \p_k)\setminus 
\cup_{1\le i\le m_0} \p_i.$$
Let $f_1:=f'_1+\alpha\in Q$. Then $\{ f_1\}$ is a strictly secant sequence contained in $Q$, and $f_1$
maps to $g_1$ in $A/I$, as desired. 

Let us assume by induction that we can lift $\{g_1,\dots, g_{r-1}\}$ to a strictly secant sequence
$\{f_1,\dots, f_{r-1}\}$ in $Q$. Let $J:= (f_1,\dots, f_{r-1})$. Apply the case $r=1$ to the ring 
$A':=A/J$, with the ideals $I':=I+J/J$ and 
$Q' := Q+ J/J$: the strictly secant sequence 
$\{ \bar{g_{r}}\}$ in $A'/I'$ lifts to a strictly secant sequence $\{ \bar{f_{r}}\}$ in $Q'$, and we let 
$f_r$ denote a lift in $Q$ of $\bar{f_{r}}$. Then 
$\{f_1,\dots, f_{r-1}, f_r\}$ in $Q$ is the desired strictly secant sequence lifting ${\bf g}$.

Now we complete ${\bf g}$ into a strict system of parameters
in $(Q+I)/I$ and lift it to a strictly secant sequence
$f_1, \dots, f_{n}$ in $Q$ (with $n:=\dim A/I\ge r$). 
It is easy to check that $\m = \sqrt{I \cap Q + (f_1,
\dots, f_n)}$.
Then the image of $I\cap Q$ in $A/(f_1,\dots, f_n)$ contains a power
of the maximal ideal, and we can use the case `$I=(0)$' applied to 
the ring $A/(f_1,\dots, f_n)$ to find that there exist $f_{n+1}, \dots, f_d\in
I\cap Q$ whose images in $A/(f_1,\dots, f_n)$ form a 
strict system of parameters. Then 
$\{f_1,\dots, f_d\}$ is a strict system of parameters of $A$
as desired. 
\qed 

\medskip 

The following theorem, whose proof relies heavily on Proposition \ref{eq-sp},  is the main result of this section.
We use Theorem \ref{local-mv} to prove the Generic Moving Lemma \ref{mvsimple}, which  in turn  is used to prove Theorem \ref{dx-ds}.
 
\begin{theorem}\label{local-mv} 
Let $(A, \m)$ be a noetherian local ring of dimension 
$d\ge 1$, and let $Q$ be a $\m$-primary ideal of $A$. 
Let $F$ be a closed subset of $\Spec A$ with $\dim F<d$. 
Then there exist integral closed subschemes $C_1,\dots, C_n$ of 
$\Spec A$, of dimension $1$, and invertible rational functions
$\varphi_i\in k(C_i)^*$ such that: 
\begin{enumerate}[{\rm (i)}]
\item If $F \neq \emptyset$, then for all $i\le n$, $C_i\cap F = \{\m\}$. 
\item $e(Q)=\sum_{1\le i\le n} \ord_{C_i}(\varphi_i)$.
\end{enumerate}
In particular, $e(Q)[\m]$ is rationally trivial on $\Spec A$. 
\end{theorem}

 We will need the 
following two facts.

\begin{emp} \label{e-irr} Let $(A, \m)$ be a noetherian local ring,
let $Q$ be a $\m$-primary ideal of $A$. 
Let $\p_1,\dots, \p_t$ be the minimal prime ideals of $A$ such that
$\dim A/\p_i=\dim A$. Then 
$$ e(Q, A)= 
\sum_{1\le i\le t} \length(A_{\p_i})e((Q+\p_i)/\p_i, A/\p_i)$$
(\cite{Bou8-9}, VIII, \S7, n$^{\circ}$ 1, Proposition 3). 
In particular, $ e(A)= \sum_{1\le i\le t} \length(A_{\p_i})e(A/\p_i)$.
\end{emp}

\begin{emp} \label{e-f-fp} (Associative Law for Multiplicities.) Let $(A, \m)$ be a noetherian local ring of dimension $d$.
Let $f_1,\dots, f_d$ be a system of parameters. 
Denote  by ${\bf f}$ the sequence $f_1,\dots, f_d$
and by $({\bf f})$ the ideal generated by the $f_i$'s. Fix $s \leq d$. 
Let $\p_1,\dots, \p_t$ be the minimal prime ideals of $A$ over $(f_1,\dots,f_s)$
such that $\Ht(\p_i)=s$ and
$\dim A/\p_i=d-s$. (When $s=0$, we interpret $(f_1,\dots,f_s)$ to be the ideal $(0)$.)
Then
$$ e(({\bf f}),A)= \sum_{i=1}^t e((f_1,\dots,f_s)A_{\p_i}, A_{\p_i}) e((f_{s+1},\dots,f_d)+\p_i/\p_i, A/\p_i)$$
(see \cite{Lec}, or \cite{Nor}, Theorem 18, page 342, or \cite{HIO} (1.8), or \cite{Mat2}, exer. 14.6, page 115).
\end{emp}

\begin{emp} \label {proof.local-mv}
{\it Proof of \ref{local-mv}.} 
It suffices to prove the theorem when $F\ne \emptyset$.
Let $I$ be a proper ideal of $A$ such that $F=V(I)$. By hypothesis,
$r:=\dim F<d:=\dim A$. 

Proposition \ref{eq-sp} (1) shows that the multiplicity $e(Q)$
can be expressed in terms of the multiplicity of ideals generated 
by (strict) systems of parameters.  Thus, we are reduced to proving the theorem in the case $Q$ is 
generated by a system of parameters.

Proposition \ref{eq-sp} (2) shows that it is enough to  consider
the case where $Q$ is generated by a strict system of parameters 
${\bf f}=\{ f_1, \dots, f_d\}$
of $A$ such that, when $r\ge 1$, 
the image of $\{f_1, \dots, f_r \}$  in $A/I$ is a strict 
system of parameters of $A/I$. 
Let $\p_1, \dots, \p_t$ be the minimal prime
ideals of $A$ over $(f_1,\dots,f_{d-1})$ of height $d-1$ with $\dim A/\p_i=1$. 
The formula in \ref{e-f-fp} gives
$$e(Q,A)=
\sum_{1\le i \le t}  d_i e(f_dA, A/\p_i), 
$$ 
with $d_i:= e((f_1,\dots, f_{d-1})A_{\p_i},A_{\p_i})$.
As $A/\p_i$ is integral, $e(f_dA, A/\p_i)=\ord_{A/\p_i}(f_d)$
where $\ord_{A/\p_i}(f_d)$ denotes the 
order of the image of $f_d$ in $A/\p_i$. 
Let $C_i:=V(\p_i)$, $1\le i\le t$,
and let $\varphi_i:=f_d^{d_i}|_{C_i}$. 
Then 
$$e(Q,A)=\sum_{1\le i\le t}\ord_{C_i}(\varphi_i),$$
as desired. When $\dim F=r>0$,  $\dim A/(I,f_1,\dots,f_r) =0$ by construction.
Thus, $\dim (C_i\cap F) = \dim V(I+\p_i) = 0$ for all $i=1,\dots, t$.
That $e(Q)[\m]$ is rationally equivalent to $0$ follows
immediately from the definition recalled at the beginning of  section \ref{section2}.
\qed 
\end{emp}

\medskip
Part (3) of our next proposition is a slight strengthening of a well-known theorem
(\cite{B-H}, Corollary 4.5.10, or 
\cite{Z-S}, VIII, \S 10, Theorem 22 when $M=A$). 
Parts (1) and (2) are used in \ref{local-mv}, and Part (4) will be used in \ref{gamma-1}.

\begin{proposition}\label{eq-sp} Let $(A, \m)$ be a noetherian
local ring and 
let $M$ be a non-zero finitely generated $A$-module.
Let $Q$ be an
$\m$-primary ideal of $A$.
\begin{enumerate}[{\rm (1)}]
\item There exist finitely many 
strict systems of parameters ${\bf f}^{\alpha}, \alpha \in {\mathcal A},$
contained in $Q$, such that $e(Q, M)= \sum_{\alpha \in {\mathcal A}} \pm e(({\bf f}^{\alpha}), M)$.
\item Fix a proper ideal $I$ of $A$ with $\Ann_A(M)\subseteq I$. Let
$r:=\dim V(I)$. When $r>0$, we can choose the strict system of parameters 
$\{ f_1^{\alpha}, \dots, f_{\dim A}^{\alpha}\}$ in {\rm (1)} such that for each $\alpha \in {\mathcal A}$,
the image of 
the sequence $\{ f_1^{\alpha}, \dots, f_r^{\alpha}\}$ in $A/I$ 
forms a strict system of parameters in $A/I$.
\item If $A/\m$ is infinite, then $e(Q, M)=e(({\bf f}), M)$ for
some strict system of parameters ${\bf f}$ contained in $Q$, satisfying, when applicable, 
the property in {\rm (2)}. 

\item
Let $\p$ be a prime ideal of $A$ with $ \Ht(\p) \ge 1$. Let $Q_0$ be a $\p A_\p$-primary ideal of $A_\p$.
Then there exist finitely many strictly secant sequences 
${\bf f}^{\alpha}$ in $\m$, $\alpha \in {\mathcal A}$, whose
images in $A_\p$ are contained in $Q_0$ and such that
$e(Q_0, A_\p)=
\sum_{\alpha \in {\mathcal A}} \pm e(({\bf f}^{\alpha}), A_\p)$.
\end{enumerate}
\end{proposition}

\proof (1) and (2). Let $J:=\Ann_A(M)$.  
By definition, $\dim M:=\dim A/J$.
We can consider $M$ as an $A/J$-module.
We have $e_A(Q, M)=e_{A/J}((Q+J)/J, M)$ because $JM=0$. Suppose
that we can construct strict systems of parameters 
${\bf g}^{\alpha}$ in $(Q+J)/J$ such that 
$e_{A/J}((Q+J)/J, M)$ is a combination of the $e_{A/J}(({\bf g}^{\alpha}), M)$'s
and such that the image of $\{g_1^{\alpha},\dots, g_r^{\alpha}\}$ in
$A/I$ is a strict system of parameters for all $\alpha$. 
By Lemma \ref{strict-secant}, 
we can lift each ${\bf g}^{\alpha}$
and complete it into a strict system of parameters 
${\bf f}^{\alpha}=\{ f_1^{\alpha}, \dots, f_{\dim A}^{\alpha}\}  $ in $Q$ with $f_i^{\alpha}\in J$ for all
$i > \dim M$. Then 
$e_{A/J}(({\bf g}^{\alpha}), M)=e_A(({\bf f}^{\alpha}), M)$, and (1) and (2) hold. 
Therefore, it remains to prove (1) and (2) when  $\dim M=\dim A$. 

Assume that $d:=\dim M= \dim A$. We proceed  by induction on $d$. 
If $d=0$, then a strict system of parameters ${\bf f}$ for $A$ 
is empty, and we set ${\bf f}=(0)$.
Then $e(Q, M)=\length(M)= e((0), M)$, and no additional condition
is required in (2) since $\dim A/I=0$.

Assume that $d \geq 1$. Let $P$ be
 the set consisting of the
minimal prime ideals of $A$, 
 of the associated primes of $M$, and if $\dim A/I>0$,  of  the minimal prime ideals of $A$ over $I$. 
 Our hypotheses imply that if $\m \notin {\rm Ass}(M)$, then $\m \notin P$.
Suppose that $\m \in {\rm Ass}(M)$. Then there exists an exact sequence
$(0) \to M' \to M \to M'' \to (0)$ with $M' $ isomorphic to $A/\m$.
Since $\dim(M)>0$, we find that for any $\m$-primary ideal $Q_0$ of $A$, 
$e(Q_0,M)=e(Q_0,M'')$ (\cite{Bou8-9}, VIII.46, Prop. 5,  
or use the proof\footnote{It should be noted that the definition
of $e(Q,M)$ taken in \cite{Mat2}, page 107, (or in \cite{SH}, 11.1.5) is different from the definition 
taken at the beginning of this section and in \cite{Bou8-9},
VIII.72. Thus the statement of \cite{Mat2}, 14.6,  cannot be applied directly.} of \cite{Mat2}, 14.6).
It follows that it suffices to prove our statement for modules $M$ with $\m \notin {\rm Ass}(M)$.
In particular, we can assume that $\m \notin P$. 

We apply Lemma \ref{induction} with this set $P$ of primes.
Pick some $s \geq s_0$ and $x_s \in Q^s$, $x_{s+1} \in Q^{s+1}$ as in \ref{induction}.  
By our choice of $P$, both $\{ x_s\}$  and $\{ x_{s+1}\}$ are
strictly secant sequences in $Q$  (\ref{induction}(c)).
Suppose $d=1$. 
Then $e(Q, M)=e((x_{s+1}), M)-e((x_s),M)$ (\ref{induction}(b)). 
When $r=1$, our choice of $P$ shows that the images of $\{ x_s\}$  and $\{ x_{s+1}\}$ in 
$A/I$ are again strict systems of parameters in $A/I$. 
Hence, (1) and (2) are true when $d=1$. 

Suppose now that $d\ge 2$ and that (1) and (2) hold for $d-1$. 
Lemma \ref{induction}(b)
shows that $$e(Q, M)=e(Q, M/x_{s+1}M)-e(Q, M/x_sM).$$ 
Let $B:=A/x_sA$. Then $\dim B = d-1$ and $M/x_sM$ is a $B$-module of dimension
$d-1$ (\ref{induction}(a)). Let $IB$ denote the ideal $(I+x_sA)/x_sA$ of $B$.
By our choice of primes in $P$, we find that  $\dim V(IB)=\max\{ 0, r-1\}$.
We also have $e_A(Q, M/x_sM)=e_B(QB, M/x_sM)$. By induction hypothesis, 
there exist finitely many strict systems of parameters ${\bf g}^{\alpha}$ in $QB$ with 
$$e(QB, M/x_sM)= \sum \pm
e(({\bf g}^{\alpha}), M/x_sM)$$
and 
such that, if $\dim V(IB) \geq 1$, 
the image  of $g_1^{\alpha}, \dots, g_{r-1}^{\alpha}$
in $B/IB$ is a strict system of parameters. Let 
${\bf f}^{\alpha}$ be any lifting of ${\bf g}^{\alpha}$ in $Q$. 
Then $\{ x_s, {\bf f}^{\alpha}\}$ is a 
strict system of parameters in $Q$, and the image in $A/I$  of the 
first $r$  elements of $\{ x_s, {\bf f}^{\alpha}\}$ is a strict system of parameters in $A/I$. 
By our choice of primes in $P$, $x_s$ is not a zero divisor in $M$.
It follows from \cite{Mat2}, 14.11  that\footnote{
Due to the two different definitions of $e(Q,M)$, the proof of \cite{Mat2}, 14.11, needs to be slightly 
adjusted to prove the second equality.}
$$e(({\bf g}^{\alpha}), M/x_sM) = e(({\bf f}^{\alpha}), M/x_sM)
= e((x_s, {\bf f}^{\alpha}), M).
$$
Repeating the same argument with $x_{s+1}$ instead of $x_s$ allows us to express the multiplicity 
$e(Q, M/x_{s+1}M)$ in a similar way. Since $e(Q, M)=e(Q, M/x_{s+1}M)-e(Q, M/x_sM)$,
(1) and (2) follow. 

(3) When $A/\m$ is infinite, 
instead of 
choosing  two elements $x_s$ and $x_{s+1}$ as we did above,
Part (3) is proved by modifying the proof of (1) and (2), using 
only the element $x_1$ whose existence is asserted
in \ref{induction}(d), with $e(Q,M)= e(Q, M/x_1M)$.
We leave the details to the reader.

(4) 
Let $R$ be the preimage of $Q_0$ in $A$ under the natural map $A \to A_\p$.
Clearly, $R\subseteq \p$.
Consider the graded rings
$$\Gr(A):=\oplus_{n\ge 0} R^n/R^{n+1}  \mbox{ and }
\Gr(A_\p):=\oplus_{n\ge 0} Q_0^n/Q_0^{n+1},$$ 
and the natural homomorphism of graded rings
$\rho : \Gr(A)\to \Gr(A_\p)$.
Let $\q'_1, \dots, \q'_k$ be the
associated prime ideals of $\Gr(A_\p)$ not containing $\Gr_+(A_\p)$.
Then $\{ \rho^{-1}(\q'_j), j=1,\dots, k\}$ is a set of homogeneous prime ideals of
$\Gr(A)$ not containing $\Gr_+(A)$.
Let $\q_1,\dots, \q_m$ be the minimal prime ideals of $A$.
We are going to associate below to each $\q_i$ a 
homogeneous prime ideal $\widetilde{\q_i}$ 
of
$\Gr(A)$ which does not contain $\Gr_+(A)$. We will then apply 
the Prime Avoidance Lemma \ref{avoid} 
to the ideal $I:= \Gr_+(A)$ in the ring $\Gr(A)$ with the set of homogeneous ideals 
$ \{\rho^{-1}(\q'_j), j=1,\dots, k\} \cup \{ \widetilde{\q_i}, i=1,\dots, m\}$. 

For any prime ideal $\q$ of $A$, let
$\q^{hom}:=\oplus_{n\ge 0} (\q\cap R^n)/(\q\cap R^{n+1})$, which we view in a natural way as an ideal in $\Gr(A)$.
This ideal is homogeneous, but may not be prime. We claim that when $\p  \not \subseteq \q$, then $\sqrt{\q^{hom}}$ does not contain $\Gr_+(A)$.
Indeed, if $\sqrt{\q^{hom}}$  contains $\Gr_+(A)$, then $ \q^{hom}$  contains a power of $\Gr_+(A)$. So for some $n >0$, 
$R^n = R^{n+1} + (\q \cap R^n)$. Hence, modulo $\q$, $\overline{R}^n = \overline{R}^{n+1}$. Since $\cap_{m \geq n} \overline{R}^m = (0) = \overline{R}^n$, we find that $R^n \subseteq \q$, and passing to radicals, we find that $\p \subseteq \q$, a contradiction.
Now by standard results, $\sqrt{\q^{hom}}$ is homogeneous, and it is the intersection of all the prime ideals minimal above it, which are also homogeneous.
For each $i=1, \dots, m$, we let $\widetilde{\q_i}$ denote one of the minimal prime ideals of  $\sqrt{\q_i^{hom}}$ that does not contain $\Gr_+(A)$.

We conclude from  Lemma \ref{avoid}, 
applied to the ideal $I:= \Gr_+(A)$ in the ring $\Gr(A)$ with the set of homogeneous ideals 
$ \{\rho^{-1}(\q'_j), j=1,\dots, k\} \cup \{ \widetilde{\q_i}, i=1,\dots, m\}$,
that there exists $s_0 \geq 0$ such that, for all
$s\geq s_0$, there exists $x_s\in R^s$ whose class in $R^s/R^{s+1}$
does not belong to $\cup_j \rho^{-1}(\q'_j)\cup (\cup_i \widetilde{\q_i})$.
In particular, $x_s\notin\q_i$ for all $i =1, \dots, m$,
which by definition implies that $\{ x_s \}$ is strictly secant in $A$.
Moreover, let $\xi$ denote  the image of $x_s$ in $Q_0^s/Q_0^{s+1}$. Then, by construction,
$\xi \notin \q'_j$, for all $j=1,\dots, k$.

As in \ref{induction}, we use this latter fact to be able to apply \cite{Bou8-9}, VIII.79, 
Prop.\ 9 and Lemma 3, to verify a hypothesis needed 
in \cite{Bou8-9}, VIII.77, Prop.\ 8.b). 
It follows immediately from VIII.77, Prop.\ 8.b), (i) and (ii), because $\Ht(\p) \geq 1$,
 that 
$$se(Q_0, A_{\p})=\left\{ \begin{array}{l} e(Q_0, A_{\p}/x_sA_\p) \text{ if } \dim A_\p\ge 2, \\
e(( x_s), A_{\p}) \text{ if } \dim A_\p = 1.
\end{array} \right. $$

We can now prove (4) when $\Ht(\p)=1$. Fix $s \geq s_0$ as above, and consider $x_s$ and $x_{s+1}$. We find that 
$$e(Q_0, A_{\p})= e(( x_{s+1}), A_{\p}) - e(( x_s), A_{\p}),$$
with both $\{ x_s \}$ and $\{ x_{s+1} \}$ strictly secant in $A$, as desired.
When $\Ht(\p) \geq 2$, we proceed by induction on $\Ht(\p)$.
Since $\{x_s\}$ is strictly secant in $A$, it is secant in $A_\p$, so that
$\dim A_\p/(x_s) = \Ht(\p/(x_s)) = \Ht(\p) -1$. Let $B:= A/(x_s)$. We apply the induction 
hypothesis to the ring $B$, prime ideal $\p B$ and $\m B$-primary ideal $Q_0B$. The details 
are left to the reader.
\qed

\begin{lemma}\label{induction} 
Let $(A, \m)$ be a noetherian local ring. Let $Q$ be a proper ideal of $A$. 
Let $M$ be a finitely generated $A$-module with 
$\dim M\ge 1$ and $M/QM$ of finite length.  
Let $P$ be a finite set of prime ideals of $A$ with $\m \notin P$. 
Then there exists $s_0\ge 1$ such that  for all $s\ge s_0$, there exists
$x_s\in Q^s$ with the following properties: 
\begin{enumerate}[\rm (a)] 
\item $\dim (M/x_s M)=\dim M -1$. 
\item  $se(Q, M)= \left\{ \begin{array}{l} e(Q, M/x_sM)   \text{ if } \dim M\ge 2, \text{ and}  \\
e(x_sA, M) \text{ if } \dim M=1.
\end{array} \right.$
\item  $x_s\notin\p$ for every $\p\in P$.
\item If $A/\m$ is infinite, then there exists an element $x_1 \in Q \setminus Q^2$ with the above properties.
\end{enumerate}
\end{lemma} 

\proof Consider the graded ring 
$\Gr(A):=\oplus_{n\ge 0} Q^n/Q^{n+1}$. For any $A$-module $N$, 
the group  $\Gr(N) := \oplus_{n\ge 0} Q^nN/Q^{n+1}N$ is a $\Gr(A)$-module in a natural way. 
Consider then the $\Gr(A)$-graded
module $L:=\Gr(M)\oplus (\oplus_{\p\in P}\Gr(A/\p))$. 
Let $\p_1, ..., \p_r$ be the 
associated prime ideals of  $L$ 
not containing $\Gr_+(A)$. These are homogeneous ideals of $\Gr(A)$ (\cite{Mat}, (10.B)).  
By Lemma \ref{avoid} 
applied to $I=\Gr_+(A)$, 
there exists $s_0\ge 1$ such that 
for all $s\ge s_0$, $\Gr_s(A)$ is not contained in $\cup_i \p_i$. 

For $s \geq s_0$, let $x_s\in Q^s$ whose class $\xi$ in $\Gr_s(A)$ does not belong
to $\cup_i \p_i$. Let $\varphi : L \to L$ 
be the multiplication-by-$\xi$ map. 
By \cite{Bou8-9}, VIII.79, 
Prop.\ 9 and Lemma 3, 
$\Ker\varphi$ has finite length over $\Gr(A)$ and, hence, over $A/Q$.
This is one of the hypotheses needed to now apply \cite{Bou8-9}, VIII.77, Prop.\ 8.b).
Moreover, our hypothesis that $\m \notin P$ shows that $ \dim A/\p \geq 1 $ for all $\p \in P$.
Since we also assumed that $\dim M \geq 1$, we can  now
use \cite{Bou8-9}, VIII.77, Prop.\ 8.b), to obtain that 
$\dim M/x_sM = \dim M-1$ (proving (a)), and for each $\p \in P$,
$\dim (A/(x_s, \p))=\dim (A/\p) -1$. In particular, it follows that $x_s\notin \p$, proving (c).  
Part (b) follows immediately from VIII.77, Prop.\ 8.b), (i) and (ii).
Remarque 4) on \cite{Bou8-9}, VIII.79, shows that 
when $A/\m$ is infinite, we can apply the above proof to an element $x_1\in Q \setminus Q^2$.
\qed  

\medskip
The following Prime Avoidance Lemma for graded rings is needed in the proofs of \ref{eq-sp} and \ref{induction}.
(For related statements, see \cite{SH}, Theorem A.1.2., or \cite{Bou}, III, 1.4, Prop. 8, page 161.)
\begin{lemma}\label{avoid} 
Let $B=\oplus_{s\ge 0} B_s$ be a graded ring.
Let $\p_1,\dots,\p_r$ be homogeneous prime ideals of $B$ not containing $B_1$.
Let $I=\oplus_{s\ge 0} I_s$ be a homogeneous ideal 
such that $I\not\subseteq \p_i$ for all $i\le r$. 
Then there exists an integer $s_0 \geq 1$ such
that for all $s\ge s_0$, 
$I_s\not\subseteq \cup_{1\le i\le r} \p_i$.    
\end{lemma}

\proof We proceed by induction on $r$. If $r=1$, choose
$t\in B_1\setminus \p_1$ and a homogeneous element
$\alpha\in I\setminus \p_1$, say of degree $s_0$.
Then $t^{s-s_0}\alpha\in I_s \setminus \p_1$ for all $s\ge s_0$, as desired.
Let $r\ge 2$ and suppose that the lemma is true for $r-1$. We can suppose
that $\p_i$ is not contained in $\p_r$ for all $i\ne r$, so that 
$I\p_1\cdots\p_{r-1} \not\subseteq \p_r$. Similarly, 
we can suppose that $\p_r$ is not contained in $\p_i$ for all 
$i\ne r$, so that 
$I\p_r\not\subseteq (\p_1\cup \ldots\cup\p_{r-1})$. 
Hence, we can apply the case $r=1$ and the induction hypothesis 
to obtain that there exists $s_0$ such that for all $s\ge s_0$, 
there are homogeneous elements 
$f_s\in I\p_1\cdots\p_{r-1} \setminus \p_r$ and 
$g_s\in I\p_r\setminus (\p_1\cup \ldots\cup\p_{r-1})$
of degree $s$. It is easy to check that 
$f_s+g_s\in I_s\setminus \cup_{1\le i\le r}\p_i$, as desired. 
\qed

\begin{example}  
Let $(A, \m)$ be a noetherian local ring, and let $Q$ be an
$\m$-primary ideal of $A$. Proposition \ref{eq-sp} (3), 
states that  when $A/\m$ is infinite, there exists a system of parameters ${\bf f}$ contained in $Q$ such that  $e(Q, A)=e(({\bf f}), A)$).
We show below that when  $A/\m$ is finite, such
a system of parameters ${\bf f}$ contained in $Q$ may not exist.

Consider the ring $A:= {\mathbb F}_2[[x,y]]/(xy(x+y))$, with $\m = (x,y)$. It is clear that   $e(\m,A)=3$.
Since this ring has dimension $1$, the multiplicity $e((f), A)$ is equal to the length of $A/(f)$.
It is shown in \cite{HS}, 3.2, that there exists no regular element $f \in \m$ such that $\m^3 \subseteq (f)$.
This implies that $A/(f)$ cannot have length $3$. 

Note that in the ring $B:= {\mathbb F}_4[[x,y]]/(xy(x+y))$,
the element $z:=x-ty$, with $t \in {\mathbb F}_4 \setminus {\mathbb F}_2$, is such that $B/(z)$ has length $3$.
\end{example}
\end{section}

\begin{section}{Two local invariants} \label{mv-mul}

Let $(A, \m)$  be a noetherian local ring 
of
positive dimension. 
We introduce in this section two new invariants of $A$, $\gamma(A)$ in \ref{gammaA}, 
and $n(A)$ in \ref{def-na}. We show in Theorem \ref{thm.na} below that 
$n(A) = \gamma(A)$  when $A$ satisfies some natural hypotheses such as being reduced, excellent, and equidimensional. 
In \ref{cor.compute-n},   we will relate $n(A)$
to a resolution of singularities 
$Y \to \Spec A$ when $A$ is universally catenary.

\begin{emp} \label{gammaA} Let $(A, \m)$ be a noetherian local ring. 
Consider the set 
$\mathcal E$ of all Hilbert-Samuel multiplicities $e(Q,A)$, for all $\m$-primary ideals $Q$ of $A$.
Let $\gamma(A)$ denote {\it the greatest common divisor of the elements of $\mathcal E$}.
Clearly, $\gamma(A) $ divides $ e(\m)$.

Proposition \ref{eq-sp} (1) shows that the greatest common divisor of the integers $e(Q)$, 
taken only over the subset of all
ideals $Q$ generated by strict systems of parameters, is equal to $\gamma(A)$.
Theorem \ref{local-mv} implies that $\gamma(A)[\m]$ is rationally equivalent to zero in $\Spec A$.
It is obvious that $\gamma(A)=1$ when $A$ is regular. When $\dim A=0$, it follows from the definitions
that $\gamma(A) = \ell_A(A)$.

Let $d:=\dim A \geq 1$, and let $U \subseteq \Spec A$ be any dense open subset. 
Fix $d' \in [0,d-1]$, and let 
$$ U{(d')}:= \{ \p \in U \mid  \Ht(\p)= d', \  \dim(A/\p) =d-d'\}. $$
Let $$g(U):= \gcd_{\p\in U{(d')}} \{ \gamma(A_\p)\gamma(A/\p)\}.$$
If $V \subseteq U$ is another dense open subset,
then $g(U) $ divides $ g(V)$. 
When $d'=0$, it follows immediately from \ref{e-irr} and the definitions
that $g(U) $ divides $ \gamma(A)$. We generalize this statement in our next proposition; 
the case $d'=d-1$ will be used in the proof of Theorem \ref{thm.na}.
\end{emp}

\begin{proposition} \label{gamma-1}  
Let $(A, \m)$ be a noetherian local ring of 
dimension $d\ge 1$, and let $U \subseteq \Spec A$ be any dense open subset. 
Fix  $d' \in [0,d-1]$. 
Then $U{(d')}\ne \emptyset$ and 
$$\gamma(A)=\gcd_{\p\in U{(d')}} \{ \gamma(A/\p)\gamma(A_\p)\}.$$ 
\end{proposition}

\proof Let us first prove that $g(U) $ divides $ \gamma(A)$. 
As mentioned above, the case $d'=0$ follows from \ref{e-irr}, so we now assume 
that $d'\ge 1$. It suffices to prove the divisibility for any open dense $V \subseteq U$, so 
shrinking $U$ if necessary, 
we can suppose that the complement $V(I)$ of 
$U$ in $\Spec A$ has dimension $d-1$. Moreover, $\Ht(I)>0$ since $U$ is dense.
 
Let ${\bf f}=f_1,\dots, f_d$ be any strict system of 
parameters of $A$ such that $f_1,\dots, f_{d-1}$ induces a strict system of 
parameters of $A/I$ (\ref{strict-secant}). In particular, the image of 
$\{ f_1,\dots, f_{d'} \}$ in $A/I$ is strictly secant. 
Therefore, any minimal prime ideal of $A$ over 
$(f_1,\dots, f_{d'})$ has height $d'$, and any minimal prime ideal 
of $A/I$ over the ideal $I+ (f_1,\dots, f_{d'})/I$ has height $d'$ (\ref{height}).
Let $\p$ be a minimal prime ideal of $A$ over $(f_1,\dots, f_{d'})$.
Then $\p \in U$, since otherwise, $I \subseteq \p$ and
$d'=\Ht(\p)\ge \Ht(\p/I) + \Ht(I) =d'+ \Ht(I)> d'$. 
Since $\dim A/(f_1,\dots, f_{d'}) = d-d'$, there exists such a $\p$ with 
$\dim A/\p=d-d'$, and then $\p\in U{(d')}$. 
Using \ref{e-f-fp}, we find that
$g(U) $ divides $ e(({\bf f}), A)$.
Hence, Proposition \ref{eq-sp} (2) implies that
$g(U) $ divides $ \gamma(A)$.

Let us now prove that $\gamma(A) $ divides $ g(U)$. It suffices to prove this divisibility
when $U=\Spec A$. We 
start with the case $d'=d-1$. Fix $\p\in U{(d-1)}$, so that $\Ht(\p)= d-1$,
and $\dim A/\p =1$. We need to show that $\gamma(A) $ divides $ \gamma(A_\p)\gamma(A/\p)$.
Proposition \ref{eq-sp} (1) shows that $\gamma(A/\p)$ is equal to the greatest common 
divisor of the integers $\ord_{A/\p}(\phi)$, where $\phi\in\Frac(A/\p)^*$.
When $d=1$, then $\gamma(A_{\p}) =\length(A_{\p})$. (Recall the convention in \ref{e-f-fp}
that if $d=1$, then 
$e((f_1,\dots, f_{d-1}), A_{\p})=\length(A_{\p})$.)
Then when $d\geq 1$, Proposition \ref{eq-sp}(4)
shows that $\gamma(A_\p)$ can be computed as the greatest common divisor
of the integers $e((f_1,\dots, f_{d-1}), A_{\p})$, 
where $f_1,\dots, f_{d-1}$ is a strictly secant sequence in $\p$.
Therefore, it is enough to 
show that $\gamma(A)$ divides $\ord_{A/\p}(\phi)e(({\bf f}'), A_\p)$ 
for all $\phi\in\Frac(A/\p)^*$ and for all strictly secant sequences
${\bf f}'=f_1,\dots, f_{d-1}$ contained in $\p$.

Fix $\phi\in\Frac(A/\p)^*$ and  
a strictly secant sequence
${\bf f}'=f_1,\dots, f_{d-1}$ contained in $\p$. 
By construction, $\p_1:=\p$ is a minimal prime ideal of $A$ over $({\bf f}')$.
Let $ \p_2, \dots, \p_m$ denote  the other minimal primes over $({\bf f}')$. 
Use the isomorphism $\Frac(A/\sqrt{({\bf f}')}) \to \oplus_i \Frac(A/\p_i)$
to find an invertible 
rational function 
$\phi'\in\Frac(A/\sqrt{({\bf f}')})$ which restricts to $\phi $ in $\Frac(A/\p)$ and 
to  $1$ in $\Frac(A/\p_i)$
for all $i\ge 2$. Let $a, b\in\m\setminus \cup_{i\ge 1}\p_i$ 
be such that $a/b$ maps to $\phi'$ in $\Frac(A/\sqrt{({\bf f}')})$. 
The sequences $\{{\bf f}', a\}$ and $\{{\bf f}', b\}$ are 
strict systems of parameters in $\m$. Then 
\ref{e-f-fp} gives
$$ e(({\bf f'}, a), A)=\sum
e(({\bf f'}), A_{\p_i})\ord_{A/\p_i}(a),$$
where the sum is over the indices $i$ such that $
\dim A/\p_i =1$.
We proceed similarly for   $b$ to find that 
$$e(({\bf f}', a), A)-e(({\bf f}', b), A)
=\ord_{A/\p}(\phi)e(({\bf f'}), A_\p),$$ 
because $a = b$ in $A/\p_i$ if $i\ge 2$. Hence, 
$\gamma(A) $ divides $ \ord_{A/\p}(\phi)e(({\bf f'}), A_\p)$ and the
proposition is proved when $d'=d-1$.

The proposition is now  true when $d=1$. 
Suppose that $d \geq 2$, and 
 proceed by induction on $d$. 
Fix $d'\le d-2$. We showed above that $(\Spec A){(d')} \neq \emptyset$, so let 
$\p\in (\Spec A){(d')}$, with $\Ht(\p)=d'$, and $\dim A/\p = d-d'$. 
Since $(\Spec A/\p){(d-d'-1)} \neq \emptyset$, let 
 $\q\supseteq \p$ be any prime ideal of $A$ such that $\Ht(\q/\p)=d-d'-1$
and $\dim A/\q=1$. 
It follows that  $\Ht(\q)=d-1$, so that we obtain from the above considerations 
that 
$$\gamma(A) \text{ divides }   \gamma(A/\q)\gamma(A_{\q}).$$
Since $\dim A_{\q} =d-1$, we can apply the induction hypothesis 
to $A_{\q} $; consider the ideal $\p A_{\q}$, with $\Ht(\p A_{\q})= d'$ and 
$\dim A_\q/\p A_{\q}) =(\dim A_{\q}) - d'$. It follows that
$$\gamma(A_{\q}) \text{ divides } \gamma(A_{\q}/\p A_{\q}) \gamma((A_{\q})_{\p A_{\q}}).
$$
Let $B:=A/\p$ and $\bar{\q}:=\q B$. With this notation, the above two divisibility conditions give
$$\gamma(A) \text{ divides } \gamma(B/{\bar{\q}})\gamma(B_{\bar{\q}})\gamma(A_{\p}).$$
Varying $\bar{\q}$ in $(\Spec B){(\dim B -1)}$,
 and using the case $d'=d-1$ 
established above, we get that $\gamma(B) = \gcd_{\bar{\q}} \gamma(B/\bar{\q}) \gamma(B_{\bar{\q}})$,
so that  $\gamma(A) $ divides $ \gamma(B)\gamma(A_{\p})$.
Since $B=A/\p$, the result follows.
\qed

\begin{emp} 
Recall that the generic point of an irreducible 
component of a closed subset $F$ of a scheme $X$ is called a {\it maximal point} of $F$.
For convenience in the statement of our next definition, 
we denote by $\max(F)$ the set of the maximal points of $F$.

Let $(A, \m)$ be a noetherian local ring of
positive dimension.
Let $U$ 
be any dense open 
subset of $\Spec A$. 
Denote by $x_0$ the unique closed point of $\Spec A$.
Let $\mathcal N_U$ denote the set of all 
integers $n$ occurring as the
order 
$\ord_{x_0}(f)$
of some rational function $f\in \mathcal K_C^*(C)$  
on any reduced closed curve $C$ in $\Spec A$ with 
$\max(C)\subseteq U$.

{\it We claim that $\mathcal N_U$ is in fact an ideal in ${\mathbb Z}$.}
Indeed, 
suppose that for a given open subset $U$ in $\Spec A$, 
there exist reduced closed one-dimensional subschemes 
$C$  and $C'$ in $\Spec A$, and $f \in \mathcal K_C^*(C)$, $f' \in \mathcal K_{C'}^*(C')$, 
such that $\max(C)\subseteq U$, $\max(C')\subseteq U$, and $n=\ord_{x_0}(f)$, $n'=\ord_{x_0}(f')$.
Let us use \ref{fact-1}, (1) and (2), to build  functions  on $C \cup C'$. 
First extend $f$ to a function $F$ on $C \cup C'$, with $F=1$ on every component of $C'$ 
which is not a component of $C$. Similarly, extend $g$ to $G$ on $C \cup C'$, with $G=1$ on every component of $C$ 
which is not a component of $C'$.
Then $F^aG^b$ is such that $\ord_{x_0}(F^aG^b)=an+bn' \in \mathcal N_U$.

Denote by $n(U, A)$, or $n(U, \Spec A)$, 
{\it the greatest common divisor of the positive elements of $\mathcal N_U$}.
Clearly, if $V \subseteq U$ is dense and open in $\Spec A$,
then $n(U,A) $ divides $ n(V,A)$. 

\end{emp}
\begin{emp} \label{def-na}
Let $A$ be a noetherian local ring of
positive dimension as above.
Consider the ideal $\mathcal N := \cap_U \mathcal N_U$,
where $U$ runs through all dense open subsets of  $\Spec A$.
Theorem \ref{local-mv} immediately implies that $\mathcal N \neq (0)$. 
Let  $n(A)$ denote {\it the greatest common divisor of the positive elements of $\mathcal N$}.
The integer $n(A)$ is the positive generator of the ideal $\mathcal N$. 
When $\dim A=0$, we let $n(A):=1$.

By definition, $n(U,A) $ divides $n(A)$ for any dense open subset $U$ of $\Spec A$.
In fact, $n(A)$ is the smallest positive integer $n$ such that, for every dense open set $U$ of $\Spec A$,
there exists a reduced closed curve $C$ in $\Spec A$ with 
$\max(C)\subseteq U$ and a rational function $f\in \mathcal K_C^*(C)$ such that $n=\ord_{x_0}(f)$.
In view of Theorems \ref{local-mv} and \ref{mv}, we call $n(A)$ the {\it moving multiplicity} of $A$.

\begin{lemma} \label{lem-n(U,A)} Let $A$ be a noetherian local ring of
positive dimension.
There exists a dense open subset $U_0$ of $\Spec A$ such that $n(U_0, A)=n(A)$.
In particular, for all dense open subsets $V$ contained in $U_0$,
 $n(V, A)=n(A)$. 
\end{lemma}

\proof Since the ring ${\mathbb Z}/{\mathcal N}$ is Artinian, 
we can find finitely many dense open sets $U_i$ such that 
$\mathcal N = \cap_i \mathcal N_{U_i}$.
We let $U_0:= \cap_i U_i$. 
\qed 
\end{emp}

Theorem \ref{local-mv} motivated our definition of $n(A)$.
It follows from 
this theorem 
that   
$n(A) $ divides $ e(Q)$ for all $\m$-primary ideals $Q$ of $A$. Hence, 
$$n(A) \text{ divides } \gamma(A).$$
In particular, if $A$ is regular, then $n(A)=1$, but the converse is false. 
We now show that under some natural hypotheses, such as $A$ being reduced, excellent, and equidimensional of positive dimension,
then $n(A)=\gamma(A)$ (see also \ref{rem.add}). 

\begin{theorem} \label{thm.na} Let $A$ be a noetherian equidimensional local ring of positive dimension. 
Suppose that $A$ is catenary and the regular locus of $\Spec A$ 
contains a dense open subset of $\Spec A$. 
Then $n(A)=\gamma(A)$.
\end{theorem}

\proof It suffices to show that $\gamma(A) $ divides $ n(A)$. Let $U$ be any dense open subset of the regular 
locus of $\Spec A$.
By hypothesis, there exist a reduced curve  $C$ in  $\Spec A$ with
$\max(C)\subseteq U$ and $f\in \mathcal K_C^*(C)$ such that
$n(A)=\ord(f)$. Let $C_1,\dots, C_n$ be the irreducible components
of $C$. Each $C_i$ corresponds to a prime ideal $\p_i$ of $A$.
Since the generic point of $C_i$ belongs
to $U$, $A_{\p_i}$ is regular, so that $\gamma(A_{\p_i})=1$.
Because $A$ is catenary and equidimensional, we find that $\Ht(\p_i) = \dim A-1$.  
Proposition \ref{gamma-1} implies then that $\gamma(A)$ divides $\gcd_i(\gamma(A/\p_i))$. 
Using \ref{fact-1} (1) applied to ${\mathcal O}_C(C)$, we find that 
$\gamma(A)$ divides $\ord(f)$ if $\gamma(A)$ divides 
$\ord(f|_{C_i})$  for each $i$.
Since $A/\p_i$ is an integral domain, $\ord(g):= \ell((A/\p_i)/(g))=e((g),A/\p_i)$ for any non-zero $g \in A/\p_i$.
It follows from its definition that $\gamma(A/\p_i)$  divides $\ord(g)$. 
Hence, $\gamma(A) $ divides $ n(A)$. 
\qed

\bigskip
By construction, $n(A)=n(A/\sqrt{0})$.  Let 
$\widehat{A}$ denote the completion of $A$ with respect to its  maximal ideal. 
Then it is clear from the definition that $\gamma(A) = \gamma(\widehat{A})$.
Note however that in general,
$\gamma(A)\neq \gamma(A/\sqrt{0})$ and, as we show in our next example, $n(A) \neq n(\widehat{A})$ in general.

\begin{example} 
There exists a local noetherian domain $A$ of dimension $1$
 such that its completion $ \widehat{A}$   has a single non-zero minimal prime ideal $\mathfrak P= \sqrt{(0)}$,
and such that $\widehat{A}^{\;  red}:=\widehat{A}/\mathfrak P$ is a discrete valuation ring (see, e.g., \cite{Bennett2}, (3.0.1)).
The ring $A$ 
satisfies the hypotheses of Theorem \ref{thm.na},
so 
$$n(A) = \gamma(A)= \gamma(\widehat{A})=\length({\widehat{A}}_{\mathfrak P} ) \gamma(\widehat{A}^{\; red})= \length({\widehat{A}}_{\mathfrak P} ),$$
while $n(\widehat{A}) = n(\widehat{A}^{\; red})=1 $.
Hence,  $n(A) > n(\widehat{A})$.
We provide  in   \ref{rem.non-uc2}
an example of a noetherian local ring $A$ with henselization $A^h$ such that  $n(A) > n(A^h)$.
\end{example}

\end{section}

\begin{section}{A generic Moving Lemma}

\begin{emp} \label{nx}
We extend to 
schemes  the definition of the moving multiplicity in \ref{def-na} as follows. 
Let $X$ be a noetherian scheme 
and let $x_0\in X$ be a point with $\dim\cO_{X,x_0}\ge 1$. 
Define 
$$n_X(x_0):=n(\cO_{X,x_0}).$$
The main result in this section is Theorem \ref{mv} below, which details one of the most useful uses of the 
invariant $n_X(x_0)$. 
\end{emp}

We start with two preliminary propositions. The first one is proved in \cite{Boutot}, II.9.3, for
affine noetherian schemes.  Recall that an $\FA$-scheme $X$ is a scheme  
such that every finite subset of $X$ is contained in  an 
affine open subset of $X$ (\ref{def-FA}).

\begin{proposition}[{\bf Moving Cartier divisors}] \label{Cartier-mv} 
Let $X$ be a  noetherian $\FA$-scheme. Let $D$ be a Cartier divisor on 
$X$ and let $\mathcal F$ be a finite subset of $X$. Then there exists 
a Cartier divisor $D'$ on $X$, linearly equivalent to 
$D$ and such that $\operatorname{Supp}(D')\cap {\mathcal F}=\emptyset$. 
\end{proposition}
\proof Let $x_1,...,x_m$ be closed points of $X$ such that 
every point of $\mathcal F$ and of $\Ass(X)$ specializes to some of the $x_i$'s.
In particular, for any open subset $V$ of $X$ containing $\{x_1,...,x_m\}$,
the natural restriction map $\mathcal K_X(X)\to \mathcal K_X(V)$ is an isomorphism 
(see \cite{Liubook}, 7.1.15).

Let  $U$ be an affine open subset containing $x_1,...,x_m$.
Recall that by definition, $\cO_X(D)$ is an invertible subsheaf of $\mathcal K_X$.
The canonical homomorphism 
$$\cO_X(D)(U)\to \oplus_{1\le i\le m} \cO_X(D)\otimes k(x_i)$$
is surjective (Chinese Remainder Theorem), so there exists $f\in H^0(U, \cO_X(D)|_U)\subseteq 
\mathcal K_X(U)$ such that for all $i=1,\dots, m$, the image of $f$ in 
$\cO_X(D)\otimes k(x_i)$ 
is a basis. Then $f_{x_i}$ is a basis of $\cO_X(D)_{x_i}$ for all $i=1,\dots, m$. 
This implies that there exists an open subset $V \subseteq U$ containing $\mathcal F$ and $\Ass(X)$
such that $f_x$ is a basis of $\cO_X(D)_x$ for all $x\in V$. 
It follows that $f \in \mathcal K^*_X(V)$.
Extend $f$ to $g\in \mathcal K_X^*(X)$ using the isomorphism $\mathcal K_X(X)\to \mathcal K_X(V)$.
Then $D-\dv(g)$ has support disjoint from $\mathcal F$. 
\qed 

\medskip
Let $X$ be a noetherian scheme, and let $F$ be a closed subset
of $X$. For convenience, we will say that a cycle $Z$ on  $X$ \emph{generically avoids $F$} 
if no irreducible cycle occurring in 
$Z$ is contained in $F$. In particular, no irreducible component
of $\Supp Z$ is contained in $F$. 

\begin{proposition} \label{pro.lemma}
Let $X$ be a noetherian $\FA$-scheme.
Let $F$ be a closed subset of $X$ of positive codimension in $X$.
Let $x_0\in X$, 
and $j : \Spec\cO_{X,x_0}\to X$ be the
canonical injection. Denote again by $x_0$ the 
closed point of $\Spec\cO_{X,x_0}$, and let $F' := j^{-1}(F)$.
Suppose that there exist 
integral closed subschemes $C_1,\dots, C_r$ 
of $\Spec\cO_{X,x_0}$ of dimension $1$, elements $f_i\in k(C_i)^*$ for $i =1,\dots, r$, and an integer $n \geq 1$, such that $C_i \cap F' = \{x_0\}$ 
and 
$n[x_0]=\sum_i [\dv(f_i)]$ in ${\mathcal Z}(\Spec\cO_{X,x_0})$.

Then the cycle 
$n[\overline{\{ x_0 \}}]$ in  ${\mathcal Z}(X)$ is rationally equivalent in $X$ 
to a cycle which generically avoids $F$. 
More precisely, when $x_0 \in F$, let $\ol{C}_i$ be the scheme-theoretic closure
of $j(C_i)$ in $X$.
Then 
$n[\overline{\{ x_0 \}}]$ is rationally equivalent in $\cup_i \overline{C_i}$ 
to a cycle $Z$ 
which generically avoids $F \cap (\cup_i \overline{C_i})$, and such that each irreducible cycle occurring in $Z$ is of codimension $1$ in some $\overline{C_i}$.  
\end{proposition}
\proof
 Since $\ol{C}_i$ is the scheme-theoretic closure
of $j(C_i)$ in $X$ and $C_i \cap F' = \{x_0\}$, the closed subset $\ol{C}_i$ is not contained in $F$.
For each $i=1,\dots, r$, there exists by hypothesis a function $g_i \in k(\ol{C}_i)^*$
defined on a dense open subset of $\ol{C}_i$, whose stalk in $\mathcal O_{\ol{C}_i,x_0}$ is  $f_i$, 
and such that $[\dv_{\ol{C}_i}(g_i)]= \ord_{x_0}(f_i) [\overline{\{ x_0 \}}] + Z_i$
for some cycle $Z_i$ on $\ol{C}_i$ whose support does not contain $\overline{\{ x_0 \}}$.
To conclude the proof, it is enough to show that for all $i=1,\dots, r$, 
$Z_i$ is rationally equivalent on $\ol{C}_i$ 
to a cycle which generically avoids $F\cap \ol{C}_i$. Then  
$n[\overline{\{ x_0 \}}]$ is rationally equivalent on $X$ 
to a cycle which generically avoids $F$, as desired. 

Dropping the subscript $i$ for ease of notation, we are in the following
situation. On the noetherian integral $\FA$-scheme $\ol{C}$, we have 
a cycle $Z$, all of whose irreducible components are of codimension $1$ in 
$\ol{C}$ and whose generic points are all contained in a dense open subset $U$ of $\ol{C}$.
Moreover, $Z|_U=[\dv(g)]$. 
 We claim that  
 $Z$ is rationally equivalent on $\ol{C}$ to a cycle which generically avoids a proper closed subset 
$F \cap \ol{C}$ of $\ol{C}$. 

Indeed, let $\mathcal F$ be the set of points of $F \cap \ol{C}$ of 
codimension $1$ in $\ol{C}$.  Let $V:=\ol{C} \setminus \Supp Z$ and $W := U \cup V$. Then $W$ is open, 
and  $Z|_W =[D]$, where $D$ is the 
 Cartier divisor on $W$ given by the charts $\{ (U,g), (V,1) \}$. By construction, $W$ 
 contains all codimension $1$ points of $\ol{C}$, and so contains $\mathcal F$.
Since $W$ is an $\FA$-scheme
(\ref{FA-g} (2)), we can use 
Proposition \ref{Cartier-mv} to find an invertible 
rational function $f\in k(W)^*$ such that 
$\Supp [D+\dv(f)]$ does not meet $\mathcal F$. 
Let $h\in k(\ol{C})^*$ be the unique 
rational function on $\ol{C}$ extending $f$, and consider
$Z':=Z+[\dv(h)]$. 
Then $Z'|_W=[D+\dv(f)]$ and,  hence, $\Supp(Z')$ does not meet $\mathcal F$ 
either. In other words, the irreducible components of 
$\Supp(Z')$ are not contained in $F$. As $Z'$ only involves irreducible
cycles of codimension $1$, it generically avoids $F$. 
\qed 

\begin{theorem}[{\bf Generic Moving Lemma}] \label{mv} 
Let $X$ be a noetherian $\FA$-scheme.
Let $F$ be a closed subset of $X$ of positive codimension in $X$.
Let $x_0\in X$, and consider the cycle $[\overline{\{ x_0 \}}]$  
in  ${\mathcal Z}(X)$.  
Then 
$n_X(x_0)[\overline{\{ x_0 \}}]$ is rationally equivalent in $X$ 
to a cycle 
which generically avoids $F$. 
\end{theorem}

\proof 
The theorem is obvious when $x_0 \notin F$. So assume now that 
$x_0 \in F$. Since ${\rm codim}(F,X) >0$, none of the irreducible components of $F$
are irreducible components of $X$. Hence, 
the preimage $F'$ of $F$ under the natural map  $\Spec\cO_{X,x_0} \to X$ has dimension 
smaller than $\dim\cO_{X,x_0}$.  In particular, $\dim\cO_{X,x_0} \geq 1$,
and we can apply 
Lemma \ref{lem-n(U,A)}. Denote again by $x_0 $ the closed point of $\Spec\cO_{X,x_0}$.
Let then $U_0$ be a dense open set of 
$\Spec\cO_{X,x_0} $ contained in $\Spec\cO_{X,x_0} \setminus F'$, and such that $n(U_0,\cO_{X,x_0}) = n(\cO_{X,x_0}) = n_X(x_0)$ (\ref{lem-n(U,A)}). 

We can thus find
integral closed subschemes $C_1,\dots, C_r$ 
of $\Spec\cO_{X,x_0}$ of dimension $1$, and $f_i\in k(C_i)^*$, such that $C_i \cap F' = \{x_0\}$ 
and $n[\ol{\{x_0\}}]=\sum_i [\dv(f_i)]$ in ${\mathcal Z}(\Spec\cO_{X,x_0})$, with $n=n_X(x_0)$.
The theorem follows from Proposition
\ref{pro.lemma}. \qed

\medskip
The proof of Theorem \ref{dx-ds} uses 
only the following version of Theorem \ref{mv}, whose proof does not require
the definition and main properties of the invariant $n_X(x_0)$ discussed in the previous section.

\begin{theorem}  \label{mvsimple} 
Let $X$ be a noetherian $\FA$-scheme.
Let $F$ be a closed subset of $X$ of positive codimension in $X$.
Let $x_0\in X$, and consider the cycle $[\overline{\{ x_0 \}}]$  
in  ${\mathcal Z}(X)$.  
Let $Q$ be a $\m_{X,x_0}$-primary ideal of 
$\cO_{X,x_0}$.
Then 
$e(Q)[\overline{\{ x_0 \}}]$ is rationally equivalent in $X$ 
to a cycle 
which generically avoids $F$.  
\end{theorem}
\proof The proof of this version is completely analogous to the proof of \ref{mv}.
The theorem is obvious when $x_0 \notin F$. So assume now that 
$x_0 \in F$. Since ${\rm codim}(F,X) >0$, none of the irreducible components of $F$
are irreducible components of $X$. Hence, 
the preimage $F'$ of $F$ under the natural map  $\Spec\cO_{X,x_0} \to X$ has dimension 
smaller than $\dim\cO_{X,x_0}$.  In particular, $\dim\cO_{X,x_0} \geq 1$,
and we can apply Theorem \ref{local-mv}. Denote again by $x_0 $ the closed point of $\Spec\cO_{X,x_0}$.
We can thus find
integral closed subschemes $C_1,\dots, C_r$ 
of $\Spec\cO_{X,x_0}$ of dimension $1$, and $f_i\in k(C_i)^*$, such that $C_i \cap F' = \{x_0\}$ 
and $n[\ol{\{x_0\}}]=\sum_i [\dv(f_i)]$ in ${\mathcal Z}(\Spec\cO_{X,x_0})$, with $n=e(Q)$.
The theorem follows from Proposition
\ref{pro.lemma}.
\qed

\begin{remark}
In Theorem \ref{mv}, a multiple of the cycle $[\ol{\{x_0\}}]$
can be moved, but in general the irreducible cycle $[\ol{\{x_0\}}]$ itself cannot be moved.
Indeed, consider  for instance the singular projective curve $X$ over 
${\mathbb R}$ defined by the equation $x^2 + y^2 =0$ in 
${\mathbb P}^2_{\mathbb R}$. 
The singular point $x_0:=(0:0:1)$ is the unique rational point of $X$,
and all closed points of $X\setminus \{ x_0 \}$ are smooth and have
degree $2=[\mathbb C : \mathbb R]$. Therefore the $0$-cycle
$[x_0]$ cannot be rationally equivalent to a $0$-cycle with
support in $X^{\reg}$. 
\end{remark}

\begin{corollary} \label{mv-degree} 
Let $X$ be a scheme of finite type over a field $k$. Let $U$ be any
dense open subset of $X$. Let $x_0\in X$ be a closed point. 
Let $Q$ be a $\m_{X,x_0}$-primary ideal of 
$\cO_{X,x_0}$. Then 
$e(Q)\deg_k(x_0)$ is the 
degree of some $0$-cycle $Z$ with support in $U$. 
If in addition  $X$ is separated, then  $Z$ can be chosen to be 
 rationally equivalent in $X$ to $e(Q)[x_0]$. It follows that when $X$ is reduced, then
$\delta(X^{\reg}/k)$ divides $ \gamma(\cO_{X,x_0})\deg_k(x_0)$.
\end{corollary}

\proof 
Let $V$ be an affine neighborhood of $x_0$ in $X$.
We let $V\to \overline{V}$ be an open dense immersion of $V$ 
into a projective variety $\overline{V}$. We apply  
Theorem \ref{mvsimple} to $\overline{V}$ and the closed subset 
$\overline{V}\setminus (U \cap V)$ of positive codimension.
Then $e(Q)[x_0]$ is rationally equivalent in $\overline{V}$ to a 
$0$-cycle $Z$ with support in $U \cap V$. Since $\overline{V}$ is 
projective, we find that $\deg(e(Q)[x_0])= \deg(Z)$ 
because any principal divisor on a projective curve 
over $k$ has degree $0$ (see, e.g., \cite{Liubook}, Corollary 7.3.18). 

Assume now that $X$ is separated. Apply Theorem \ref{mvsimple}
to find a closed curve $C$ in $X$, a function $f \in k(C)^*$, and a $0$-cycle $Z$ on $C\setminus \{x_0\}$,
such that on $C$, $e(Q)[x_0]- Z= [{\rm div}(f)]$. 
Let $U_1:= C \setminus \Supp Z$, and $U_2:= C \setminus \{x_0\}$.
Then let $D$ be the Cartier divisor given by the pairs $(U_1,f)$ and $(U_2,1)$. By construction, $[D]=e(Q)[x_0]$.

 Since $C$ is separated,
Nagata's Theorem lets us find an open embedding of $C$ into a proper curve $C'$ over $k$, and such a curve is known to be also projective.
 Extend in a natural way $D$ to a Cartier divisor $D'$ on $C'$. Let $F := X \setminus U$.
Using \ref{Cartier-mv}, we can find a Cartier divisor $D''$ linearly equivalent to $D'$
on $C'$ and such that $\Supp(D'')$ does not intersect $(C'\setminus C) \cup (C \cap F)$.
As above, $\deg(D'') = \deg(D')$ since $C'$ is a projective curve.
Since $C \to C'$ is an open immersion, we find that $D''$ restricted to $C$
is equivalent to $D$ on $C$.
\qed 
\medskip

Let $k$ be any field. An \emph{algebraic variety $X$ over $k$} is
a scheme of finite type over $k$ (not necessarily separated).
Let $\mathcal D$ denote the set of all degrees of closed points of $X$.
When $X/k$ is not empty, the {\it index} $\delta(X/k)$ of $X/k$ is the greatest common
divisor of the elements of $\mathcal D$.
 
The following proposition is only slightly more general than
\cite{CT}, page 599, or \cite{Cla}, Lemma 12. See also \cite{CM}, 1.12.

\begin{proposition} \label{pro.deltareg}
Let $X$ be a (non necessarily proper) regular non-empty algebraic variety 
over a field $k$. Then $\delta(U/k) = \delta(X/k)$ for any dense 
open subset $U$ of $X$. 
In particular, if $X_1$ and $X_2$ are two integral regular 
algebraic varieties over $k$ which are birational, then 
$\delta(X_1/k)= \delta(X_2/k)$.
\end{proposition} 

\proof  Clearly, $\delta(X/k) $ divides $ \delta(U/k)$. Let $x_0\in X$ be a 
closed point. Applying Corollary \ref{mv-degree},
$x_0$ has same degree as some $0$-cycle with support in $U$. 
Hence $\delta(U/K) $ divides $ \deg (x_0)$, so $\delta(U/K) $ divides $ \delta(X/K)$. 
\qed 

\begin{remark} \label{rem.multiplicity} 
Let $X/k$ be an integral normal algebraic variety. Let $k'$ be the 
algebraic closure of $k$ in the field of rational functions $k(X)$. 
As $k'/k$ is finite 
(hence integral), we have $k'\subseteq \mathcal O_{X}(X)$. 
Therefore, $X \to \Spec(k)$ factors as $X \to \Spec(k') \to \Spec(k)$.
We will write $X/k'$ when we regard $X$ as a variety over $k'$ 
through the morphism $X \to \Spec(k')$. We have 
$\delta(X/k)=[k' : k]\delta(X/k')$. Note that as a variety over $k'$,
$X$ is geometrically irreducible. 
\end{remark}

Let $X$ be any noetherian scheme. Let $U \subseteq X$ be an open subset. 
Consider the natural map ${\mathcal Z}(U) \to \CH(X)$, 
which sends an irreducible cycle on $U$ to its Zariski closure in $X$ 
modulo rational equivalence. Denote by $\CHU$ the cokernel of ${\mathcal Z}(U) \to \CH(X)$.
Our next proposition generalizes \cite{SS}, Proposition 7.1, 
where $X$ is assumed to be integral and regular, and 
of finite type over the
spectrum of a discrete valuation ring. 

\begin{proposition} \label{U-X} 
Let 
$U$ be a dense open subset of a noetherian $\FA$-scheme $X$. 
\begin{enumerate}[\rm (1)] 
\item If $X$ is regular, then $\CHU=(0)$. 
\item $\CHU$ is a torsion group. Let $x_0 \in X$ be a point with $\dim\cO_{X,x_0}\ge 1$. 
Then the class of $[\ol{\{x_0\}}]$ in the
cokernel 
$\CHU$ has order dividing $n_X(x_0)$,
and there exists one such open subset $U$ where the order is exactly $n_X(x_0)$.
\item $\CHU$
has finite exponent in the following situations: 
\begin{enumerate}[\rm (a)] 
\item $X$ is (quasi-)excellent. 
\item There exists a morphism of 
finite type $f : X\to S$ with $S$ noetherian regular and
$f$ either 
open or equidimensional.
\end{enumerate} 
\end{enumerate}
\end{proposition}

\proof (1) and (2) are immediate consequences of Theorem \ref{mv}. Part 
(3) results from the fact that the set $\{ e(\cO_{X,x}) \mid x\in X\}$
is bounded, as discussed in Proposition \ref{cor.mv} (3) and Remark \ref{quasi-excellent}
below. 
\qed 

\begin{remark}
If $V \subset U$ are two dense open subsets of an integral scheme $X$, then there is a natural surjective 
group homomorphism $f_{V,U}:{\mathcal A}(X,V) \to \CHU$. Let ${\mathcal O}$ denote 
the set of all dense open subsets $U$ of $X$, with the  relation $V \leq U $ if $V \supseteq U$.
Define  $\SCH(X)$ to be the inverse limit of the projective system $\{\CHU\}_{U \in {\mathcal O}}$.
Let $x_0 \in X$ be a point with $\dim\cO_{X,x_0}\ge 1$. Then  Proposition \ref{U-X}
implies that the class of $[\ol{\{x_0\}}]$ in $\SCH(X)$ has order $n_X(x_0)$, and that
when $X$ is quasi-excellent, 
$\SCH(X)$ is a torsion group of finite exponent.
In particular, since $n_X(x_0)=1$ when $x_0$ is regular, the group ${\mathcal S}(X)$ is generated by the classes of certain irreducible cycles $[\overline{\{x\}}]$ with $x$   singular on $X$.
\end{remark}

\begin{proposition} \label{cor.mv} 
Let $f: X\to S$ be a morphism of finite type over a noetherian 
scheme $S$. Then
\begin{enumerate}[\rm (1)]
\item $X$ is  the union of finitely many open subsets $X_i$,
each endowed with a quasi-finite $S$-morphism $f_i : X_i \to W_i:=\mathbb A^{d_i}_S$, 
such that every $x\in X$ belongs to some $X_i$ with 
$\dim_x X_{f(x)}=d_i$. 
\item Assume that $S$ is irreducible. Suppose in addition that either $f$ is 
locally\footnote{See \cite{EGA} IV.13.2.2, with a correction 
in \cite{EGA} (${\bf Err}_{\rm IV}$, {\bf 34}) on pp. 356-357 of no. {\bf 32}.}  equidimensional and $S$ is universally 
catenary, or that $f$ is open. 
Then, with $x \in X_i$ as in (1),  
$\dim\cO_{X,x}= \dim\cO_{W_i, f_i(x)}$. 
\item Suppose that $S$ is regular and that $f$ is either open 
or locally equidimensional. 
Then the set of multiplicities $\{ e(\cO_{X,x}) \mid x\in X \}$ is bounded. 
\end{enumerate}
\end{proposition}

\proof (1) By \cite{EGA}, IV.13.3.1.1, for each $x \in X$, there 
exists an open neighborhood $U_x$ of $x$ and a quasi-finite 
$S$-morphism $U_x \to \mathbb A^d_S$, where $d=\dim_{x} X_{f(x)}$.
For any integer $n\ge 0$, let 
$F_n:=\{ x\in X \mid \dim_x X_{f(x)}\ge n\}$ and
$G_n=F_n\setminus F_{n+1}$. By a
theorem of Chevalley, $F_n$ is closed (\cite{EGA}, IV.13.1.3), 
and $F_n=\emptyset$ if $n\ge n_0$ for 
some $n_0$ (\cite{EGA}, IV.13.1.7). For all $n\le n_0$, 
$G_n$ is quasi-compact and,  hence, can be covered by finitely many 
open subsets $\{ U_{n,j}\}_j$ belonging to the covering $\{U_x, x \in X\}$.
Since $x\in X$ belongs to $x\in G_d$ for
$d=\dim_x X_{f(x)}$,  $x$ belongs to some $U_{d,j}$. 

(2) Suppose that $x \in X_i$, with $f_i : X_i \to W_i:=\mathbb A^{d_i}_S$
such that $\dim_x X_{f(x)}=d_i$. Clearly, 
$\dim \cO_{W_i, f_i(x)}= \dim \cO_{X_{f(x)}, x} + \dim\cO_{S,f(x)}$.
We thus need to show that 
$$\dim\cO_{X,x} = \dim \cO_{X_{f(x)}, x} + \dim\cO_{S,f(x)}.$$ 
That this equality holds when $f$ is flat is remarked in 
\cite{EGA}, IV.13.2.12 (i), and proved for $f$ open in \cite{EGA}, IV.14.2.1.
 The other case  follows from \cite{EGA}, IV.13.3.6. 

(3) Since $S$ is regular, it is the disjoint union of finitely many 
integral regular schemes, and we are reduced to consider the case where $S$ is integral. 
Let $f_i$ be as in (1). For every $i$, there exists an open immersion of $u: X_i \to Z_i$ and a finite
morphism $g_i : Z_i\to W_i$ with $g_i \circ u = f_i$ (Zariski's Main Theorem \cite{EGA} IV.8.12.6). 
Our hypotheses allow us to apply (1) and (2), and 
every $x\in X$ belongs to some $X_i$ with $\dim\cO_{X,x}=\dim\cO_{W_i, f_i(x)}$.
For such an $x$, 
we find that there exists an affine open neighborhood of $f_i(x)$ such that 
$\cO_{W_i}(U)$ injects into $\cO_{Z_i}(f_i^{-1}(U))$, since both of these rings have the same dimensions.
Since $g_i$ is finite, there exists a 
surjective homomorphism 
$\cO_{W_i}^{m_i}\to g_{i*}(\cO_{Z_i})$. Then, for all $x\in X_i$ such
that $\dim\cO_{X,x}=\dim\cO_{W_i, f_i(x)}$, 
we have $e(\cO_{X,x})\le m_i e(\cO_{W_i, f_i(x)})$,  and since $W_i$ is regular by hypothesis, $e(\cO_{W_i, f_i(x)})=1$
(use \cite{Z-S}, VIII, \S 10, Corollary 1 to Theorem 24, along 
with the discussion (3) following Corollary 1). Therefore, the set of the
multiplicities $e(\cO_{X,x})$ is bounded by $\max_i \{ m_i \}$. 
\qed 

\begin{remark} \label{quasi-excellent} 
  Let $X$ be a noetherian excellent 
scheme of finite dimension. Then only a finite number of polynomials appear 
as Hilbert-Samuel polynomials of $\cO_{X,x}$ when $x$ varies in $X$.
In particular, 
the set of the multiplicities $\{ e(\cO_{X,x}) \mid x \in X\}$ is
bounded. This statement is found in \cite{Bennett}, Remark III.1.3, page 83.
Let us explain now how  the latter conclusion still holds if one only assumes that 
$X$ is a noetherian quasi-excellent scheme. 

Recall that an excellent scheme is a quasi-excellent and universally catenary scheme.
 This latter hypothesis is not needed in \cite{Bennett}, Remark III.1.3:
 One finds in the discussion following  \cite{Bennett}, II (2.2.1), page 34,
that in the use of normal flatness, the needed hypothesis is that the regular locus
of a closed subset $Y$ of $X$ is open in $Y$. This hypothesis is satisfied when $X$ is quasi-excellent. 

Let $X$ be a noetherian quasi-excellent  
scheme. The results of \cite{Bennett}
 show that the set of the multiplicities $\{ e(\cO_{X,x}) \mid x \in X\}$ is
bounded, even when $X$ is not assumed to have finite dimension.
Indeed, first note that there is a finite decreasing sequence of closed subsets
 $X=F_0 \supset F_1 \supset \ldots \supset F_n = \emptyset$
 such that $F_i \setminus F_{i+1}$ with the reduced  structure is regular and 
 $X$ is normally flat along $F_i \setminus F_{i+1}$. This uses ``noetherian induction'' and one does not 
 need to argue using the finite dimensionality of $X$ as in \cite{Bennett}, Remark III.1.3. The relevant statement 
 occurs in \cite{Bennett}, page 28, just before Theorem (2): For a local noetherian ring ${\mathcal O}$ with Hilbert-Samuel function
 $H^{(1)}_{\mathcal O}$ (see page 26) and
with  $P$ a prime ideal 
of coheight $c$ and ${\mathcal O}/P $ regular,
$H^{(1)}_{\mathcal O} =  H^{(1+c)}_{{\mathcal O}_P} $ if and only if $   \Spec {\mathcal O}$ normally flat along 
$\Spec {\mathcal O}/P$ at the closed point. The direction that we need is stated in \cite{Bennett}, (2.1.2) Chapter 0, page 33. This implies that if a noetherian scheme $X$ is normally flat along a connected locally closed regular scheme $Y$, then the multiplicity of ${\mathcal O}_{X,y}$ $(y \in Y)$ is constant on $Y$.
 \end{remark}
 
 \begin{example}
 A noetherian domain $A$ of dimension $1$  with the following properties is exhibited in \cite{HL}, 3.2:
$\Spec A$ has infinitely many singular maximal ideal $\m_i$, $i \in I$, and such a domain can be found 
even when, for each $i$, the $efg$-numbers at $\m_i$ are specified in advance
(where $\m_iB= \prod_{j=1}^{g(i)}\n_{ij}^{e(\n_{ij}/\m_i)}$ 
in the integral closure $B$ of $A$, and  $f(\n_{ij}/\m_i):=[B/\n_{ij}: A/\m_i]$). 
Moreover, for each $i$, $B_{\m_i}$ is a finitely generated $A_{\m_i}$-module (\cite{HL}, 3.3).
Choosing such an example with $g(i)=1 $ for all $i$ and $\lim_i f(\n_{i1}/\m_i)= \infty$
produces an example of a noetherian scheme  $X$ of dimension $1$ 
such that the set $\{ e(\cO_{X,x}) \mid x\in X\}$ is 
not bounded. Indeed, since $A$ is a domain of dimension $1$, there exist non-zero elements $a_i,b_i\in A$ such that $ e(\m_i,A)= e((a_i), A) - e((b_i), A)$ (\ref{local-mv}).
It follows from \ref{fact-1} (3) that both $e((a_i), A)$ and $ e((b_i), A)$ are divisible by $f(\n_{i1}/\m_i)$.
\end{example}

\end{section}
\begin{section}{A different perspective on the index} \label{perspective}

Let $X$ be a noetherian scheme. Let $x_0\in X$ be a point with $\dim\cO_{X,x_0}\ge 1$. 
We study in this section the invariant $n_X(x_0)$ introduced in \ref{nx}. 
In particular, let $f : Y\to X$ be a 
morphism of finite type, and set $E:=Y \times_X \Spec k(x_0)$.
We relate in  \ref{compute-n} the invariant $n_X(x_0)$ with the index of the scheme $E/k(x_0)$.
This leads us in Corollary \ref{newdescriptionIndex} to provide a new way of computing the index of a regular 
closed subvariety $X$
of a projective space using data pertaining only to the singular vertex of a cone over $X$. 

\begin{theorem}\label{compute-n} 
Let $X$ be a noetherian scheme. 
Let $f : Y\to X$ be a 
morphism of finite type
such that the generic point of every irreducible component of  $Y$ maps 
to the generic point of an irreducible component of $X$. 
Let $x_0\in X$ be a point with $\dim\cO_{X,x_0}\ge 1$, and set $E:=Y \times_X \Spec k(x_0)$. 
Assume that $E \neq \emptyset$.
\begin{enumerate}[{\rm (a)}] 
\item Assume that $f$ is birational and proper. Then $\delta(E/k(x_0)) $ divides $ n_X(x_0)$. 
\item Assume that $f$ is birational and finite. Then
$$ \gcd \{ n_Y(y)[k(y) : k(x_0)] \mid y \in f^{-1}(x_0), y 
\ {\rm  closed} \}  \mbox{ \it divides } n_X(x_0).$$
\item Assume that $\cO_{X,x_0}$ is universally catenary. 
Then 
$$n_X(x_0)\mbox{ \it divides } \gcd \{ n_Y(y)[k(y) : k(x_0)] \mid y \in f^{-1}(x_0), y 
\ {\rm  closed} \}.$$
\end{enumerate}
\end{theorem}

\proof 
After the base change $\Spec\cO_{X, x_0}\to X$, we can suppose that 
$X$ is local with closed point $x_0$.

(a)
Let us show that $\delta(E/k(x_0))$ divides $ n_X(x_0)$. 
Let $U$ be any dense open subset of $X$ such that 
$f^{-1}(U)\to U$ is an isomorphism. Since $U$ is dense, there exists, by definition of $n({\mathcal O}_{X,x_0})$,
a one-dimensional reduced closed subscheme $C$ of $X$
such that $\max(C)\subseteq U$ and a rational function
$g\in \mathcal K^*_C(C)$ such that $n_X(x_0)[x_0]=[\dv(g)]$.
Let $\tilde{C}$ be the strict transform of $C$ in $Y$. Then we 
have a finite birational morphism $\pi : \tilde{C}\to C$, and 
$$n_{X}(x_0)[x_0]= 
\sum_{y\in \pi^{-1}(x_0)} \ord_y(g)[k(y) : k(x_0)][x_0]$$ 
(use \ref{fact-1} (3)). 
As $\pi^{-1}(x_0)\subseteq E$, we conclude that $\delta(E/k(x_0))\mid n_X(x_0)$.

(b) We will proceed as in the proof of (a), after 
carefully choosing the initial dense open subset $U$ of $X$.
First, we may assume that $X$ and $Y$ are affine.
Let $y_1,
\dots, y_n$ denote the preimages of $x_0$.
For each $i=1,\dots, n$, choose a dense open subset $V_i$ of $Y$ such that
$n(V_i \cap \Spec \cO_{Y,y_i}, \cO_{Y,y_i}) = n_Y(y_i)$ (see \ref{lem-n(U,A)}).
Since $f$ is birational, let $V \subseteq  \Spec B$ be a dense open subset such that $f_{|V} $ is an isomorphism.
Then $W:=(\cap_i V_i) \cap V$ is a dense open subset of $Y$, contained in each $V_i$, 
and such that $n(W \cap \Spec \cO_{Y,y_i}, \cO_{Y,y_i}) = n_Y(y_i)$. We let $U:=f(W)$, a dense open subset in $X$ 
 with $f^{-1}(U)=W$ by construction.
As in the proof of (a),
we find that $$n_{X}(x_0)[x_0]= 
\sum_{y\in \pi^{-1}(x_0)} \ord_y(g)[k(y) : k(x_0)][x_0].$$ 
By our choice of $W$, we find that for each $y\in \pi^{-1}(x_0)$, $n_{Y}(y)$ divides $\ord_y(g)$,
and the result follows.

(c) Recall that we assume $X=\Spec \cO_{X,x_0}$. Let $y_0$ be a closed point 
in $E$. We  claim that $n_X(x_0) $ divides $ n_Y(y_0)[k(y_0) : k(x_0)]$.
Let $W$ be an open affine neighborhood of $y_0$ in $Y$.
Since $W$ is of finite type over the affine scheme $X$, 
we can find a projective scheme $g:Z \to X$ and an open embedding $i:W \to Z$
such that $g \circ i = f_{\mid W}$. Since $n_Y(y_0) = n_W(y_0)= n_Z(y_0)$,  it suffices to prove 
our claim in the case where $f$ is projective.

Using \ref{lem-n(U,A)}, we may find  a dense open subset $U$ of $X$ which does not contain $x_0$, and such that $n_X(x_0)= n(U,\cO_{X,x_0})$.
Let $V:=f^{-1}(U)$. Then 
$V$ is dense in $Y$. Let $s: \Spec \cO_{Y,y_0} \to Y$ denote the natural morphism,
and let $V':= s^{-1}(V)$. Since $V'$ is dense in  $\Spec \cO_{Y,y_0}$, 
we can find finitely many closed integral subschemes ${C}_i$ of $\Spec \cO_{Y,y_0}$ with generic points in $V'$, 
and an invertible rational function $g$ on $\cup_i C_i$ such that $\ord_{y_0}(g) = n_Y(y_0)$. 
Since $Y$ is noetherian and $\FA$, we can apply Proposition \ref{pro.lemma}
and find that 
$n_Y(y_0)[ y_0 ]$ is rationally 
equivalent to a cycle $Z$ on $Y$ such that 
$\max(\Supp Z)\subseteq V$. 
More precisely,
there exist finitely many closed integral subschemes $\overline{C}_i$ in $Y$ such that $ y_0 $ has codimension $1$ 
in   $\overline{C}_i$, such that the generic point of $\overline{C}_i$ is in $V$ for all $i$, 
and 
such that $n_Y(y_0)[y_0]$ is rationally equivalent on $S:= \cup_i \overline{C}_i$ to a cycle 
$Z$ with $\max(\Supp Z)\subseteq V$, and such that each irreducible cycle occurring in $Z$
is of codimension $1$ in some $\overline{C}_i$. 

Let $\overline{D}_i$ denote the schematic image of $\overline{C}_i$ in $X$.
Then $\overline{D}_i$ is universally catenary.
We can thus use the dimension formula (see \cite{EGA}, IV.5.6.5.1) for the morphism  $\overline{C}_i \to \overline{D}_i$,
$$\dim\cO_{\overline{C}_i, z_0}+{\rm trdeg} (k(z_0)/k(x_0)) = \dim\cO_{\overline{D}_i,x_0} + {\rm trdeg} (k(\eta_i)/k(\xi_i)),$$
where $z_0$ is any closed point of $\overline{C}_i$, and $\eta_i$ and $\xi_i$ denote respectively the generic point of $\overline{C}_i$ and $\overline{D}_i$.
Since $\dim\cO_{\overline{C}_i, y_0}=1$ by construction, and ${\rm trdeg} (k(y_0)/k(x_0))=0$ by hypothesis, 
we find that $\dim\cO_{\overline{D}_i,x_0} =1$  and $ {\rm trdeg} (k(\eta_i)/k(\xi_i))=0$.
It follows that $\dim\cO_{\overline{C}_i, z_0}=1$ for all $z_0$ closed in $\overline{C}_i$, so that 
 $\dim(\overline{C}_i)=1$. Since $x_0$ is the only closed point of $\overline{D}_i$, we find that $\dim(\overline{D}_i)=1$.
Since $\dim(\overline{C}_i)=1$, we find that $\max(\Supp Z)$ consists in a set of closed points of $Y$.
Since $V$ does not contain any closed point of $Y$ by construction, we find that $n_Y(y_0)[ y_0 ]$ is rationally 
trivial on $Y$.

The morphism $f_{|S}: S \to f(S)$ is proper, with $f(S) $ universally catenary. 
Hence, the cycle $f_*(n_Y(y_0)[y_0])$ is rationally trivial on $f(S)$ and, thus, on $X$
(use \cite{Th}, 6.5 and 6.7).
Since the irreducible components of $S$ have dimension $1$, since the generic points of $f(S)$ belong to $U$ by construction, 
and since $f_*(n_Y(y_0)[y_0])=n_Y(y_0)[k(y_0) : k(x_0)][x_0]$ is rationally trivial on $f(S)$, we find that by definition,
$n(U, \cO_{X,x_0}) $ divides $ n_Y(y_0)[k(y_0) : k(x_0)]$,
and the statement of (c) follows.
\qed 

\begin{remark} \label{rem.non-uc}
The hypothesis that $X$ is universally catenary is needed in \ref{compute-n} (c).
Indeed, consider  the finite birational morphism $\pi: Y \to X$ described
in Example \ref{non-uc}. The scheme $Y$ is regular, and $X$  is catenary but not universally catenary.
The preimage of the point $x_0 \in X$ consists in the two regular points $y_0$ and $y_1$ in $Y$,
with $[k(y_0):k(x_0)]=d$, and $[k(y_1):k(x_0)]=1$.
It follows from the discussion in \ref{non-uc} that $n_X(x_0)$ is divisible by $d$. Thus, when $d>1$, the morphism $\pi$ 
provides an example where   
$$\gcd \{ n_Y(y)[k(y) : k(x_0)] \mid y \in f^{-1}(x_0), y 
\ {\rm  closed} \}= \delta(E/k(x_0))<n_X(x_0).$$ 
\end{remark} 

\begin{corollary}\label{cor.compute-n} 
Let $X$ be a universally catenary noetherian scheme. Let $x_0\in X$ be a point
with $\dim\cO_{X,x_0}\ge 1$. 
Let $f : Y\to X$ be a 
proper birational morphism such that 
$f^{-1}(x_0)$ is contained in the regular locus of $Y$. Let $E:=Y \times_X \Spec k(x_0)$. Then 
 $n_X(x_0)=\delta(E/k(x_0))$.
 
In particular, if $X$ is an integral excellent scheme of dimension $1$, and
if $f : Y\to X$ is the normalization morphism, then 
$n_X(x_0)=\gcd\{ [k(y) : k(x_0)] \mid y\in f^{-1}(x_0)\}$. 
\end{corollary}
\proof The statement follows immediately from the previous theorem, since $n_Y(y) =1 $ for any 
 $y \in f^{-1}(x_0)$ because $y$ is regular on $Y$ by hypothesis.
\qed

\begin{corollary} \label{newdescriptionIndex}
Let $K$ be any field. Let $V/K$ be a  regular 
closed integral subscheme 
of ${\mathbb P}^n/K$, 
and denote by $W/K$ a 
cone over $V/K$ in ${\mathbb P}^{n+1}/K$.  Let $w_0$ denote the vertex of  $W$.
Then $\delta(V/K) = n(\cO_{W,w_0})=\gamma(\cO_{W,w_0})$.
\end{corollary}
\proof Let $Z \to W$ denote the blow-up of the vertex $w_0$ of the cone $W$.
It is well-known that the exceptional divisor $E$ of $Z \to W$ is isomorphic to $V/K$.
When $V/K$ is regular, 
the points of $E$ are regular on $Z$. 
Theorem \ref{thm.na}, along with  \ref{cor.compute-n}, shows that  $\delta(V/K)  = n(\cO_{W,w_0}) = \gamma(\cO_{W,w_0})$. \qed

\medskip

Keep the notation of Corollary \ref{newdescriptionIndex}. Consider the  sets  
$${\mathcal D}(V/K):= \{ \deg_K(P), P \text{ closed point of } V\}$$ and $$
{\mathcal E}(\cO_{W,w_0}) := \{ e(Q, \cO_{W,w_0}), Q {\rm \ is \ } {\rm primary}  \}.$$
Corollary \ref{newdescriptionIndex} shows that 
$\gcd(d, d \in {\mathcal D}) = \gcd(e, e \in {\mathcal E})$.
 Let us note here that the sets $\mathcal D$ and $\mathcal E$ can be very different.
For instance, it is known that $e(\m_{W,w_0}) = \deg(V)$, and this integer  is the minimal element 
in $\mathcal E$. At the same time, it may happen that $V(K) \neq \emptyset$, so that
$1 \in {\mathcal D}(V)$. 

If $Q$ is a primary ideal of $\cO_{W,w_0} $ generated by a system of parameters, 
then any positive multiple of $e(Q, \cO_{W,w_0})$ belongs to ${\mathcal E}$ (\cite{SH}, 11.2.9).
When $K$ is infinite, this statement remains true for all $\m_{w_0}$-primary ideals of $\cO_{W,w_0} $
(use \cite{Z-S}, Thm. 22).
In view of the property of the set ${\mathcal D}$ in Proposition \ref{numberfields}, one may wonder whether an analogous property holds for the 
set ${\mathcal E}$.   

In the next proposition, we call a field $K$  {\it Hilbertian} if every separable Hilbert set of $K$ is non-empty.
This is the definition adopted  in \cite{F-J}, 12.1, but not the one used in \cite{LangDG}, page 225, for instance.
Both definitions agree when $\charr(K)=0$, such as in the case of number fields, and the reader 
will see that the proof of \ref{numberfields} is simpler in this case\footnote{The case of inseparable extensions
requires careful consideration (see for instance the Notes on page 230 of \cite{F-J}). We also   note that 
 Proposition 5.2 in Chapter 9 of \cite{LangDG}, page 240, requires further hypotheses to ensure that the element $y$ in its proof exists.}.
A scheme $V/K$ of finite type is called {\it generically smooth} if it contains a dense open subset $U/K$ which is smooth over $K$.

\begin{proposition} \label{numberfields}
Let $K$ be a Hilbertian 
field. Let $V/K$ be an irreducible regular generically smooth algebraic variety
of positive dimension.
Then there exists $n_0>0$ such that $$\{ n \delta(V/K), 
n \geq n_0 \} \subseteq {\mathcal D}(V/K).$$
\end{proposition}

\proof Let $K'$ denote the subfield of elements of $  K(V)$ algebraic over $K$. Then $V/K'$
is geometrically irreducible and $\delta(V/K) = [K':K]\delta(V/K')$
(see \ref{rem.multiplicity}). It is thus sufficient to prove the proposition
when $V/K$ is geometrically irreducible. Choose an affine
open subset $U$ of $V$. Then $\delta(U/K) = 
\delta(V/K)$ (\ref{pro.deltareg}). We can thus find finitely many closed points $P_1,\dots, P_r$ in $U$
such that $\gcd(\deg_K(P_i),i=1,\dots, r)=\delta(V/K)$. In case $\charr(K)>0$, we use \ref{sep.index}
to show that we can assume that each  point $P_i$ has its residue field $K(P_i)$ separable over $K$.
Since $U/K$ is quasi-projective, we can use \cite{MB2}, 2.3,
and find a geometrically 
integral curve $C/K$ in $U$ which contains $P_1,\dots, P_r$
in its 
smooth locus $C^{sm}$ 
(the proof of \cite{MB2}, 2.3, uses  Bertini's Theorem in \cite{J}, 
where the only hypothesis on $K$ is that it is infinite). Let $C'/K$ denote a 
regular 
projective curve containing
$C^{sm}$ as a dense open subset. Clearly, $\delta(C'/K) = \delta(C^{sm}/K) $ 
and $\delta(C^{sm}/K)$ divides $\gcd(\deg_K(P_i),i=1,\dots, r)=\delta(V/K)$.

Since $C' \setminus C^{sm}$ is a finite set, it suffices to prove the proposition 
in the case where $V/K$ is the smooth locus of a regular  
projective geometrically integral curve $\overline{V}$, of arithmetical genus $g(\overline{V})= \dim H^1(\overline{V}, {\mathcal O}_{\overline{V}})$.

Let again 
$P_1,\dots, P_r$ in $V$
be such that $\gcd(\deg_K(P_i),i=1,\dots, r)=\delta(V/K)$.
Clearly, every large enough integer multiple  $j$ of $\delta(V/K)$ can
be written as 
$j= \sum_{i=1}^r x_i \deg_K(P_i)$ with $x_i  \geq 0$. Suppose that $j \geq 2g(\overline{V}) -1 + \max_i(\deg_K(P_i))$.
Then by the Riemann-Roch Theorem there exists 
a function $f \in H^0(\overline{V}, \sum_{i=1}^r x_i P_i)$ which does not belong to 
$H^0(\overline{V}, {\mathcal O}(D))$ for any effective
$D<(\sum_{i=1}^r x_i P_i)$.
It follows that $f$ defines  a morphism $f:\overline{V} \to {\mathbb P}^1$ over $K$ of degree equal to $j$. 
If the induced extension $K(\overline{V})/K({\mathbb P}^1)$ is separable, then our assumption that
$K$ is Hilbertian implies that
$f$ has  irreducible fibers, and so, that points
of degree $j$ on $V$ exist. 

If $\charr(K)>0$, we can modify the argument as follows.
Consider a closed point $P_0$ of $\overline{V}$ with residue field $K(P_0)$ separable over $K$, 
and such that $P_0$ is distinct from $P_1,\dots, P_r$.
Then every large enough integer multiple  $j$ of $\delta(V/K)$ can
be written as 
$j= \deg_K(P_0) + \sum_{i=1}^r x_i \deg_K(P_i)$ with $x_i  \geq 0$. 
As before, when $j \geq 2g(\overline{V}) -1 + \max(\deg_K(P_i), i=0,\dots, r)$, we use the Riemann-Roch theorem to  define  a morphism $f:\overline{V} \to {\mathbb P}^1$ over $K$ of degree equal to $j$.
This morphism now has the property that the induced extension $K(\overline{V})/K({\mathbb P}^1)$ is separable, and we can conclude 
using the fact that $K$ is Hilbertian, as above.
 \qed

\begin{remark} The statement \ref{numberfields} does not hold if $V/K$ is not assumed to be regular.
Indeed, consider the curve $X/{\mathbb Q}$ defined by $x^2+y^2=0$. Then $\delta(X/{\mathbb Q})=1$, and 
$(0,0)$ is the only point of odd degree.
\end{remark}

\begin{remark} \label{algorithm}
When  $X/F$ is a smooth proper geometrically connected curve of genus $g >1$ over any field $F$,
slightly more can be said about the set $\mathcal D(X/F)$.
Since the canonical class can be represented by an effective divisor,
the curve $X/F$ has at least one point of degree at most $2g-2$.
Let $Q$ be a point of minimal degree on X.

{\it We claim that there exists a divisor 
$\sum_{i=1}^s a_i P_i$ of degree $\delta(X/F)$ 
such that   for all} $i=1,\dots, s$, $$\deg(P_i) < g+ \deg(Q) \leq 3g-2.$$
Indeed,
let $D:=\sum_{i=1}^s a_i P_i$ be a divisor
of degree $\delta(X/F)$ such that $\max_i (\deg(P_i))$
is minimum among all such divisors. We may clearly choose such a divisor $D$
such that the degrees of the closed points in the support of $D$ are pairwise distinct. 
After reordering if necessary, we may assume that $\deg(P_1) < \ldots < \deg(P_s)$.
If $\deg(P_s) \geq \deg(Q) + g$, then by the Riemann-Roch theorem,
the divisor $P_s-Q$ is equivalent to an effective divisor $\sum_{j=1}^k
b_j Q_j$, with $b_j >0$ and $\deg(Q_j) <\deg(P_s)$.
Replacing $P_s$ by $Q+ \sum_{j=1}^k b_j Q_j$ in $D$ 
contradicts the minimality of $\max_i (\deg(P_i))$.

{\it Suppose now that there exists an algorithm which determines, given an irreducible variety $V/F$,
whether $V/F$ has an $F$-rational point. 
Then there exists an algorithm which determines the index $\delta(X/F)$ of a smooth proper geometrically irreducible curve
$X/F$ of genus $g$}. Indeed, for a given $d \geq 1$, consider  the quotient $X^{(d)}/F$ 
of the $d$-fold product $X^d/F$ by the action 
of the symmetric group $S_d$ acting by permutation on $X^d$. 
Then a $F$-rational point on $X^{(d)}$ corresponds to a point of $X$ defined 
over an extension $L/F$ of degree $[L:F]$ dividing $d$.
To compute $\delta(X/F)$, it suffices 
to determine whether $X^{(d)}/F$ has a $F$-rational point for $d=1,\dots, 3g-3$.
\end{remark}

\begin{example} It goes without saying that most often,
the explicit determination of the set ${\mathcal D}(X/F)$ is a very difficult problem.
Consider for instance the Fermat curve $X_p/{\mathbb Q}$ given by the equation
$x^p + y^p = z^p$, with $p >3$ prime.
This curve has obvious points of degrees $1$, $p-1$, and $p$. Points of degree $2$ on the line
$x+y=z$ were noted already by Cauchy and Liouville (\cite{Tze}, Introduction). 
Indeed, a point $(x:y:1)$ on the intersection of  $X_p$ with the line $x+y=z$ satisfies:
$$  x^p + (1-x)^p - 1  =  x(x - 1)(x^2 - x + 1)^b   E_p(x),$$
with  $E_p(x) \in   {\mathbb Z}[x]$, and $b=1$ if $p \equiv 2 \pmod{3}$ and $b=2$ if $p \equiv 1 \pmod{3}$.
Mirimanoff conjectured in 1903 that $E_p(x)$ is irreducible over ${\mathbb Q}$ for all primes $p\geq 11$.
 Klassen and Tzermias
conjecture in \cite{KT} that any point $P$ on $X_p/{\mathbb Q}$ of degree at most $p-2$ lies on the line $x+y=z$.
Putting these two conjectures together when $p\geq 11$, we would find that 
$${\mathcal D}(X_p/{\mathbb Q}) = \{  1,2, \deg(E_p(x)),p-1,p, \dots \}.$$
\end{example}

\begin{example} 
Let $A$ be any noetherian local ring of  dimension $d>0$.
Set $X = \Spec A$, with closed point $x_0$. Our next example shows that $n_X(x_0)$ is not bounded
when $d$ is fixed.

Let $k$ be a field which has finite extensions of any given degree.
Let $Y=\Spec B$ be an affine normal integral algebraic 
variety of dimension $d>0$ over   $k$. Let $y_0\in Y$ be a closed point corresponding
to a maximal ideal $\m$ of $B$. Let $r\ge 1$, let $A:=k+\m^r$, and let 
$X = \Spec(A)$. Let $\pi : Y\to X$ be the induced morphism. 
Then $\pi$ is finite birational, $x_0:=\pi(y_0)\in X(k)$, and
$\pi: Y\setminus \{ y_0 \} \to X\setminus \{ x_0\}$ is an isomorphism.
Theorem \ref{compute-n} 
shows that 
$n_X(x_0)= n_Y(y_0)[k(y_0): k]$.
A straightforward computation shows that 
$e(\cO_{X,x_0})=e(\cO_{Y, y_0})[k(y_0): k]r^{\dim B}$. 

\end{example}   

\begin{remark} Let $(A, \m)$ be a noetherian local ring of
positive dimension. Let $s(A)$ denote the smallest positive integer $s$
such that 
there exist a reduced closed one-dimensional subscheme 
$C$ in $\Spec A$, and $f\in \mathcal K_C^*(C)$, 
with $s=\ord_{\m}(f)$. In other words, $s(A)$ is the order
of the class of $[\m]$ in the Chow group ${\mathcal A}(\Spec A)$, and in the notation of \ref{def-na}, $s(A)= n(\Spec A, A)$.
It is clear that $s(A) $ divides $ n(A)$, and we note in the example below that 
it may happen that $s(A) < n(A)$.

Let $k:={\mathbb R}$. Consider the projective plane curve $C/k$ given by the equation $y^2z+x^2z+x^3=0$.
Let $X/k$ denote the affine cone over $C$ in ${\mathbb A}^3/k$. Let $c_0 \in C$ denote the singular 
point corresponding to $(0:0:1)$. The preimage of $c_0$ in the normalization of $C$ 
consists of a single point with residue field ${\mathbb C}$. Hence, it follows from \ref{cor.compute-n}  that 
$n_C(c_0)=2$. Let now $x_0:=(0,0,1)\in X(k)$.
Using a desingularization of $x_0$ and \ref{cor.compute-n}, we find that $n_X(x_0)=2$.
We claim that $s(\cO_{X,x_0})=1$. Indeed, the ring $ \cO_{X,x_0}$ is the localization at $(x,y,z-1)$
of the ring $k[x,y,z]/(y^2z+x^2z+x^3)$. The ideal $(x,y)$ defines a closed subscheme of $\Spec \cO_{X,x_0}$
on which $\ord_{x_0}(z-1)=1$. 
\end{remark}

We conclude this section with a further study of the integer $n_X(x_0)$ which will not be used in the remainder of this article. 
\begin{proposition}\label{compute-n2} 
Let $X$ be a noetherian scheme and let $x_0\in X$ be a point with $\dim\cO_{X,x_0}\ge 1$. 
Let $\Gamma_1,\dots, \Gamma_r$ be the irreducible components
of $X$ passing through $x_0$, each endowed with the structure of integral 
subscheme. Then 
$$n_X(x_0) = \gcd \{ n_{\Gamma_i}(x_0) \mid 1\le i\le r\}.$$ 
\end{proposition}

\proof We can suppose that $X$ is local with closed point $x_0$.
Let $U$ be a dense open subset of $X$ and let $i\le r$. 
Then $U\cap\Gamma_i$ is dense in $\Gamma_i$. So 
$n_{\Gamma_i}(x_0)=\ord(f)$ for some 
$f\in \mathcal K_{C}^*(C)$ where $C$ is a reduced curve in 
$\Gamma_i$, with $\max(C)\subseteq U\cap \Gamma_i\subseteq U$. 
Hence. $n_X(x_0)$ divides $ n_{\Gamma_i}(x_0)$ for all $i\le r$. 

For each $i$, fix a dense open subset $U_i$ of 
$\Gamma_i\setminus \cup_{j\ne i}\Gamma_j$ 
such that $n_{\Gamma_i}(x_0)=n(U_i, \Gamma_i)$ (see \ref{lem-n(U,A)}). Let $U$ be  a dense open subset 
of $X$ such that $n_X(x_0)=n(U, X)$. Replacing each $U_i$ by $U_i \cap U$ and $U$ by 
$\cup_i (U_i \cap U)$ if necessary,
we can suppose
that $U:=\cup_i U_i$ satisfies $n_X(x_0)=n(U, X)$. 
Let $C$ be a reduced curve in $X$ such that $\max(C)\subseteq U$
and $n_X(x_0)=\ord(f)$ for some $f\in \mathcal K_C^*(C)$. 
Let $C_i=\Gamma_i\cap C$. Then 
$n_X(x_0)=\sum_i \ord(f_i)$ 
where $f_i$ is the image of $f$ in $\mathcal K_{C_i}^*(C_i)$. As
$\max(C_i)\subseteq U_i$, we have 
$n_{\Gamma_i}(x_0)=n(U_i, \Gamma_i) \mid \ord(f_i)$. Therefore
$n_X(x_0)$ is divisible by the $\gcd$ of the $n_{\Gamma_i}(x_0)$'s.
\qed 

\begin{remark} \label{rem.add}
Let $(A, \m)$ be a noetherian local ring
of dimension $d\ge 1$.
Let $\p_1,\dots, \p_r$ be the minimal prime ideals of $A$, and assume that 
$\dim A/\p_i=d$ if and only if $i \leq s$. Proposition 
\ref{gamma-1} implies that 
$$\gamma(A)=\gcd \{ \length(A_{\p_i})\gamma(A/\p_i) \mid 1\le i\le s\}.$$
Proposition \ref{compute-n2} implies that 
$$n(A) = \gcd \{ n(A/\p_i) \mid 1\le i\le r\}.$$ 
We produced in \ref{thm.na} a class of rings where $n(A) = \gamma(A)$. 
The above formulas indicate that the hypotheses in \ref{thm.na} that $A$ be equidimensional ($s=r$)
and  that $\p_i \in {\rm Reg}(A)$ 
for all $i=1,\dots, r$, are 
optimal.
\end{remark}

Our next proposition strengthens Corollary \ref{mv-degree}
when $n_X(x_0)<\gamma(\cO_{X,x_{0}})$.

\begin{proposition} \label{pro.deltareg2}
Let $X/k$ be a reduced scheme of finite type over a field $k$
and let $x_0 \in X$ be a closed point with $\dim\cO_{X,x_0}\ge 1$. 
Then 
$\delta(X^{\reg}/k)$ divides $ n_X(x_0)\deg_k(x_0)$. 
\end{proposition}
\proof Let $\Gamma_1,\dots, \Gamma_r$ be the irreducible components
of $X$ passing through $x_0$, each endowed with the structure of integral 
subscheme. Since $\delta(X^{\reg}/k) $ divides $ \delta(\Gamma_i^{\reg}/k)$ for all $i=1,\dots, r$,
we find from \ref{compute-n2} that it suffices to prove the statement when $X$ is irreducible.
When $X$ is irreducible, $n_{X}(x_{0}) = \gamma(\cO_{X,x_{0}})$ (\ref{thm.na}),
and the results follows immediately from Corollary \ref{mv-degree}.
\qed

\begin{proposition} \label{compute-n22(a)}  Let $A$ be noetherian excellent local ring (e.g., the local ring at some point of a scheme of finite type over a field) with $\dim(A) \geq 1$. 
If the residue field $k$ of $A$ is algebraically closed, then $n(A)=1$.
\end{proposition}
\proof
Let $X:= \Spec A$, with closed point $x_0$. Proposition \ref{compute-n2} shows that it suffices to prove the statement when $X$ is integral. 
Let $U$ be a dense open subset of $\Spec\cO_{X,x_0}$. 
There exists an integral curve $C$ in $X$ that contains $x_0$ 
and some point of $U$ (Theorem \ref{local-mv}). 
Consider the normalization map $\tilde{C} \to C$. Applying \ref{cor.compute-n} 
to this map and using the fact that $k$ is algebraically closed, we find that 
$n_C(x_0)=1$. By definition, 
there exists $f \in {\mathcal K}_C(C)^*$ such that $\ord(f)=1$. It follows from the definitions
that $n_X(x_0)=1$. \qed

\begin{proposition}\label{compute-n22} 
Let $X/k$ be a scheme of finite type over a field $k$.
Let $x_0 \in X$ be a $k$-rational closed point with $\dim(\cO_{X,x_0})\geq 1$. 
For any extension $F/k$, denote by $x_{0,F}$ 
the unique preimage of $x_0$
under the natural base change map $X_{F} \to X$.
Then there exists a finite extension $F/k$
such that $n_{X_F}(x_{0,F}) = 1$. 
Moreover, any such extension $F/k$ has degree 
divisible by $n_X(x_0)$.
\end{proposition}

\proof
Assume first that $X/k$ is geometrically integral. 
We find then that for any extension $F/k$, $n_{X_F}(x_{0,F}) = \gamma(\cO_{X_F,x_{0,F}})$ (\ref{thm.na}).
Let $\overline{k}/k$ denote an algebraic closure of $k$.
Since \ref{compute-n22(a)} shows that $n_{X_{\overline{k}}}(x_{0,\overline{k}})=1$, 
we can find ideals of definition $Q_1,\dots, Q_r \subset \cO_{X_{\overline{k}},x_{0,\overline{k}}}$
such that $\gcd_i(e(Q_i, \cO_{X_{\overline{k}},x_{0,\overline{k}}}))=1$.
For each $i=1,\dots,r$, choose a system of generators for $Q_i$, 
and denote by $F/k$ the extension of $k$ generated by the coefficients 
of the rational functions needed to define all these generators.
Thus, the generators of $Q_i$ are elements of $\cO_{X_{F},x_{0,F}}$, and generate
in this ring an ideal of definition that we shall denote by $P_i$. 
To conclude the proof, it suffices to note that $e(Q_i, \cO_{X_{\overline{k}},x_{0,\overline{k}}})
= e(P_i, \cO_{X_{F},x_{0,F}})$.

Assume now that $X/k$ is not geometrically integral. First make a finite extension $k'/k$ such that
each irreducible component of $X_{k'}/k'$ is geometrically integral. Let $x_0'$ be the preimage of $x_0$ under $X_{k'} \to X$.
For each irreducible component $\Gamma/k'$ of $X_{k'}/k'$ passing through $x_0'$, we can find using our earlier argument 
an extension $F_{\Gamma}/k'$ such that 
the moving multiplicity of the unique preimage of $x_0'$  under $\Gamma_{F_{\Gamma}} \to \Gamma$ is equal to $1$. 
Let $F/k$ 
be one of the extensions $F_{\Gamma}$.  
It follows then from Proposition \ref{compute-n2} that the moving multiplicity of the preimage of $x_0'$  under $X_{F} \to X_{k'}$ is equal to $1$.  

Let now $F/k$ such that $n_{X_F}(x_{0,F})=1$. Consider the finite flat 
base change $f:X_F \to X$, with $f(x_{0,F}) = x_0$ and $k(x_0)=k$.
Theorem \ref{compute-n}  (c) shows that $ n_X(x_0) $ divides $[F:k] n_{X_F}(x_{0,F}) $.
It follows that the degree of the extension $F/k$ is divisible by $ n_X(x_0) $.
\qed 

 \medskip
 As the last paragraph of the above proof indicates, when $n_X(x_0)>1$ and $F/k$ is separable and non-trivial, the \'etale morphism $f:X_F \to X$, with $f(x_{0,F}) = x_0$,
is such that $n_{X_F}(x_{0,F}) < n_X(x_0)$.  This example shows that the hypothesis on the residue 
fields of the points $x_0$ and $y_0$ is necessary in (b) of our next proposition. 

\begin{proposition} \label{pro.smooth} 
Let $X$ be a 
noetherian scheme, and let 
$f : Y\to X$ be a smooth morphism.  
Let $x_0\in X$ with $\dim \cO_{X,x_0} \geq 1$ and  
let $y_0 $ be any point in the fiber $ f^{-1}(x_0)$. 
Then
\begin{enumerate}[\rm (a)]
\item  $n_Y(y_0)$ divides $ n_X(x_0)$. 
\item Assume that $X$ is universally catenary.  If  
$k(x_0) \to k(y_0)$ is an isomorphism, then $n_Y(y_0)=n_X(x_0)$. 
\item Assume that $X$ is universally catenary. Let 
$g : Z\to X$ be any smooth morphism with   
  $z_0 \in g^{-1}(x_0)$ such that $k(z_0)$ is isomorphic over $k(x_0)$ to $k(y_0)$.  
Then $n_Z(z_0)=n_Y(y_0)$.  
\end{enumerate}
\end{proposition}

\proof (a) 
Since by definition $n_X(x_0)= n(\cO_{X,x_0})$, we can assume that $X$ is local 
with closed point $x_0$ by replacing $f$, if necessary, 
by $Y \times_X \Spec \cO_{X,x_0} \to \Spec \cO_{X,x_0}$.
We can also assume that $Y\to X$ is of finite type.  

Since $f$ is smooth at $y_0$, replacing $Y$ by an open neighborhood of 
$y_0$ in $Y$ if necessary, we can suppose that $f$ factors as the composition of an 
\'etale map $Y \to {\mathbb A}^d_X$ followed by the canonical projection 
${\mathbb A}^d_X \to X$. 
Thus it is enough to prove (a) separately for \'etale morphisms
and for the projections $\mathbb A^d_X \to X$ for all $d\ge 1$. 

We now make one further reduction, to the case where $X$ is local integral of dimension $1$.
This paragraph and the next only assume that $f$ is flat.
Consider the natural morphism $i: \Spec \cO_{Y,y_0} \to Y$. 
Using Lemma \ref{lem-n(U,A)}, we can find a dense open set $U'$ of $Y$  
such that $n(\cO_{Y,y_0})=n(i^{-1}(U'), \cO_{Y,y_0})$.
Note as in \ref{lem-n(U,A)} that for any dense open subset $U_1 \subset U'$, we have
$n(\cO_{Y,y_0})=n(i^{-1}(U_1), \cO_{Y,y_0})$. We now choose a dense open subset 
$U \subseteq U'$ such that $U$ is \emph{fiberwise empty or dense in $Y$ over $X$}, that is, such that 
for all $x \in X$, $U_x $ is either empty  or   dense in $Y_x$. Indeed, first 
use \cite{EGA}, IV, 9.7.8 (applied to each irreducible component of $X$), to 
find a dense open subset $X_1$ of $X$
such that for each $x \in X_1$, say belonging to an irreducible component $\overline{ \{\eta\}}$ of $X$, the number of geometric irreducible components of $Y_x$ is equal to the number of geometric irreducible components of $Y_{\eta}$.
Note now that $f$ is flat, so the image $f(U')$ of $U'$ in $X$ is open. Apply then again \cite{EGA}, IV, 9.7.8,
to the morphism $f_{\mid U'}: U' \to f(U')$ and find an appropriate dense open subset $V_1$ of $f(U')$
with similar conditions on the fibers. It follows that 
$W:=V_1 \cap X_1$ is a dense open subset of $X$, and we take $U:=f_{\mid U'}^{-1}(V_1 \cap X_1) \subset U'$. Note that $f(U)=W$.

Recall that $n(U, \cO_{Y,y_0})= n_Y(y_0)$, and that $n(W, \cO_{X,x_0})$ divides $n_X(x_0)$.  
To prove that $n_Y(y_0)$ divides $n_X(x_0)$, 
it suffices to prove that $n(U, \cO_{Y,y_0})$ divides $n(W, \cO_{X,x_0})$.
Thus, it is enough, given any integral
curve $C$ in $X$ containing $x_0$ and whose generic point $\eta$ belongs to $ W$, and given any non-zero 
regular function $g\in \cO(C)$, to find a
one dimensional closed subscheme $D$ of $\Spec \cO_{Y,y_0}$, 
contained in $Y\times_X C$ with maximal points in $U$,  
and an element $g_1\in \mathcal K_D^*(D)$ such that 
$\ord_{y_0}(g_1)$ divides $\ord_{x_0}(g)$. 
As $U \times_X \eta$ is dense in $Y \times_X \eta$, 
we find that $U\times_X C$ is a dense open subscheme of $Y\times_X C$. 
Therefore, to prove that $n_Y(y_0)$ divides $n_X(x_0)$, 
it is enough to prove that 
$n_{Y\times_X C}(y_0)$ divides $n_{C}(x_0)$; indeed, then there exists 
such a $D$ and $g_1$ with $\ord_{y_0}(g_1)= n_{Y\times_X C}(y_0)$, and this latter integer divides 
$n_{C}(x_0)$, which by definition itself divides $\ord_{x_0}(g)$.

Returning to the case where $f$ is smooth as in the proposition, 
we find that it suffices to prove (a) when $X$ is local integral of dimension $1$ with generic point $\eta$.
Now suppose that $f$ is \'etale. We keep the notation $U$ and $W$ of the previous paragraphs.
Suppose given a
nonzero regular function $g$ on   $X$ ($X$ can now play the role of $C$). 
Let $D:= \Spec \cO_{Y,y_0}$.
As $U_{\eta}$ is dense in the discrete set $Y_{\eta}$, it is equal to $Y_{\eta}$. Hence, 
the maximal points of $D$ belong to $U$. Let $g_1:=f^*(g)$. 
Then $\ord_{y_0}(g_1)=\ord_{x_0}(g)$ because $f$ is \'etale, and (a) is proved in this case.  

Suppose now that $f$ is a projection morphism $\mathbb A^d_X\to X$. An easy
induction in $d$ reduces the problem to the case $d=1$. Consider now  
$f:Y:=\mathbb A^1_X \to X$. Assume first that $y_0$ is a closed point of $Y$.
Write $X=\Spec R$, with $R$ a local domain with maximal ideal $\m $. Fix $g \in R \setminus \{0\}$.
Let $F(T)\in R[T]$ be a monic
polynomial which lifts a generator of the maximal ideal of 
$k(x_0)[T]$ corresponding to $y_0\in Y_{x_0}$. Let $S:= R[T]/(F(T))$, and $D:=\Spec S$. 
Then $S$ is finite and flat over $R$, with $S/\m S = k(y_0)$ a field, so that $S$ is local.
Let $g_1=g\in S$. The element $g_1$ is regular  by flatness and we have 
$$\ord_{y_0} (g_1)=\length_{S} (S/(g_1)) =
\length_{S} (S\otimes_R R/gR)= 
\length_{R} (R/gR)=\ord_{x_0}(g),$$   
where the third equality holds by flatness and because 
$S/\m S=k(y_0)$ is a field (see \cite{FulBook}, A.4.1,
or \cite{Liubook}, 
Exercise 7.1.8(b)). To conclude the proof, we need to further specify the element $F(T)$, so that $D$ is such that
$D_\eta\subset U$. For this we modify a given $F(T)$ by an element $\lambda\in \m$
as follows.  Write $Y_{\eta}\setminus U_{\eta}$
as $V(\phi(T))$ for some non-zero $\phi(T)\in \Frac(R)[T]$. Since $\m $
is infinite, 
we can find $\lambda\in \m$ such that
$\gcd(F(T)+\lambda, \phi(T))=1$ in $\Frac(R)[T]$. Then we take the curve $D$ defined
by $ F(T)+\lambda$. 

Let us consider now the case where $y_0$ is the generic point of the fiber of $ \mathbb A^1_X \to X$
above $x_0$. As before, let $R:= \cO_{X,x_0}$, with maximal ideal $\m$. Then $\cO_{Y,y_0}= R[T]_{\m R[T]}=:R(T)$. 
Let $I$ be an ideal in $R$. Then $I$ is $\m$-primary if and only if $IR(T)$ is $\m R(T)$-primary, 
and when $I$ is $\m$-primary, then $e(I,R) = e(IR(T),\m R(T))$ (as stated for instance in \cite{SH}, 8.4.2).
It follows that $\gamma(R(T))$ divides $\gamma(R)$. 
Since $n(R(T))$ divides $\gamma(R(T))$, we find that $n(R(T))$ divides $n(R)$ when $n(R) = \gamma(R)$. 
When $R$ is a local   domain of dimension $1$, the hypotheses of \ref{thm.na} are satisfied 
and $n(R) = \gamma(R)$.

(b) When $X$ is universally catenary, we can apply  Theorem \ref{compute-n} (c) (since $f$ is open, the hypothesis on generic points is satisfied) 
and obtain that  $n_X(x_0)$ divides $n_Y(y_0)[k(y_0):k(x_0)]$. This
gives the divisibility relation $ n_X(x_0) \mid n_Y(y_0)$ when $[k(y_0):k(x_0)]=1$.

(c) Consider the cartesian diagram
$$
\xymatrix{
W \ar[r]^{\beta}  \ar[d]^{\alpha} & Y \ar[d] \\
Z \ar[r] & X. \\
}
$$
Let $w$ denote a   point of $W:= Z \times_X Y $ above $z_0$ and $y_0$ inducing isomorphisms $k(z_0) \to k(w)$ and $k(y_0) \to  k(w)$.
The morphisms $\alpha$ and $\beta$ are smooth,
and it follows from (b) applied to $\alpha$ and $\beta$ that $ n_Z(z_0)= n_W(w)= n_Y(y_0)$. 
\qed

\begin{example} \label{rem.non-uc2}
We show in this example that there exist a local ring $\cO$ and an ind-\'etale local extension $\cO^h$ 
such that $\Spec \cO^h \to \Spec \cO$ induces an isomorphism on closed points, and 
such that $n(\cO^h) < n(\cO)$.
In our example, the local ring $\cO$ is not universally catenary, 
while the ring $\cO^h$ is the henselization of $\cO$, and is universally catenary.

Consider  again the finite birational morphism $\pi: Y \to X$ described
in Example \ref{non-uc}. The scheme $Y$ is regular, and $X$  is catenary but not universally catenary.
The preimage of the point $x_0 \in X$ consists in the two regular points $y_0$ and $y_1$ in $Y$,
with $[k(y_0):k(x_0)]=d$, and $[k(y_1):k(x_0)]=1$.
It follows from the discussion in \ref{non-uc} that $n_X(x_0)$ is divisible by $d$.
Let $\cO := \cO_{X,x_0}$, and consider the henselization $\cO^h$ of $\cO_{X,x_0}$.
It is known that a noetherian local  Henselian ring of dimension $2$ is universally catenary (see, e.g, \cite{Rat}, 2.23 (i), with for
instance, \cite{SH}, B.5.1).
We let $x_0'$ denote the closed point of $\Spec \cO^h$.
We have the following diagram with cartesian squares:
$$
\xymatrix{
Y_{\cO^h}   \ar[d]^{\pi'} \ar[r] & Y_{\cO}   \ar[d] \ar[r] & Y \ar[d]^{\pi} \\
\Spec \cO^h \ar[r] &  \Spec \cO \ar[r] & X. \\
}
$$
We use \ref{compute-n}, (b) and (c), on $\pi'$ to conclude that
$$\gcd \{ n_{Y_{\cO^h}}(y)[k(y) : k(x_0')] \mid y \in \pi'^{-1}(x_0'), y 
\ {\rm  closed} \}  =n(\cO^h) =1.$$  
Thus, when $d>1$,  we find that  $1=n(\cO^h) < n(\cO)$.

\end{example} 
\begin{remark}
It is easy to construct examples of flat morphisms $f:Y \to X$ with $y_0 \in Y$
such that $n_{Y}(y_0) > n_X(f(y_0)) $. Indeed, start with an integral scheme $Y$ of finite type 
over a field $k$ with a closed point $y_0 \in Y$ such that $n_{Y}(y_0) > 1$.
Choose any non-constant function $f:Y \to {\mathbb P}^1_k$. This morphism is then flat, 
and since ${\mathbb P}^1_k$ is regular, $n_{{\mathbb P}^1_k}(f(y_0))=1$.
\end{remark} 

\end{section}

\begin{section}{Index of varieties over a discrete valuation field} 
\label{valuation-field}

\begin{emp} \label{defX/S}
Let $K$ be the field of fractions of a discrete valuation 
ring ${\cO_K}$, with maximal ideal $(\pi)$ and residue field $k$. Let $S:=\Spec(\cO_K)$. Let $\X$ 
be an integral scheme, and let $f:\X \to S$ be a flat, separated, surjective morphism of finite type. 

Since $f $ is flat, $\dv(\pi)$ is a Cartier divisor on $\X$,
and we denote its associated cycle by $[\dv(\pi)] = \sum_{i=1}^n r_i \Gamma_i$.
Each  
$\Gamma_i$ is an integral variety over 
$k$, of multiplicity $r_i$ in $\X_k$.
The generic fiber of $\X/S$ is denoted by $X/K$. 

For any irreducible $1$-cycle $C$ (endowed with the reduced induced 
structure) on $\X$ and for any Cartier divisor $D$ on $\X$ whose support does
not contain $C$, the restriction $D|_C$ is again a Cartier divisor on $C$.
The associated cycle $[D|_C]$ is supported on the special fiber of $\X\to S$. 
Writing $[D|_C]= \sum_{x {\rm \ closed}} (D.C)_x[x]$, then its degree 
over $k$ is $\deg_k [D|_C]= \sum_x (D.C)_x \deg_k(x)$. 
Note that if $D$ is effective, 
then for any closed point $x\in C\cap \Supp D$, 
$(D.C)_x=\length(\cO_{C, x}/\cO_{C}(-D)_x)\ge 1$. 

Assume that $C \to S$ is finite.   
Then $C$ is the closure in $\X$ of a closed point $P \in X$, and we have 
$$\deg_K(P) = \deg_k(\dv(\pi)|_C).$$
Assume now that $\X$ is locally factorial (e.g., regular). Then the Weil divisors $\Gamma_i$
are Cartier divisors (again denoted by $\Gamma_i$) and 
we can write that $\dv(\pi)= \sum_{i=1}^n r_i \Gamma_i$ as Cartier divisors.
It follows that
\begin{equation}\label{degree-P}
\deg_K (P)=
\sum_{x\in \X_k\cap C}\left(\sum_{\Gamma_i\ni x} r_i(\Gamma_i.C)_x\deg_k (x)\right). 
\end{equation}
In particular, 
$$\gcd_i\{r_i \delta(\Gamma_{i}/k)\} \text{ \ \rm divides \ }  
   \overline{\delta}(X/K),$$
where $\overline{\delta}(X/K)$ denotes the greatest common divisor 
of the integers $\deg_K(P)$, with $P \in X$ closed, and whose closure in $\X$
is finite over $S$. This statement is sharpened 
in our next theorem.

When $f: \X \to S$ is proper, the closure of 
any closed point $P \in X$ is finite and flat over $S$ and, thus, in this case, 
$\gcd_i\{r_i \delta(\Gamma_{i}/k)\} $ divides 
$\delta(X/K)=\overline{\delta}(X/K)$. 
\end{emp}

\begin{theorem} \label{dx-ds}  
Let $f:\X \to S$ be as above, 
with $\X$ regular.
Let $X/K$ denote the  generic fiber of $\X/S$. 
\begin{enumerate}[\rm (a)]
\item 
Then $\gcd_i\{r_i \delta(\Gamma^{\reg}_{i}/k)\}$ divides $\deltadelta(X/K)$.
\item  When $\cO_K$ is Henselian, then 
$\deltadelta(X/K) = \gcd_i\{r_i \delta(\Gamma^{\reg}_{i}/k)\}$.
\end{enumerate}
\end{theorem}

\proof For ease of notation, we will write 
$\gcd(\X_k):= \gcd_i\{r_i \delta(\Gamma^{\reg}_{i}/k)\}$.

(a) Let $P$ be a closed point of $X$
whose closure in $\X$ is finite over $S$, and
let us show that $\gcd(\X_k)$ divides $\deg_K(P)$. 
If
$\overline{\{ P \}}\cap \Gamma_i\subseteq \Gamma_i^{\reg}$ for 
all $i\le n$,  then Formula (\ref{degree-P}) above shows that $\gcd(\X_k) $ divides $ \deg_K(P)$.  
In general, though,  $\ol{\{P\}}$ may intersect the singular
locus of some $\Gamma_i$. If that is the case, then we can end the proof in two different ways.

The first method relies on Theorem \ref{mv-slocal}, and assumes in addition that $\X$ is $\FA$
when $\cO_K$ is not Henselian.
Indeed, \ref{mv-slocal}
shows that there exists an affine open subset $V$ of $\X$ which contains 
the $1$-cycle  
$\ol{\{ P \}}$ and a  $1$-cycle $C$  rationally equivalent to $\ol{\{ P \}}$ in  $V$, and whose support
is proper over $S$ and does
not intersect the singular locus $F$ of $(\X_k)_\mathrm{red}$. 
Then $P$ is rationally equivalent on $V_K$ to $C_{\mid V_K}$, whose support is a union
of closed points of $X$. We claim that  $\deg_K(P)= \deg_K C_{\mid X}$.
Indeed, since $V$ is affine,
we can consider an open embedding $V \to \Y$ over $S$ where $\Y/S$ is projective. 
Theorem \ref{mv-slocal} shows that $\ol{\{ P \}}$ and $C$ are closed and rationally 
equivalent in $\Y$. Then 
$\deg_K(P)= \deg_K C_{\mid \Y_K}=\deg_K C_{\mid V_K}=\deg_K C_{\mid X}$. The above discussion shows that
each point in $\Supp C_{\mid X}$ has degree
divisible by $\gcd(\X_k)$, so that   $\gcd(\X_k) $ divides $ \deg_K(P)$, as desired. 

The second method relies in the end on Theorem \ref{mvsimple} and its corollary \ref{mv-degree}.
It consists in a succession of blowing-ups, starting with the blowing-up of specializations of $P$ in $\X$ 
 as in Lemma \ref{ds-blowup} (3), to produce 
a new regular model $\Y$, where the specializations of $P$ belong to 
the regular locus of $(\Y_k)_{\mathrm{red}}$. 
Our initial  discussion above implies that
$\gcd(\Y_k)$ divides $\deg_K(P)$. Then Lemma \ref{ds-blowup} (2) 
shows that  $\gcd(\Y_k)= \gcd(\X_k)$, and (a) is proved. 

\medskip

(b) Let $\Gamma_i^0 := \Gamma_i^{\reg} \setminus \cup_{j\neq i}\Gamma_j$.
Proposition \ref{pro.deltareg} shows that $\delta(\Gamma_i^{\reg}/k) = \delta(\Gamma_i^0/k)$. 
We use then Proposition \ref{lift} (3) on each closed point of $\Gamma_i^0$ to find that $\deltadelta(X/K) $ divides
$ \gcd(\X_k)$. 
\qed 

\begin{lemma} \label{ds-blowup} 
Let $f:\X \to S$ be as in {\rm \ref{defX/S}}, 
and assume $\X$ regular. Let  $d:=\dim X$. Let $x_0\in \X_k$ be a closed
point and let $\Gamma_1,\dots, \Gamma_s$ denote the irreducible components
of $\X_k$ which contain $x_0$. Let $e_i$ denote the 
Hilbert-Samuel multiplicity of  $x_0$ on $\Gamma_i$.   Consider the blowing-up 
$\tilde{\X}\to \X$ of $\X$ along the  reduced closed subscheme 
$\{ x_0 \}$.
\begin{enumerate}[{\rm (1)}] 
\item The scheme $\tilde{\X}$ is regular, and the  exceptional divisor $E_1$ in $\tilde{\X}$ is isomorphic to
$\mathbb P^d_{k(x_0)}$. The multiplicity 
$r(E_1)$ of $E_1$ in $\tilde{\X}_k$ is $\sum_{i=1}^s r_ie_i$. 
\item We have $\gcd(\tilde{\X}_k)=\gcd(\X_k)$. 
\item Let $P$ be a closed point of $X$. Then there exists a finite
sequence 
$$ \X_m \to \X_{m-1} \to \cdots \to \X_0=\X $$
such that each morphism $\X_i \to \X_{i-1}$ is the blowing-up of a 
 closed point in the special fiber, and such that the closure of 
$P$ in $\X_m$ intersects  $(\X_m)_k$ only in regular points of 
$((\X_m)_k)_{\red}$.
\end{enumerate}
\end{lemma}

\proof (1) As we blow-up a regular scheme along the regular center
$\Spec k(x_0)$, we have $\tilde{\X}$  regular and $E_1\simeq \mathbb P^d_{k(x_0)}$ (\cite{Liubook}, 8.1.19). 
In the regular local ring ${\mathcal O}_{\X,x_0}$, 
factor a uniformizing element $\pi$ of $\cO_K$ 
as $\pi= u g_1^{r_1}\cdot \ldots \cdot g_s^{r_s}$, 
where $g_i$ is a local equation in $\X$ of the component $\Gamma_i$ at 
$x_0$, and $u$ is a unit. 
It is not hard to see that the Hilbert-Samuel multiplicity 
of $x_0$  on  $\Gamma_i$ is the positive 
integer $e_i$ such that $g_i \in ({\m}_{\X,x_0})^{e_i}
\setminus ({\m}_{\X,x_0})^{e_i+1}$. Since the associated graded ring 
$\oplus_{q\ge 0} ({\m}_{\X,x_0}^{q}/{\m}^{q+1}_{\X,x_0})$ 
is a polynomial ring over $k(x_0)$, we find that $\pi \in ({\m}_{\X,x_0})^{\sum_i r_ie_i}
\setminus ({\m}_{\X,x_0})^{\sum_i r_ie_i+1}$. Let $\xi$ denote the generic 
point of $E_1$. Since $\tilde{\X}$  is regular, the local ring $\cO_{\tilde{\X},\xi}$
is a discrete valuation ring with normalized valuation $v_{\xi}$, 
and the multiplicity of $E_1$ in $\tilde{\X}_k$ is  $v_{\xi}(\pi)$. 
We leave it to the reader to check that 
if $f \in ({\m}_{\X,x_0})^{r}
\setminus ({\m}_{\X,x_0})^{r+1}$, then $v_{\xi}(f) =r$. 

(2) Let 
 $\tilde{\Gamma}_i$ denote the strict transform of $\Gamma_i$ in $\tilde{\X}$.
The blowing-up  $\tilde{\X} \to \X$ restricts to a morphism 
$\tilde{\Gamma}_i \to \Gamma_i$ which is nothing 
but the blowing-up of $\Gamma_i$ along $x_0$. Hence, 
the varieties $\tilde{\Gamma}_i/k$ and $\Gamma_i/k$ are birational
and $\delta(\tilde{\Gamma}_i^{\reg}/k) = \delta(\Gamma_i^{\reg}/k)$ (\ref{pro.deltareg}).
Recall that 
$$
\gcd(\tilde{\X}_k)= \gcd_{1\le i\le n}\{ r_i\delta(\tilde{\Gamma}_i^{\reg}/k), 
r(E_1)\delta(E_1/k)\}. $$ 
We thus find that  $\gcd(\tilde{\X}_k) $ divides $ \gcd(\X_k)$. 
Clearly, $\delta(E_1/k)= [k(x_0):k] = \deg_k (x_0)$.
Corollary \ref{mv-degree} implies that for 
every $\Gamma_i$ passing through $x_0$, $e_i\deg_k(x_0)$ is the degree of 
some $0$-cycle supported in $\Gamma_i^{\reg}$ and, thus, 
$e_i\deg_k(x_0)$ is divisible by $\delta(\Gamma_i^{\reg}/k)$.
Therefore $\gcd(\X_k) $ divides $\gcd(\tilde{\X}_k)$.

(3) Let $P$ be a closed point on $X$. We describe below how to obtain 
a new regular model using a sequence of blowing-ups along closed points 
$$ \X_m \to \X_{m-1} \to \cdots \to \X_1\to \X$$
such that the specializations of $P$ in $\X_m$ are regular points
in $((\X_{m})_k)_{\red}$. 

Suppose $\ol{\{P\}}$ meets the singular locus $F$ of $(\X_k)_{\red}$.
Let $\Y\to \X$ be the blowing-up of $\X$ along $\ol{\{P\}}\cap F$. 
Then $\Y\to \X$ is also obtained by successively blowing-up points
of $\ol{\{P\}}\cap F$. If $P$ specializes only to regular points of
$\Y_k$, then the process stops. Otherwise, 
let $y$ be a specialization of $P$ in $\Y_k$ belonging to the singular 
locus of $(\Y_k)_{\red}$. 
Then $y$ belongs to an exceptional divisor $E_1$, of multiplicity $r(E_1)$, 
and to at least one strict 
transform $\tilde{\Gamma}_i$. Using these two facts and 
Formula (\ref{degree-P}), we find that 
$$\deg_K(P)> r(E_1).$$ 
Consider the blowing-up $\mathcal Z\to \Y$ of $\Y$ along 
$\ol{\{ P \}}$ intersected with the singular locus of $(\Y_k)_{\red}$. 
The exceptional divisor $E_2$ above $y$ has multiplicity $r(E_2)$ in 
$\mathcal Z_k$, and since $y\in E_1\cap \tilde{\Gamma}_i$,
Part (1) of this lemma implies that 
$r(E_2) > r(E_1)$.  Repeating the above argument on $E_2$ shows that
$\deg_K(P)> r(E_2)$. Therefore, the process must stop after
at most $1+\deg_K (P)$ steps. 
\qed 

\medskip 
Let $F$ be any field, and let $W/F$ be a scheme of finite type. Recall that a  
closed point $P\in W$ is called {\it separable} 
if the residue field extension $F(P)/F$ is a separable extension.

\begin{proposition} \label{lift} 
Let $f:\X \to S$ be as in {\rm \ref{defX/S}}. 
Let $x_0\in \X_k$ be a closed point, regular in $\X$. 
Let $\Gamma_1,\dots,\Gamma_s$ be the irreducible components of $\X_k$ 
passing through $x_0$. Let 
$e_i:=e(\cO_{\Gamma_i, x_0})$ denote the 
 Hilbert-Samuel multiplicity of $x_0 $ on  $\Gamma_i$. 
Then 
\begin{enumerate}[{\rm (1)}]
\item The Hilbert-Samuel multiplicity $e(\cO_{\X_k, x_0})$ of $x_0$ on $\X_k$  is equal to 
$\sum_{i=1}^s r_ie_i$. 
\item Any closed point $P\in X$ such that $x_0\in \overline{\{ P\}}$ has
$\deg_K(P)\geq e(\cO_{\X_k, x_0})\deg_k (x_0)$. 
\item Assume that $\cO_K$ is Henselian. If $k$ is infinite,
or if $x_0$ is a regular point of $(\X_k)_{\red}$, 
then there exists a closed 
point $P\in X$ such that $x_0\in \overline{\{ P\}}$ and 
$$\deg_K(P)=e(\cO_{\X_k, x_0}) \deg_k (x_0).$$ 
If $k$ is finite, then there exists a $0$-cycle on $X$ of degree 
$e(\cO_{\X_k, x_0}) \deg_k (x_0)$ and  such that  each closed point in  
its support specializes to 
$x_0$ in $\X$. 

If $X/K$ is generically smooth, then the point $P$ and the support of the $0$-cycle
can be chosen to be separable over $K$.
\end{enumerate}
\end{proposition}

\proof (1) This is the same computation as in Lemma \ref{ds-blowup} (1),
or use \ref{e-irr}. 

(2) Let $\tilde{\X}\to \X$ be the blowing-up along $x_0$ with
exceptional divisor $E_1$. We saw in \ref{ds-blowup} (1) that 
$E_1$ is isomorphic to $\mathbb P^d_{k(x_0)}$ and has multiplicity 
$e(\cO_{\X_k, x_0})$ in $\tilde{\X}_k$. 
As $P$ has a specialization in $\tilde{\X}$ belonging to $E_1$, 
Formula (\ref{degree-P}) before Theorem \ref{dx-ds} shows that 
$\deg_K(P)\geq e(\cO_{\X_k, x_0})\deg_k (x_0)$. 

(3) Let us suppose first that $x_0$ is a regular point of 
$(\X_k)_{\red}$. 
Since $\X^{\reg}$ is open in $\X$ (\cite{EGA} IV.6.12.6 (ii)), 
we may without loss of generality assume   that $\X= \Spec A$ is affine, irreducible, and regular,
with irreducible special fiber. Since $\cO_{X,x_0}$ is factorial, we may if necessary replace $\X$ by an open dense subset and assume 
that the uniformizer $\pi$ of $\cO_{K}$ factors as $\pi= u t^e$ in $A$,
with $t\in A$, $u\in A^*$ and $e:=e(\cO_{\X_k, x_0})$. 
By hypothesis, there exists a system of generators $\{f_1,\dots, f_d\}$ of the maximal ideal $\m$ of $A/(t)$ corresponding to $x_0$
with $d= \dim(A/(t))$. 

For $i=1,\dots, d$, let  $g_i \in A$ be any lift of $f_i$. Let $T_0$ be the closed subscheme of 
$\X$ defined by $T_0:=\Spec A/(g_1,\dots,g_d)$. Consider the induced morphism $\varphi: T_0 \to S$. Clearly, $\varphi^{-1}((\pi))= \{x_0\}$. 
Let $T_1 \subseteq T_0$ denote the open subset (see \cite{EGA}, IV.13.1.4) consisting of all the points of $T_0$ where $\varphi$ is quasi-finite. It is clear that $x_0 \in T_1$, and $x_0$ is in fact a regular point of $T_1$ since 
$A/(g_1,\dots,g_d,t)$ is a field. 
Using \cite{EGA} IV.6.12.6 (ii) again, we find that the regular locus of $T_1$ is open, and we can thus if necessary 
replace $T_1$ by a open subset containing $x_0$ and assume that $T_1 $ is regular. Since $\cO_K$ is Henselian, 
we may use \cite{BLR}, 2.3, Proposition 4 (e), and obtain that there exists an open neighborhood $T$ of $x_0$ in $T_1$
such that $T \to S$ is finite.  Since both $T$ and $S$ are regular and $T\to S$ is finite, we find that $T \to S$ is also flat.
The generic point $P$ of $T= \overline{\{P\}} $ has thus degree $e \deg_k(k(x_0))$ over $K$, as desired.

Let us now assume in addition that $X/K$ is generically smooth. We claim that   
we can find a lift $P \in X$ of $x_0$ with $K(P)/K$ separable of degree $e \deg_k(k(x_0))$ over $K$.
Keep the assumptions in the first paragraph of the proof of (3) above. 
If $d=0$,  then the generic 
smoothness of $X/K$ implies that the fraction field of $A$ is separable over $K$, and our claim is true.

Assume now that $d \geq 1$, and let us show that the
lifts $g_i$ can be chosen such that the generic point of the 
associated $T$ is separable over $K$. To start, we 
show that there exists a lift $g_1$ of $f_1$ such that the closed irreducible 
subscheme $V(g_1)$ of $\X$ 
is flat of finite type over $S$, with generic 
fiber generically smooth. We then conclude the proof  by induction on the dimension.

Note that there always 
exists $h \in A$ such that the differential $dh$ in $ \Omega^1_{A \otimes K/K}$ 
is not zero. Start with any lift $g \in A$ of $f_1$. If the differential $dg$ is zero, replace $g$ by the lift $g+\pi h$
and assume now that $dg \neq 0$.  Let $Z$ 
be the proper closed subset of $\X$ where $\X \to S$ is not smooth. 
Then on $\X \setminus Z$, the sheaf of differentials is locally free, and we let $Z'$ be the zero locus of $dg$ in $\X \setminus Z$.
By construction, $Z'$ is closed in $\X \setminus Z$, and
$Z\cup Z'$ 
contains only finitely many irreducible closed subsets of codimension $1$ in $\X$. 

The ideals $(g + \pi^s)$ of $A_\m$, $s \in {\mathbb N}$, are infinitely  many pairwise distinct prime ideals (recall that
by construction $(g)$ is a prime ideal of $A_\m$; and the maximal ideal 
of $A_\m$ is also generated by $t,g+\pi^s,g_2,\dots, g_d)$. 
Therefore, we can choose 
$g_1:=g+\pi^s$ for some $s$ such that in $\X$, $V(g_1)$ is not contained in $Z \cup Z'$.
Since $V(g_1)$ is not contained in $Z$, we find that $V(g_1)$ intersects the smooth locus
of $\X \to S$. The generic fiber of the morphism $V(g_1) \to S$ is then smooth at a point outside $Z \cup Z'$
since $d(g_1)$ is not zero at such a point. We conclude the proof of the existence
of $P$ with $K(P)/K$ separable when $X/K$ is generically smooth by induction on the dimension.

In the general case where $x_0$ need not be regular in $(\X_k)_{red}$, consider the blowing-up $\tilde{\X}\to\X$ 
of $x_0$ as in (2). If $k$ is infinite, then there exists a $k(x_0)$-rational 
point $x_1$ in the interior of $E_1$, and so regular in $(\tilde{\X}_k)_{\red}$. 
By the above, there exists a closed point 
$P$ of $X$ such that $\deg_K(P)=e(\cO_{\X_k, x_0})\deg_k(x_0)$
and which specializes to $x_1$ in $\tilde{\X}$. Then 
$P$ specializes to $x_0$ in $\X$, as desired.

If $k$ is finite, 
then $\delta(Y/k)=1$  for any geometrically irreducible 
algebraic variety $Y$ over $k$ (\cite{L-W}, Corollary 3. See also 
\cite{CM}, 3.11).  
Therefore, there exists a $0$-cycle supported in the
interior of $E_1$, of degree $1$ over $k(x_0)$. 
Then lift this
$0$-cycle to $X$ as above. 
\qed

\smallskip
Variations on the statement of \ref{lift} (3) when $x_0$ is a regular point of $(\X_k)_{\red}$
can be found in 
\cite{C-S}, Lemme 2.3,  or \cite{BLR}, Corollary 9.1/9.
 When we relax the hypothesis that $\X$ is regular in 
Theorem \ref{dx-ds} (b), we obtain:

\begin{proposition}  Let $\cO_K$ be
Henselian. Let $f:\X \to S$ be as in {\rm \ref{defX/S}}, with $\X$ regular in codimension one.
Then 
$\delta(X^{\reg}/K) $ divides $ \gcd_i \{ r_i\delta(\Gamma_i^{\reg}/k)\}$.
\end{proposition}

\proof Indeed, $\X^{\reg}$ is open in $\X$ (\cite{EGA} IV.6.12.6 (ii)).
Since $\X$ is regular in codimension one, $\X \setminus \X^{\reg}$ is of codimension $\ge 2$. 
Thus,  $\Gamma_i^{\reg}\cap \X^{\reg}$ is not empty, 
and $\delta(\Gamma_i^{\reg}\cap \X^{\reg}/k) =\delta(\Gamma_i^{\reg}/k)$ 
for all $i$ (\ref{pro.deltareg}). Then lift closed 
points of $\Gamma_i^{\reg}\cap \X^{\reg}$ to $X^{\reg}$ as in 
\ref{lift} (3).  
\qed

\begin{remark} \label{rem.multiplicity2} 
Let $\Gamma/k$ be an integral normal algebraic variety. Let $k'$ be the 
algebraic closure of $k$ in the field of rational functions $k(\Gamma)$.
Let $e(\Gamma/k')$ be the geometric multiplicity of $\Gamma$ as $k'$-scheme 
(\cite{BLR}, 9.1/3). We claim that when $\Gamma$ is regular,   
$$e(\Gamma/k')\text{ divides } \delta(\Gamma/k').$$
Indeed, let $L/k'$ be a separable closure of $k'$. Then 
$e(\Gamma\times_{k'}L/L)= e(\Gamma/k')$, and 
$\delta(\Gamma\times_{k'}L/L) $ divides $ \delta(\Gamma/k')$. 
Let $x_0 \in \Gamma\times_{k'}L$ be a Cohen-Macaulay closed point, 

and let $f_1,\dots,f_d$ be a maximal regular sequence in  
$\m_{x_0}\mathcal O_{\Gamma\times_{k'}L,x_0}$ of length 
$\mu:=\length(\mathcal O_{\Gamma\times_{k'}L, x_0}/(f_1,\dots,f_d))$. 
By \cite{BLR}, 9.1/7 (b), 
$e(\Gamma/k')$ divides $\mu \deg_L(x_0)$.
Hence, when $\Gamma$ is regular, so is  $\Gamma\times_{k'}L$ and 
$e(\Gamma/k')$ divides $\delta(\Gamma/k')$. Note that $e(\Gamma/k')=e(\Gamma/k)$ if 
$k'/k$ is separable.

 Let now $\X \to S$ be as in Theorem \ref{dx-ds}, with $f$ proper and $\X $ regular.
 Let $k_i$ denote the algebraic closure 
of $k$ in the function field $k(\Gamma_i)$ of the component $\Gamma_i$
of $\X_k$. As noted in \ref{rem.multiplicity}, the scheme 
$\Gamma_i^{\reg}$ is defined over $k_i$, and we 
have $\delta(\Gamma_i^{\reg}/k) = [k_i:k]\delta(\Gamma_i^{\reg}/k_i)$.
The geometric multiplicity $e(\Gamma_i/k_i)$ of $\Gamma_i/k_i$ is also 
the geometric multiplicity of $\Gamma^{\reg}_i/k_i$. It follows then 
from above 
that  $e(\Gamma_i/k_i) $ divides $ \delta(\Gamma_i^{\reg}/k_i)$.
Hence, Theorem \ref{dx-ds} (a) implies that 
$$\gcd_i\{r_i[k_i:k]e(\Gamma_i/k_i)\} \text{ divides } \delta(X/K),$$
answering a question in \cite{Bosch-Liu}, 1.6. 
\end{remark} 

\begin{remark} (1) In general, if $\cO_K$ is not Henselian, 
$\gcd_i\{r_i\delta(\Gamma_i^{\reg}/k)\}$ is not 
equal to $\delta(X/K)$. This can be seen easily
when $\X/\cO_K$ is of relative dimension $0$. 
We can also consider a smooth projective conic $X$ over $\mathbb Q$
without rational point, and  with a regular proper model $\X$ over $\mathbb Z$. 
If $p$ is a prime of good reduction of $X$, then $\gcd(\X_{{\mathbb F}_p})=1$ 
because every smooth conic over ${\mathbb F}_p$
has an ${\mathbb F}_p$-rational point, 
while $\delta(X/\mathbb Q)=2$. 

(2) In general we cannot replace   
$\delta(\Gamma_i^{\reg}/k)$ by $\delta(\Gamma_i/k)$ in \ref{dx-ds}(b). For example, 
let $\cO_K=\mathbb R[[t]]$ and let $\X=\Proj(\cO_K[x,y,z]/(x^2+y^2+tz^2))$. 
Then $\X$ is regular, flat and projective over 
$\cO_K$. The special fiber $\Gamma:=\X_k$ is integral, with a singular rational
point, and its regular locus is isomorphic to $\mathbb A^1_{\mathbb C}$.
In this example, $\delta(X/K)=2$, but 
$r(\Gamma)\delta(\Gamma/\mathbb R)=1$.
\end{remark}
\begin{example} 
Let $A/K$ be a central simple $K$-algebra of dimension $n^2$.
The square root of the degree over $K$ of the skew-field $D$ such that 
$A$ is isomorphic to $M_r(D)$ for some $r\geq 1$ is called the index 
${\rm ind}(A)$ of $A$. Associated with $A$ is a Severi-Brauer variety 
$X/K$, a twisted form of ${\mathbb P}^{n-1}/K$, 
with $\delta(X/K) = {\rm ind}(A)$. 

Suppose that $K$ is a complete discrete valuation field with perfect residue field $k$. 
Let $\Lambda/{\mathcal O}_K$ be a hereditary order of $A$.
Associated with $\Lambda$ is a model ${\mathcal X}/{\mathcal O}_K$ of $X/K$ called an Artin model
(\cite{Fro}, 2.1). Artin (\cite{Art}, 1.4) shows that when $\Lambda$ is a maximal order, 
the model $\X$ is regular. The special fiber $\X_k$ is described in some cases in 
2.4 and 2.5 of \cite{Fro}. In particular, the special fiber in the model 
described in \cite{Fro} 2.5 contains only irreducible components $\Gamma$ of 
multiplicity $1$ such that $\delta(\Gamma/k) =1$, 
and $\delta(\Gamma^{\reg}/k)= \delta(X/K)$.
\end{example}

\begin{emp} Let $W$ be a non-empty scheme of finite type over a field $F$. 
Let $\mathcal D$ be the set of all degrees of closed points of $W$.
Denote by $\nu(W/F)$ the smallest integer in $\mathcal D$. 
Clearly, $\delta(W/F) $ divides $\nu(W/F)$. 

Let $f:\X\to S$ be as in \ref{defX/S}, and assume $\X$ regular and $f$ proper. 
Then Proposition \ref{lift} (2) shows that
$$ \nu(X/K)\ge 
\min_{x_0 \text{ closed in } \X_k}\left\{e(\cO_{\X_k, x_0}) \deg_k(x_0)\right\},$$
with equality if in addition $\cO_K$ is Henselian and $k$ is infinite (\ref{lift} (3)). 

Lang and Tate asked in 
\cite{LT}, page 670, whether $\nu(W/F)= \delta(W/F)$ when $W$ is a 
homogeneous space 
for an abelian variety $A/F$.
Recall that non-empty scheme $W/F$ is a homogeneous space for $A/F$ if 
$W/F$ is endowed with  a transitive action of $A/F$, in the sense that
the natural morphism $A \times_F W \to W \times_F W$, which sends
$(a,w) $ to $(a\cdot w, w)$, is surjective as a morphism of fppf-sheaves.

Let us say that a field $k$ is \nice \
if every homogeneous space $X/k$ for any abelian variety $B/k$ has a $k$-rational point (\cite{Cla2}, 1.2).
A finite field is \nice \ (\cite{Lang}, Theorem 2), and it follows from the definitions that a pseudo-algebraically closed field $k$ is \nice. 
\end{emp}

\begin{proposition}
Let $\cO_K$ be Henselian with 
a \nice \ perfect residue field. Let $A$ be  an abelian variety over 
$K$ having good reduction. Let  $X/K$ be a 
homogeneous space for 
$A$. Then $\delta(X/K)=\nu(X/K)$. 
\end{proposition}

\proof As $A$ is an abelian variety, $X$ is a principal homogeneous space for a quotient 
$B$ of $A$. Let $\mathcal B/\cok$ be the N\'eron model of $B$ over $\cok$. 
Then $\mathcal B$ is an abelian scheme (\cite{S-T}, Corollary 2 to 
Theorem 1). We know (\cite{LLR}, Proposition 8.1) that $X$ has 
a regular proper model $\X/\cok$ such that 
$\X_k=rV$ with $V$ proper, smooth over $k$ (because $k$ is perfect), and $V/k$ is a 
homogeneous space for $\mathcal B_k$. Hence, $V(k)\ne\emptyset$,  so that $\delta(V/k)=\nu(V/k)=1$. 
Using Proposition \ref{lift}, we find that $\nu(X/K)=r\nu(V/k)=r$. 
Theorem \ref{dx-ds} shows that $\delta(X/K)=r\delta(V/k)=r$. 
\qed 
 
\end{section}
\begin{section}{The separable index}

Let $k$ be any field, and let $X/k$ be a scheme. 
The set of separable closed points of $X$ can be empty, even when $X/k$ is regular. This is the case for instance
for $X=\Spec L$, 
where $L/k$ is a non-trivial purely inseparable extension. On the other hand, when $X/k$ is smooth, 
the set of separable closed points of $X$ is dense in $X$ (\cite{BLR}, 2.2/13).

When $X/k$ is generically smooth and non-empty, define the {\it separable index} $\delta_{sep}(X/k)$ of $X/k$ to be 
the greatest common divisor of the degrees of the separable closed points of $X$.
Clearly, $\delta(X/k) $ divides $ \delta_{sep}(X/k)$, and the question of whether $\delta(X/k) $ and $ \delta_{sep}(X/k)$
are always equal was raised by
Lang and Tate in 
\cite{LT}, page 670. In this section, we answer this question positively. The case where $X/k$ is a smooth projective curve
was treated already in \cite{Har}, Theorem 3.

Let us note first that our work in the previous section lets us answer 
the question positively under the following additional hypotheses.

\begin{corollary} \label{dx-ds.cor}  
Let $K$ be the field of fractions of a Henselian discrete valuation ring $\cO_K$. 
Let $f:\X \to S$ be as in {\rm \ref{defX/S}}, with $f$ proper and flat, and $\X$ regular.

Let $X/K$ denote the  generic fiber of $\X/S$,
assumed to be generically smooth. 
Then 
$$\delta(X/K) = \delta_{sep}(X/K)=\gcd_i\{r_i \delta(\Gamma^{\reg}_{i}/k)\}.$$
\end{corollary}
\proof Follows immediately from Theorem \ref{dx-ds} (b) and Proposition \ref{lift} (3). \qed

\medskip
The proof of the next theorem is independent of the results in the previous sections of the paper.
\begin{thm} \label{sep.index} Let $X$ be a regular and generically smooth non-empty scheme of finite type over a
field $k$. Then $\delta(X/k)=\delta_{sep}(X/k)$. 
\end{thm}

\proof Let $X^{sm}$ denote the dense open subset of $X$ where $X/k$ is smooth.
It follows immediately from Proposition \ref{pro.deltareg} and the fact that $X$ is regular
that $\delta(X/k) = \delta(X^{sm}/k)$. Thus, it suffices to prove the theorem when $X/k$ is smooth.

Assume that $X/k$ is smooth. It suffices to prove that for any closed point $x_0\in X$,
there exists a separable $0$-cycle on $X$ of  degree $\deg_k(x_0)$. 
We proceed by induction on $d:=\dim_{k(x_0)} \Omega^1_{k(x_0)/k}$. 
If $d=0$, then $x_0$ is already separable. If $d=1$, 
then by Lemma \ref{embed}, $x_0$ belongs to a closed curve  $C$ in
$X$ which is smooth at $x_0$, and we can conclude using Lemma \ref{curve}.
 
Now fix $d\ge 2$, and suppose that the statement holds for any 
closed point $y$ in any smooth variety over $k$ such that 
$\dim_{k(y)}\Omega^1_{k(y)/k}\le d-1$. Replacing $X$ by the smooth locus of the closed subvariety $Y$ 
whose existence is proved in Lemma \ref{embed}, 
we can suppose that $\dim_{x_0} X=d$. Since $x_0$ is smooth in $X$, 
we can replace $X$ by an affine open subset if necessary, and assume 
that there exists an \'etale morphism $f : X\to \mathbb A^d_k$. Let 
$\pi : {\mathbb A}^{d}_k\to {\mathbb A}^{d-1}_k$ be the projection 
to the first $d-1$ coordinates. Let $z_0:=\pi(f(x_0))$. Then 
$C:=X\times_{{\mathbb A}^{d-1}_k} \Spec k(z_0)$ is a closed subscheme of $X$
containing $x_0$. It is a smooth curve over $k(z_0)$ 
because it is \'etale over the affine line 
$\pi^{-1}(z_0)={\mathbb A}^1_{k(z_0)}$. Therefore, we can use Lemma \ref{curve}, and obtain that 
there exists a $0$-cycle $D$ on $C$, separable over $k(z_0)$, and of 
degree $\deg_{k(z_0)}(x_0)$. For any point $y$ in the support of $D$, we have 
$$\dim_{k(y)}\Omega^1_{k(y)/k}=\dim_{k(z_0)}\Omega^1_{k(z_0)/k}
\le \dim {\mathbb A}^{d-1}_k=d-1.$$ 
Thus, we can apply the induction hypothesis to each such point $y$. It follows that there exists a separable $0$-cycle on $X/k$
of degree $\deg_k(x_0)$. 
\qed

\begin{lemma} \label{embed} 
Let $X$ be a smooth algebraic variety over a field $k$, and let  
$x_0\in X$. 
Then $x_0$ is contained in a closed subvariety $Y$ of $X$ of pure dimension $\dim_{k(x_0)} \Omega^{1}_{k(x_0)/k}$ 
which is smooth at $x_0$.
\end{lemma}

\proof Let $\m_0$ be the maximal ideal of $\cO_{X,x_0}$.
When $f \in \m_0$, denote by $\bar{f}$ its class in $\m_0/\m_0^2$. The
closed immersion $\Spec k(x_0)\to \Spec \cO_{X,x_0}$ induces the canonical second fundamental exact sequence 
$$ \m_0/\m_0^2 \overset{\delta}{\longrightarrow} \Omega^1_{X,x_0}\otimes k(x_0) \overset{\mu}{\longrightarrow }
\Omega^1_{k(x_0)/k}\longrightarrow  0 $$
(\cite{Mat}, Theorem 58, page 187).
Let $f_1,\ldots,f_m \in \m_0$ be such that 
$\delta(\bar{f}_1),\ldots,\delta(\bar{f}_m)$ form 
a basis of 
$\Ker(\mu)$. 
Choose an open neighborhood $U$ of $x_0$ and $g_i \in \cO_{X}(U)$ such that for all $i=1,\dots, m$,
 $f_i$ is the stalk of $g_i$ at $x_0$.
Let $Y$ be a closed subscheme of $X$, equal to $V(g_1,\ldots,g_m)$
in an open neighborhood of $x_0$.  As $\{f_1, \ldots, f_m \}$ is part of 
a regular system of parameters of $\cO_{X,x_0}$, we find that
$$\dim_{x_0}Y=\dim_{x_0}X - m = \dim_{k(x_0)} \Omega^{1}_{k(x_0)/k}.$$
Since 
$(\Omega^1_{Y/k})_{x_0}\otimes k(x_0)\simeq \Omega^1_{k(x_0)/k}$,
we find that $(\Omega^1_{Y/k})_{x_0}$ is generated by  $\dim_{x_0}Y$ elements,
and this implies that $Y$ is
smooth at $x_0$ (see, e.g., \cite{BLR}, 2.2/15). We may thus replace $Y$ by its irreducible component which contains $x_0$.
\qed 

\medskip
  
\begin{lemma}  \label{curve} 
Let $C$ be a smooth connected curve over $k$. Let 
$\overline{C}$ be a regular scheme, separated and of finite type over $k$, 
such that there exists an open immersion $C \to \overline{C}$. 
Let $z \in C$ be a closed point. Then  there exists a separable 
$0$-cycle on $C$ which is rationally equivalent to $[z]$ in 
$\overline{C}$. 
\end{lemma}
\proof
The connected component of $\overline{C}$ containing $C$ embeds as an open subscheme
of a regular compactification $C'$ of $C$. 
Lemma \ref{curve} follows from a stronger statement
(\cite{FSV}, Chapter 3, Lemma 3.16), 
 in the proof of which   the smooth compactification (which may not exist
 in general) should be replaced by a regular compactification. 
 We thank O. Wittenberg for the reference to \cite{FSV}. \qed

\end{section}

\end{document}